\documentclass[10pt,a4paper,oneside]{amsart}

\usepackage{amsfonts,amssymb}
\usepackage{amsmath,amsthm}
\usepackage{graphicx, color}
\usepackage{enumerate,array,enumitem}
\usepackage{url}
\usepackage{tikz-cd}
\usepackage[hidelinks]{hyperref}

\usepackage[all]{xy}

\usetikzlibrary{arrows.meta,calc}

\newtheorem{thmA}{Theorem}

\newtheorem{theorem}{Theorem}[section]
\newtheorem{prop}[theorem]{Proposition}
\newtheorem{lemma}[theorem]{Lemma}
\newtheorem{cor}[theorem] {Corollary}

\newtheorem*{claim*}{Claim}

\theoremstyle{remark}
\newtheorem{remark}[theorem]{Remark}

\theoremstyle{definition}
\newtheorem{definition}[theorem]{Definition}
\newtheorem{notation}[theorem]{Notation}

\newcommand{\rar}{\rightarrow}
\newcommand{\id}{\textnormal{id}}

\newcommand{\Z}{\mathbb{Z}}
\newcommand{\C}{\mathbb{C}}
\newcommand{\R}{\mathbb{R}}

\newcommand{\K}{\mathbb{K}}
\newcommand{\metricdefect}[1][h]{\nu_{#1}}
\newcommand{\opnorm}[3]{\|#1\|_{L(#2,#3)}}

\newcommand{\im}{\operatorname{im}}

\makeatletter
\newcommand{\xRightarrow}[2][]{\ext@arrow 0359\Rightarrowfill@{#1}{#2}}
\makeatother

\title[FEEC for the Hodge--Dirac and Hodge--Laplace problems]{Convergence of finite element exterior calculus for the Hodge--Dirac and Hodge--Laplace problems on closed Riemannian manifolds}
\author{Thomas A. Wasserman}
\address{Department of Mathematics and Systems Analysis, Aalto University}
\email{thomas.wasserman@aalto.fi}
\date{\today}
\begin{document}
	\begin{abstract}
	We prove that the finite element exterior calculus (FEEC) approximation to the Hodge--Dirac and
	Hodge--Laplace problems on any closed oriented Riemannian manifold converges.  Its approximate solutions, together with their exterior derivatives, converge to the continuum solutions for every source, and at the rate of the mesh size when the source and its codifferential are square integrable.

	We use Licht's estimate for the FEEC approximation of forms	with one square-integrable derivative by finite element forms.	We combine this with a new abstract error bound for approximations of Hilbert complexes. This bound uses a version of Di Pietro	and Droniou's third Strang lemma. To apply these estimates, we show that every closed manifold admits a family of arbitrarily fine triangulations whose shapes stay uniformly under control.
\end{abstract}
	\maketitle

\setcounter{tocdepth}{2}
\tableofcontents

\section{Introduction}
We show that on any closed oriented Riemannian manifold, finite element exterior calculus converges for the Hodge--Dirac and Hodge--Laplace problems, the weak, mixed formulations of the differential equations
$$
(d+d^*)u=F \quad\mbox{and}\quad \Delta u =F.
$$
We need no regularity of the source $F$, and prove convergence at the rate of the mesh size when the source satisfies a mild regularity condition.%

Finite element exterior calculus approximates differential forms and their
calculus through piecewise polynomial forms on a triangulation. For Euclidean domains, it is a well-established formalism \cite{Arnold2006,arnold2010finite}. On manifolds, Holst and Stern \cite{Holst2012} developed its error analysis abstractly by treating the difference between the discrete and the exact geometry as a variational crime. Licht \cite{Licht2023smoothed}
constructed the smoothed projections their analysis needs, proved the
approximation estimates, and bounded the error made by computing with an
approximate metric. These analyses treat only the Hodge--Laplace problem, and assume families of triangulations whose shapes are uniformly controlled, which on a general closed manifold have to be constructed. We also cover the Hodge--Dirac problem, and check that repeated Freudenthal subdivisions provide the triangulations.

The known error bounds are suited for finite element exterior calculus. The bounds we give are sharper and needed for the discrete exterior calculus application we have in mind. Our sharper bounds come from an adaptation of Di Pietro and Droniou's third Strang lemma \cite{DiPietroDroniou2018}, and do not pass through Holst and Stern's estimate.

\subsection{Main results}\label{sec: main results}
Throughout, $(M,g)$ is a closed
oriented Riemannian manifold of dimension $n$. On such a manifold we have the exterior derivative $d$, the Hodge star $\star$, the codifferential $d^*=\pm\star d\star$, the Hodge--Dirac operator $d+d^*$, and its square, the Hodge--Laplacian. Forms in the kernel of the Hodge--Laplacian are called \emph{harmonic}.

\subsubsection{Abstract error bounds}
We use the abstract framework of Hilbert complexes, chain
complexes of Hilbert spaces with closed, densely defined differentials
\cite{arnold2010finite,Holst2012}. Abstract versions of both the
Hodge--Dirac and Hodge--Laplace problems can be formulated for these, allowing a functional on the domain of $d$ as source. An \emph{approximation} of a ``continuous'' Hilbert complex consists of a ``discrete'' finite-dimensional Hilbert complex $V_h^\bullet$, a map $i_h$ into the continuous complex and a bounded projection $\pi_h$ back, both commuting with the differentials. We then approximate the continuum solutions by solutions to the problems on the discrete complexes for some choice of source. Both problems constrain the solution to be orthogonal to the harmonic forms; their discrete versions use the discrete harmonic forms. The error of the approximation is controlled by three quantities: the \emph{best approximation error}, the distance from the solution to $i_h(V_h^\bullet)$; the \emph{harmonic gap} $\mu_h$, which measures how well the approximation reproduces harmonic forms; and, when $i_h$ does not preserve inner products, the \emph{metric defect} $\metricdefect$, the norm of the bilinear form $\langle(\mathrm{id}-i_h^*i_h)\cdot,\cdot\rangle_h$ (where $i_h^*$ is the adjoint of $i_h$ and $\langle\cdot,\cdot\rangle_h$ the inner product on $V^\bullet_h$) on the domain of the differential. Our abstract bound, loosely stated, is then:
\begin{thmA}[{Theorems~\ref{thm: Hodge Dirac abstract approximation theorem} and~\ref{thm:Hodge Laplace bound}}]\label{thm: intro abstract}
	Let $V_h^\bullet$ be an approximation of a closed Hilbert complex.  The
	error in the discrete Hodge--Dirac or Hodge--Laplace problem is bounded by a multiple of the sum of the best approximation error of the solution,
	the error in the source of the discrete problem, the metric defect $\metricdefect$ times
	the norm of the source, and the harmonic gap $\mu_h$ times a best approximation
	error.
\end{thmA}

\subsubsection{Uniformly regular triangulations}
Finite element exterior calculus on $(M,g)$ requires a smooth
triangulation $K$ of $M$, whose closed simplices are smoothly embedded. A family $\mathcal K$ of such triangulations is \emph{uniformly regular} if the embeddings of the simplices
stay under control, uniformly across the family, and the mesh size $\ell_K$,
the largest diameter of a simplex of $K$, takes arbitrarily small values on
$\mathcal K$. We call the constants that control the embeddings of the
simplices of $K$ the \emph{regularity} of $K$. We show that a standard
construction gives such families on every closed manifold.
\begin{thmA}[{Theorem~\ref{thm: freudenthal subdivisions uniformly regular}}]\label{thm: intro triangulations}
	Every closed Riemannian manifold carries a uniformly regular family of
	smooth triangulations: the repeated Freudenthal subdivisions of any
	smooth triangulation of $M$.
\end{thmA}

\subsubsection{Convergence of finite element exterior calculus}\label{sec: results convergence}
We apply Theorem~\ref{thm: intro abstract} to the finite element exterior calculus approximations to square-integrable forms, on a uniformly regular family. On a smooth triangulation $K$ the finite element complex $V_K^\bullet$ is built from the polynomial spaces of Arnold, Falk and Winther \cite{Arnold2006}. Its inclusion $i_K$ into the square-integrable forms commutes with the differentials. Licht's smoothed projection \cite{Licht2023smoothed} serves as a bounded commuting projection $\pi_K$ back. The complex $V_K^\bullet$ carries the inner products of the piecewise flat metric $g_K$. This metric makes each simplex of $K$ a flat simplex with edge lengths those of $(M,g)$. Since $g_K$ differs from $g$, the inclusion $i_K$ is not an isometry.

Along a uniformly regular family, the difference between $g_K$ and $g$ is at
most of order $\ell_K$, and so is the metric defect.
The best approximation errors in Theorem~\ref{thm: intro abstract} are
bounded by an approximation estimate of Licht \cite{Licht2023smoothed}: every form with one square-integrable derivative is
approximated to order $\ell_K$ by the finite element space in the norm of
square-integrable forms. Somewhat informally stated, our main result is:
\begin{thmA}[{Theorem~\ref{thm: feec convergence} and Corollary~\ref{cor: feec L2 convergence}}]\label{thm: intro convergence}
	Let $\mathcal K$ be a uniformly regular family of smooth triangulations of
	$M$, each carrying a finite element complex with the inner products of its piecewise flat metric.  Solve the Hodge--Dirac or
	the Hodge--Laplace problem with a source $F$ and
	its discrete counterpart with source $F\circ i_K$.  Then the image under $i_K$ of the
	discrete solution converges to the solution, together with its exterior derivative, as
	$\ell_K\to0$.  If $F$ is the inner product with a square-integrable form $\gamma$ whose codifferential $d^*\gamma$ is square integrable too, then the convergence has rate $\ell_K$.
\end{thmA}

\subsection{Application to discrete exterior calculus}\label{sec: dec application}
In a companion paper \cite{WassermanDEC} we use our abstract bounds to establish convergence of  discrete exterior calculus on surfaces. Discrete exterior calculus \cite{Hirani2003} is another approach to approximating differential forms and their
calculus. It represents forms by simplicial cochains. Like finite element
exterior calculus with the piecewise flat metric, it uses the Riemannian metric
only through the lengths of the edges of the triangulation. On the lowest-order
finite element spaces, the Whitney forms, the two methods differ only in their
inner products.

This difference is where our results are needed. The error bounds of Holst and Stern
\cite{Holst2012} need the discrete inner products to be close to the exact ones
on all discrete forms, measured in the square-integrable norm.
Theorem~\ref{thm: intro abstract} needs this only when discrete forms are
measured in the \emph{graph norm} (the sum of the square-integrable norms of the form and its exterior derivative), through the metric defect. Finite element exterior
calculus satisfies both conditions, but discrete exterior calculus satisfies
only the second.

Discrete exterior calculus needs well-centred triangulations, in which each
simplex contains its circumcentre. This is one reason we use Freudenthal's subdivision for
Theorem~\ref{thm: intro triangulations}. On a surface it copies the shape of
each initial triangle, so well-centredness is decided by the shapes of the initial triangles, measured in the metric pulled back to the polyhedron.

\subsection{Background and related work}\label{sec: related work}

\subsubsection{Hilbert complexes and error bounds}\label{sec: related abstract}
Hilbert complexes were introduced by Br\"uning and Lesch to treat elliptic
complexes abstractly \cite{BruningLesch1992}, and Arnold, Falk and Winther made
them the framework of finite element exterior calculus
\cite{Arnold2006,arnold2010finite,Arnold2018}. They gave the mixed formulation
of the Hodge--Laplace problem used here, and Leopardi and Stern that of the
Hodge--Dirac problem \cite{LeopardiStern2016}, each with an analysis of the
error of approximations by subcomplexes. Holst and Stern extended the analysis
of the Hodge--Laplace problem to approximations by a complex with inner
products of its own \cite{Holst2012}, motivated by the FEEC approximation to forms on Riemannian manifolds. They compare the discrete problem with an intermediate problem, posed on the image
of $i_h$ with the continuous inner products. We bound the error directly, with a version of
the third Strang lemma of Di Pietro and Droniou \cite{DiPietroDroniou2018}.

Theorem~\ref{thm: intro abstract} recovers the bound of Holst and Stern, with
a sharper metric term. Holst and Stern measure the failure of $i_h$ to
preserve inner products by the norm of the bilinear form
$\langle(\mathrm{id}-i_h^*i_h)\cdot,\cdot\rangle_h$ with respect to the
square-integrable norm; the metric defect is its norm with respect to the graph
norm, so it is never larger. Our hypotheses are also weaker: once the metric defect
is small, we need no bound on $i_h$ for the square-integrable norms. For the
Hodge--Dirac problem, which Holst and Stern do not treat,
Theorem~\ref{thm: intro abstract} extends the estimate of Leopardi and Stern to
approximations that do not respect the inner product.

\subsubsection{Triangulations}\label{sec: related triangulations}
Our estimates hold along uniformly regular families of triangulations.
Licht \cite{Licht2023smoothed} outlines how to obtain such families on a manifold: refine a triangulation by a rule that
keeps the shapes of its simplices under control, and carry the refinements to
the manifold by a fixed piecewise smooth homeomorphism. Theorem~\ref{thm: intro triangulations} carries this out for Freudenthal's subdivision, starting from a smooth triangulation, which every closed manifold has by the theorem of Cairns and Whitehead \cite{Cairns1935,Whitehead1940}.

Freudenthal's subdivision \cite{Freudenthal1942} is known in the finite element
literature as red refinement \cite{Bey2000}, and is a special case of the
edgewise subdivision of Edelsbrunner and Grayson
\cite{EdelsbrunnerGrayson2000}. We use it because it is a standard refinement
in finite element computations, and because on surfaces it copies the shape of
each triangle, which the application to discrete exterior calculus needs. Whitney's standard subdivision
\cite{Whitney1957}, which for example Dodziuk uses for his convergence theorems
\cite{Dodziuk1976}, does not.

Triangulations adapted to the metric
\cite{BoissonnatDyerGhosh2018,DyerVegterWintraecken2015}, and metrics within
relative error of order $\ell_K^2$ of $g$
\cite{DyerVegterWintraecken2015,BoissonnatDyerGhoshWintraecken2023,MeyerToubiana2023},
whose construction goes back to Cheeger, M\"uller and Schrader
\cite{Cheeger1984}, come, as far as we know, without the bounds on higher
derivatives of the parametrisations that uniform regularity asks for.

\subsubsection{Approximation of forms}\label{sec: related approximation}
Approximating forms requires constructing an interpolant. The canonical
interpolant onto the finite element spaces is defined for continuous forms, and
on $k$-forms with $s$ square-integrable derivatives it is bounded only for
$s>(n-k)/2$ \cite{Arnold2006}. For rougher forms the standard tool is a
quasi-interpolant, which replaces point values by local averages
\cite{Clement1975,ScottZhang1990}. Gawlik, Holst and Licht extend these to
finite element differential forms \cite{GawlikHolstLicht2021}, see also
\cite{ErnGuermond2017}, and Licht \cite{Licht2023smoothed} carries their estimate to closed manifolds. Before Licht's estimate, results on closed manifolds of arbitrary dimension
either required more regularity, as for the Whitney interpolant of smooth forms
\cite{Dodziuk1976}, or gave no rate \cite{Christiansen2007}.

Stability of the discrete problems, uniformly in the mesh, rests on bounded
projections onto the finite element spaces that commute with the exterior
derivative \cite{arnold2010finite}; our error bounds also compare the discrete
solution with the projection of the solution. Such smoothed projections go back
to Sch\"oberl and to Christiansen
\cite{Schoberl2001,Schoberl2008,Christiansen2007} and were refined in
\cite{Arnold2006,ChristiansenWinther2008,FalkWinther2014}; we use the ones
Licht constructed on closed manifolds \cite{Licht2023smoothed}.

\subsubsection{Convergence of finite element exterior calculus}\label{sec: related convergence}
Convergence on manifolds has been proved in three settings. For the Whitney
forms with the inner product pulled back from the square-integrable forms,
Dodziuk proved that the combinatorial Hodge decomposition converges to the
Hodge decomposition \cite{Dodziuk1976}, and Dodziuk and Patodi that the Hodge
components converge at the rate of the mesh size up to a logarithmic factor
\cite{DodziukPatodi1976}, as announced in \cite{Patodi1974}.

On an approximating manifold, Holst and Stern bound the error by the best
approximation error and a term measuring how much the approximating manifold
distorts the metric. For hypersurfaces approximated by piecewise polynomial
surfaces, as in the surface finite element methods of Dziuk and Demlow
\cite{Dziuk1988,Demlow2009,DziukElliott2013}, this gives convergence at a rate,
assuming elliptic regularity \cite{Holst2012}.

On the manifold itself, Licht's smoothed projections and the Galerkin theory of
Arnold, Falk and Winther give stability and quasi-optimality for the exact
metric. Licht bounds the additional error of an approximate metric for the
Hodge--Laplace problem, of the order of the $(l+1)$-st power of the mesh size
for a metric interpolated by polynomials of degree $l$, and works out rates for
a source with Sobolev regularity, assuming full elliptic regularity and rates
for three auxiliary approximation quantities \cite{Licht2023smoothed}. Gawlik
and McKee bound the metric term of Holst and Stern's estimate when the manifold carries a Regge
metric, such as the piecewise flat metric, in place of $g$
\cite{GawlikMcKee2024}, and we use their bound. Together with Licht's smoothed
projections and estimate and Theorem~\ref{thm: intro triangulations}, this therefore also gives convergence for the Hodge--Laplace
problem with the piecewise flat metric; Theorem~\ref{thm: intro abstract} is
needed for the Hodge--Dirac problem and for discrete exterior calculus
(Section~\ref{sec: dec application}).

\subsection{Structure of the paper}
Section~\ref{sec: abstract problems} proves Theorem~\ref{thm: intro abstract}
for an abstract Hilbert complex and turns it into rates.
Section~\ref{sec: feec on manifolds} treats finite element exterior calculus on
$M$: Sections~\ref{sec: triangulations} and~\ref{sec: existence of triangulations}
prove Theorem~\ref{thm: intro triangulations}, Section~\ref{sec:FEECforms}
bounds the metric defect, proves the regularity of the solutions and derives
the approximation estimate we need from Licht's, and
Section~\ref{sec: convergence on manifolds} proves
Theorem~\ref{thm: intro convergence}.

\section{The abstract Hodge--Dirac and Hodge--Laplace problems}\label{sec: abstract problems}
Section~\ref{sec: variational problems and strang} proves a version of the third Strang lemma for variational problems whose solutions are constrained to be orthogonal to a subspace (Theorem~\ref{thm: constrained strang}). Section~\ref{sec: Hilbert complexes} recalls Hilbert complexes and poses the Hodge--Dirac and Hodge--Laplace problems on them. Section~\ref{sec: approximating hilbert complexes} defines approximations of a Hilbert complex and the error terms associated to them. Section~\ref{sec: abstract error bounds} proves the error bounds on approximations to the Hodge--Dirac (Theorem~\ref{thm: Hodge Dirac abstract approximation theorem}) and Hodge--Laplace (Theorem~\ref{thm:Hodge Laplace bound}) problems by applying Theorem~\ref{thm: constrained strang}, and turns them into rates (Section~\ref{sec: abstract rates}).

\addtocontents{toc}{\protect\setcounter{tocdepth}{0}}
\subsection*{Conventions}
Throughout, $\mathbb K$ denotes $\R$ or $\C$, fixed once and for all, and every
Hilbert space is a Hilbert space over $\mathbb K$. We denote the inner product
on a Hilbert space $V$ by $\langle\cdot,\cdot\rangle_V$ and its norm by $\|\cdot\|_V$.
For normed spaces $X$ and $Y$ and a bounded linear map $A\colon X\to Y$, we
write $\opnorm{A}{X}{Y}$ for the operator norm of $A$.
Several results we quote are stated for real Hilbert spaces; routine arguments
adapt them to complex Hilbert spaces. By a \emph{bilinear form} on
$V\times W$ we mean a map $V\times W\to\mathbb K$ that is linear in the first
argument and conjugate-linear in the second. We denote by $V^*$ the space of bounded conjugate-linear functionals on $V$, so
that for $\gamma\in V$ the map $v\mapsto\langle\gamma,v\rangle_V$ lies in $V^*$.

\addtocontents{toc}{\protect\setcounter{tocdepth}{2}}
\subsection{Variational problems and Strang-type lemmas}\label{sec: variational problems and strang}
This section recalls the basics of variational problems and proves the version
of the third Strang lemma that we need (Theorem~\ref{thm: constrained strang}).

\subsubsection{Variational problems}
The Hodge--Dirac and Hodge--Laplace problems will be posed as variational
problems.

\begin{definition}\label{def: variational problem}
    A \emph{variational problem} is a Hilbert space $V$ equipped with a bounded
    bilinear form $B\colon V\times V\to\K$ and an element $F\in V^*$. A \emph{solution} to the variational problem is $u\in V$ such that
\begin{equation}\label{eq:generic variational problem}
    B(u,v)=F(v) \quad \forall v\in V.
\end{equation}
\end{definition}

The first question when confronted with a variational problem is whether it is \emph{well-posed}, i.e. whether it has a unique solution. The following condition is sufficient for well-posedness.

\begin{definition}\label{def: infsup}
	A bilinear form $B\colon V\times V \to \K$ is called \emph{inf-sup stable} if
	there is a constant $\sigma>0$ such that, for all $u\in V$,
	\begin{equation}\label{eq: infsup}
		\sigma \|u\|_V\ \le\ \sup_{v\in V\backslash\{0\}}
		\frac{|B(u,v)|}{\|v\|_V}
		\qquad\text{and}\qquad
		\sigma \|u\|_V\ \le\ \sup_{v\in V\backslash\{0\}}
		\frac{|B(v,u)|}{\|v\|_V}.
	\end{equation}
\end{definition}

Inf-sup stability gives well-posedness, together with a bound on the solution in
terms of the data. This bound is how stability enters the error estimates
below.

\begin{theorem}[{\cite[Theorem~2.1]{Babuska1971}}]\label{thm: babuska}
	Any variational problem $(V,B,F)$ with an inf-sup stable bilinear form $B$ is well-posed. Its solution satisfies $\|u\|_V\le\sigma^{-1}\|F\|_{V^*}$, where $\sigma$ is the inf-sup constant of $B$.
\end{theorem}

Much of the finite element literature reserves the names \emph{inf-sup condition} and \emph{inf-sup stability} for the first inequality of \eqref{eq: infsup} alone, imposing the second separately where it is needed; see for instance \cite[Definition~4]{DiPietroDroniou2018}.  On its own the first inequality says that $u\mapsto B(u,\cdot)$ is injective with closed range.

\subsubsection{A Strang-type lemma}
To approximate a variational problem, a
common approach is to pick a finite-dimensional subspace $V_h$ of $V$ and to
solve the \emph{Galerkin problem}, obtained by restricting $B$ and $F$ to $V_h$.
\emph{Criminal methods} also replace $B$ and $F$ by discrete approximations
$B_h$ and $F_h$. We go one step further and allow $V_h$ not to be a subspace of
$V$. The map back from $V_h$, the reconstructor of the next definition, may then take values in a Hilbert space $E$ containing
$V$. This ensures that the constrained problems of Section~\ref{sec: constrained problems}
fit the same framework.
\begin{definition}\label{def: criminal}
    Let $(V,B,F)$ be a variational problem. A \emph{criminal approximation} of $(V,B,F)$ is a variational problem $(V_h,B_h,F_h)$ on a finite-dimensional Hilbert space $V_h$ together with bounded, not necessarily isometric, linear maps
    $\pi_h\colon V\to V_h$, the \emph{interpolator}, and $i_h\colon V_h\to E$,
    the \emph{reconstructor}, where $E$ is a Hilbert space that contains $V$
    isometrically. We will refer to $(V_h,B_h,F_h)$ and other associated objects as \emph{discrete}, and to $(V,B,F)$ as \emph{continuous}.
\end{definition}

For the Galerkin problem, the discrete solution is the best
approximation from $V_h$ up to a constant.  Strang's first and second
lemmas \cite{strang1972variational}, \cite[Theorems~4.1.1 and~4.2.2]{Ciarlet1978},
extend this to criminal methods, but the first needs $V_h\subseteq V$, and the
second needs $V$ and $V_h$ to lie in a common space to which $B_h$ extends. We use the third Strang lemma of Di Pietro and
Droniou, which needs neither.

\begin{theorem}[{\cite[Section~2.4.1]{DiPietroDroniou2018}}]\label{thm: third strang}
    Let $(V_h,B_h,F_h)$ be a criminal approximation of $(V,B,F)$ with interpolator $\pi_h$ and reconstructor $i_h\colon V_h \to E$ where $E\supseteq V$, assume that $B_h$ is inf-sup stable, and that $u$ and $u_h$ are solutions of $(V,B,F)$ and $(V_h,B_h,F_h)$ respectively. Then we have
    \begin{equation}\label{eq: strang}
    \|u-i_h u_h\|_E\leq \|i_h \pi_h u-u\|_E + \opnorm{i_h}{V_h}{E}\,\sigma_h^{-1}\|F_h(\cdot)-B_h(\pi_h u, \cdot)\|_{V_h^*},
    \end{equation}
    where $\sigma_h$ is the inf-sup constant of $B_h$.
\end{theorem}

We call the first term in \eqref{eq: strang} the \emph{approximation error}
and the second the \emph{consistency error}.  The approximation error depends
only on $u$, the interpolator and the reconstructor, not on $B_h$ and
$F_h$.  The variational crime enters through the consistency error.

\subsubsection{Constrained problems}\label{sec: constrained problems}
The solutions of the Hodge--Dirac and Hodge--Laplace problems are required to
be orthogonal to the harmonic forms. These problems are \emph{constrained}
variational problems, in which the solution $u$ is required to be orthogonal
to a given closed subspace $Q\subset V$, a constraint imposed through a Lagrange multiplier. The associated \emph{constrained variational problem} is to find $(u,p)\in V\times Q$ such that
\begin{equation}\label{eq: constrained}
	\begin{aligned}
			B(u, v) + \langle p, v \rangle_V &= F(v) &\forall v \in V,  \\
		\langle u, q \rangle_V &= 0 &\forall q \in Q,
	\end{aligned}
\end{equation}
with $B$ and $F$ as in Definition~\ref{def: variational problem}. If we let
\[
\mathcal A\big((w,r),(v,q)\big) := B(w,v) + \langle r,v\rangle_V - \langle w,q\rangle_V,
\]
then a pair $(u,p)\in V\times Q$ solves \eqref{eq: constrained} precisely when
\begin{equation}\label{eq: block constrained}
	\mathcal A\big((u,p),(v,q)\big)=F(v)\qquad\forall\,(v,q)\in V\times Q.
\end{equation}
This allows us to view the constrained problem as a variational problem on $V\times Q$.

\begin{remark}\label{rem: product norms}
	The product $V\times Q$ carries the Hilbert norm $\|(u,p)\|^2=\|u\|_V^2+\|p\|_V^2$, and Definition~\ref{def: infsup} refers to this norm.  The inf-sup conditions we quote from the literature are stated for the sum norm $\|u\|_V+\|p\|_V$ instead, but the Hilbert norm is at most the sum norm.  Hence a form that is inf-sup stable for the sum norm is inf-sup stable for the Hilbert norm with the same constant.
\end{remark}

To approximate the constrained problem \eqref{eq: constrained} by a
criminal approximation $(V_h,B_h,F_h)$ with maps $\pi_h$ and $i_h$, we need one
more datum: a subspace $Q_h\subseteq V_h$ with which the constraint is
imposed.  The discrete constrained problem is then to find $(u_h,p_h)\in V_h\times Q_h$ such that
\begin{equation}\label{eq: discrete constrained}
	\begin{aligned}
			B_h(u_h, v_h) + \langle p_h, v_h \rangle_{V_h} &= F_h(v_h) &&\forall v_h \in V_h,\\
		\langle u_h, q_h \rangle_{V_h} &= 0 &&\forall q_h \in Q_h.
	\end{aligned}
\end{equation}
As for the continuous problem, if we let
\[
\mathcal A_h\big((w_h,r_h),(v_h,q_h)\big)
:= B_h(w_h,v_h) + \langle r_h,v_h\rangle_{V_h}
- \langle w_h,q_h\rangle_{V_h},
\]
then a pair $(u_h,p_h)\in V_h\times Q_h$ solves the discrete problem \eqref{eq: discrete constrained}
precisely when
\begin{equation}\label{eq: block discrete constrained}
	\mathcal A_h\big((u_h,p_h),(v_h,q_h)\big) = F_h(v_h)
	\qquad \forall\,(v_h,q_h)\in V_h\times Q_h.
\end{equation}
The map $i_h$ need not map $Q_h$ into $Q$. Specialising
Theorem~\ref{thm: third strang} gives a version of the third Strang lemma for
constrained problems.

\begin{theorem}\label{thm: constrained strang}
	Let $(V_h,B_h,F_h)$, $\pi_h$, $i_h$ be a criminal approximation
	of $(V,B,F)$ with $i_h\colon V_h\to V$, and
	let $Q\subseteq V$ and $Q_h\subseteq V_h$ be closed subspaces as above; write $P_{Q_h}\colon V_h\to Q_h$ for the orthogonal
	projection. Assume that $\mathcal A_h$ is inf-sup stable on
		$V_h\times Q_h$ with constant $\gamma_0>0$, and let
	$(u_h,p_h)\in V_h\times Q_h$ be the unique solution of \eqref{eq: discrete constrained}. If $(u,p)\in V\times Q$ solves the constrained
	problem \eqref{eq: constrained} then
	\begin{multline}\label{eq: constrained strang bound}
		\|u - i_h u_h\|_V + \|p - i_h p_h\|_V
		\ \le\
		C \bigg(\|u - i_h \pi_h u\|_V + \|p - i_h P_{Q_h}\pi_h p\|_V \\
		+ \big\|F_h(\cdot) - B_h(\pi_h u, \cdot)
		- \langle P_{Q_h}\pi_h p, \cdot \rangle_{V_h}\big\|_{V_h^*}
		+ \sup_{0\ne q_h \in Q_h}
		\frac{|\langle \pi_h u, q_h \rangle_{V_h}|}{\|q_h\|_{V_h}}
		\bigg),
	\end{multline}
	where $C$ depends only on $\opnorm{i_h}{V_h}{V}$ and $\gamma_0$.
\end{theorem}

\begin{proof}
	We apply Theorem~\ref{thm: third strang} with $V\times Q$ in the role of $V$, $V\times V$ in the role of $E$, and $V_h\times Q_h$ in the role of $V_h$.  The interpolator is $(w,r)\mapsto(\pi_hw,P_{Q_h}\pi_hr)$.  The reconstructor is $(w_h,r_h)\mapsto(i_hw_h,i_hr_h)$, and its norm is at most $\opnorm{i_h}{V_h}{V}$.  Then
	\begin{multline*}
		\big\|(u-i_hu_h,\,p-i_hp_h)\big\|
		\le\big\|(u-i_h\pi_hu,\,p-i_hP_{Q_h}\pi_hp)\big\|\\
		+\opnorm{i_h}{V_h}{V}\,\gamma_0^{-1}\sup_{0\ne(v_h,q_h)\in V_h\times Q_h}
		\frac{\big|F_h(v_h)-\mathcal A_h\big((\pi_hu,P_{Q_h}\pi_hp),(v_h,q_h)\big)\big|}{\|(v_h,q_h)\|}.
	\end{multline*}
	The left-hand side of \eqref{eq: constrained strang bound} is at most $\sqrt2$ times the left-hand side here.  The first term on the right is at most the sum of the two norms it is built from.  Expanding $\mathcal A_h$, the numerator of the supremum is
	\[
	\big|F_h(v_h)-B_h(\pi_hu,v_h)-\langle P_{Q_h}\pi_hp,v_h\rangle_{V_h}+\langle\pi_hu,q_h\rangle_{V_h}\big| .
	\]
	Both $\|v_h\|_{V_h}$ and $\|q_h\|_{V_h}$ are at most $\|(v_h,q_h)\|$, so the supremum is at most the sum of the two consistency terms in \eqref{eq: constrained strang bound}.  Together these give \eqref{eq: constrained strang bound}, with $C$ depending only on $\opnorm{i_h}{V_h}{V}$ and $\gamma_0$.
\end{proof}

\subsection{Hilbert complexes}\label{sec: Hilbert complexes}
This section recalls Hilbert complexes, their Hodge decomposition and Poincar\'e inequality, and poses the Hodge--Dirac and Hodge--Laplace problems on them. Hilbert complexes were defined by Br\"uning and Lesch \cite[Section~2]{BruningLesch1992}; the standard reference for them is Arnold, Falk and Winther \cite[Section~3]{arnold2010finite}.

\subsubsection{Basic definitions}
A \emph{Hilbert complex} $(H^\bullet,V^\bullet,d)$ is an orthogonally
$\Z$-graded Hilbert space $H^\bullet$ together with a \emph{differential} $d$, a
closed, densely defined operator of degree one with domain $V^\bullet$, such that
$d(V^k)\subset V^{k+1}$ and $d(V^\bullet)\subset\ker d$. We write
$Z^\bullet=\ker d$ and $B^\bullet=d(V^\bullet)$. A Hilbert complex is called \emph{closed} if $B^\bullet$ is closed in $H^\bullet$, and \emph{bounded} if $d$ extends to a bounded linear map on all of $H^\bullet$.

The \emph{domain complex} of a Hilbert complex $(H^\bullet,V^\bullet,d)$ is the bounded Hilbert complex $(V^\bullet,V^\bullet,d)$, where $V^\bullet$ is equipped with the \emph{graph inner product}
$$
    \langle u, v\rangle_V = \langle u,v\rangle_H + \langle du, dv\rangle_H,
$$
where $\langle\cdot,\cdot\rangle_H$ is the inner product on $H^\bullet$. In what follows we will often suppress the subscript $H$ from the notation, and we will write $\|\cdot\|_V$ for the norm associated to the graph inner product. For $u\in Z^\bullet$ the graph norm $\|u\|_V$ and the norm $\|u\|$ agree.

\subsubsection{Hodge decomposition and Poincar\'e inequality}\label{sec: hodge decomposition}
The differential $d$ has an adjoint $d^*$ with respect to the inner product of the Hilbert complex. We will denote the domain of definition of $d^*$ by $V_\bullet$, its image by $B_\bullet$ and kernel by $Z_\bullet$.

Hilbert complexes admit a version of the classical Hodge decomposition of differential forms, known as the \emph{weak Hodge decomposition} \cite[Lemma~2.1]{BruningLesch1992}: $H^k$ is the orthogonal direct sum
$$
H^k = \overline{B^k}\oplus \mathfrak{H}^k \oplus \overline{B_k},
$$
where $\mathfrak H^k=Z^k\cap Z_k$ is the \emph{space of harmonic $k$-forms}. If the Hilbert complex is closed we get the \emph{strong Hodge decomposition}
\cite[Section~3.1.3]{arnold2010finite}
$$
H^k = B^k\oplus \mathfrak{H}^k \oplus B_k.
$$
As a direct consequence of the strong Hodge decomposition we have an
isomorphism between the cohomology of the Hilbert complex and the space
$\mathfrak H^\bullet$ of harmonic forms. This isomorphism is induced by the orthogonal projection
$P_{\mathfrak H}\colon H^\bullet\to\mathfrak H^\bullet$ restricted to
$Z^\bullet$, so it depends on the inner product. We also write $P_B$ for the
orthogonal projection onto $B^\bullet$.

Another important corollary of the strong Hodge decomposition is the \emph{Poincar\'e inequality}:
\begin{lemma}[{\cite[Section~3.1.3, (16)]{arnold2010finite}, \cite[Lemma~4]{LeopardiStern2016}}]\label{lem: Poincare inequality}
	Let $(H^\bullet,V^\bullet,d)$ be a closed Hilbert complex.  Then there
	is a constant $c_H\ge 1$, called the Poincar\'e constant, such that
	\[
	\|u\|_V\le c_H\|du\| \qquad
	\forall\, u\in V^\bullet \text{ orthogonal to } Z^\bullet .
	\]
\end{lemma}

\subsubsection{Dirac and Laplace operators}
From the differential and its adjoint we build the \emph{Hodge--Dirac
operator} $d+d^*$ and the \emph{Hodge--Laplace operator}
$$
\Delta =(d+d^*)^2= dd^* + d^*d .
$$
The Hodge--Dirac operator is defined on $V^\bullet\cap V_\bullet$, and the
Hodge--Laplace operator on the smaller domain
$\{\alpha\in V^\bullet\cap V_\bullet:d\alpha\in V_\bullet,\ d^*\alpha\in V^\bullet\}$
\cite[Section~2.3]{Holst2012}. Their common kernel is the space
$\mathfrak H^\bullet$ of harmonic forms.

Given a source $\gamma\in H^\bullet$, we are interested in the equations
$$
(d+d^*)\alpha =\gamma \quad \mbox{and}\quad \Delta \alpha =\gamma
$$
for $\alpha$ in the respective domains. These equations have a solution only if
$\gamma\perp\mathfrak H^\bullet$, and a solution is determined only up to
harmonic forms. The \emph{mixed weak formulation}, which we use, removes both
problems and avoids the domain $V_\bullet$ of $d^*$.  In it the source enters only through the functional $\langle\gamma,\cdot\rangle$ on $V^\bullet$, so any bounded functional can take its place.  We write $(V^\bullet)^*$ and $(V^k)^*$ for the duals of $V^\bullet$ and $V^k$ with their graph norms.

\subsubsection{The abstract Hodge--Dirac problem}
The mixed weak formulation of the Hodge--Dirac problem is posed on the whole graded domain $V^\bullet$.  Given $F\in(V^\bullet)^*$, it asks for $(\alpha,\kappa)\in V^\bullet\times\mathfrak{H}^\bullet$ satisfying
\begin{equation}\label{eq: Hodge Dirac problem}
	\begin{aligned}
		\langle \alpha, \lambda \rangle &= 0
		& \forall \lambda \in \mathfrak{H}^\bullet,\\
		\langle d \alpha, \beta\rangle+ \langle \alpha,d\beta\rangle
		+ \langle \kappa, \beta\rangle &= F(\beta)
		& \forall \beta\in V^{\bullet}.
	\end{aligned}
\end{equation}
We call $F$ the source.  For a square-integrable source $\gamma\in H^\bullet$, that is, $F=\langle\gamma,\cdot\rangle$, which has norm at most $\|\gamma\|$, this is the formulation of Leopardi and Stern \cite[(1)]{LeopardiStern2016}.
The first equation removes the freedom to add a harmonic form to $\alpha$.  In the second equation, $\kappa$ cancels the harmonic part of the source: testing with $\beta\in\mathfrak H^\bullet$ gives $\langle\kappa,\beta\rangle=F(\beta)$, so $\kappa=P_{\mathfrak H}\gamma$ for a square-integrable source $\gamma$. The regularity results of
Section~\ref{sec: abstract regularity} need the solution to lie in
the domain of $d^*$, and need its codifferential.
\begin{lemma}\label{lem: HD solutions in intersection}
	Let $(\alpha,\kappa)$ solve \eqref{eq: Hodge Dirac problem} with source $\gamma\in H^\bullet$. Then $\alpha\in V^\bullet\cap V_\bullet$, with $d^*\alpha=\gamma-\kappa-d\alpha$. Furthermore, if $\gamma\in V_\bullet$, then $d\alpha\in  V^\bullet\cap V_\bullet$, with $d^*d\alpha=d^*\gamma$.
\end{lemma}
\begin{proof}
	Testing the second equation of \eqref{eq: Hodge Dirac problem} with
	$\beta\in V^\bullet$ gives $\langle\alpha,d\beta\rangle=
	\langle\gamma-\kappa-d\alpha,\beta\rangle$, so $\alpha\in V_\bullet$ with
	$d^*\alpha=\gamma-\kappa-d\alpha$. If $\gamma\in V_\bullet$, testing with
	$\beta=d\theta$, $\theta\in V^\bullet$, gives $\langle d\alpha,d\theta\rangle=
	\langle\gamma,d\theta\rangle=\langle d^*\gamma,\theta\rangle$, since
	$dd\theta=0$ and $\kappa\perp B^\bullet$. So $d\alpha\in V_\bullet$ with
	$d^*d\alpha=d^*\gamma$, and $d\alpha\in V^\bullet$ as $dd\alpha=0$.
\end{proof}

As a variational problem on $V^\bullet\times\mathfrak H^\bullet$, \eqref{eq: Hodge Dirac problem} has the bilinear form
\begin{equation}\label{eq: HD form}
	\mathcal D\big((w,r),(v,q)\big):=
	\langle dw,v\rangle+\langle w,dv\rangle
	+\langle r,v\rangle-\langle w,q\rangle .
\end{equation}

Well-posedness of the Hodge--Dirac problem rests on an inf-sup estimate for $\mathcal D$ due to Leopardi and Stern. We use the estimate twice: for every source in $(V^\bullet)^*$, and for the discrete Hodge--Dirac problem, of whose Poincar\'e constant we know only a bound.

\begin{lemma}[{\cite[Theorem~6]{LeopardiStern2016}}]\label{lem: Hodge Dirac stability}
	Let $(H^\bullet,V^\bullet,d)$ be a closed Hilbert complex and let
	$c\ge 1$ be such that $\|v\|_V\le c\|dv\|$ for every $v\in V^\bullet$
	orthogonal to $Z^\bullet$.  Then the form $\mathcal D$ of
	\eqref{eq: HD form} on $V^\bullet\times\mathfrak H^\bullet$, normed by
	$\|v\|_V+\|q\|$, is bounded, and satisfies the first inequality of \eqref{eq: infsup} with a constant depending only on $c$.
\end{lemma}
Leopardi and Stern state the estimate with the Poincar\'e constant of the complex, but their proof uses only the bound $\|v\|_V\le c\|dv\|$ for $v$ orthogonal to $Z^\bullet$; by Lemma~\ref{lem: Poincare inequality}, $c=c_H$ is admissible. The second inequality of \eqref{eq: infsup} follows from the first, as $\big|\mathcal D\big((v,q),(w,r)\big)\big|=\big|\mathcal D\big((w,-r),(v,-q)\big)\big|$. So Theorem~\ref{thm: babuska} shows that the Hodge--Dirac problem is well-posed for every $F\in(V^\bullet)^*$, with the a priori estimate
\begin{equation}\label{eq: HD a priori}
	\|\alpha\|_V+\|\kappa\|\le C\|F\|_{(V^\bullet)^*} ,
\end{equation}
where $C$ depends only on $c_H$. For square-integrable sources this is \cite[Corollary~8]{LeopardiStern2016}.

\subsubsection{The abstract Hodge--Laplace problem}
The mixed formulation of the Hodge--Laplace problem in degree $k$ with source $F\in(V^k)^*$ asks for
$(\eta,\alpha, \kappa)\in V^{k-1}\times V^k\times \mathfrak{H}^k$ satisfying
\begin{equation}\label{eq: Hodge Laplace problem}
    \begin{aligned}
    \langle \alpha, \lambda\rangle &=0 &\forall \lambda\in \mathfrak{H}^k,\\
    \langle \eta, \theta\rangle - \langle\alpha,d\theta\rangle&=0
        & \forall\theta\in V^{k-1},\\
    \langle d\eta, \beta\rangle +\langle d\alpha,d\beta\rangle
        +\langle \kappa, \beta\rangle& = F(\beta)
        & \forall \beta \in V^k .
    \end{aligned}
\end{equation}
Here the second equation weakly enforces $\alpha\in V_k$.  With a square-integrable source $\gamma\in H^k$ this is the mixed formulation of Arnold, Falk and Winther \cite[(18)]{arnold2010finite}.
For our regularity results we need the domain and the
codifferential of each component of the solution. The first assertion of the
following lemma is the interpretation of the mixed formulation in
\cite[Section~3.2.1]{arnold2010finite}.
\begin{lemma}\label{lem: HL solutions in intersection}
	If $(\eta,\alpha,\kappa)$ solves \eqref{eq: Hodge Laplace problem} with source	$\gamma\in H^k$, then $\kappa=P_{\mathfrak H}\gamma$, $\alpha\in V^k\cap V_k$ with $d^*\alpha=\eta$, $\eta\in V^{k-1}\cap V_{k-1}$ with $d^*\eta=0$, and $d\alpha\in V^{k+1}\cap V_{k+1}$ with $d^*d\alpha=\gamma-\kappa-d\eta$. Furthermore, if $\gamma\in V_k$, then
	$d\eta\in V^k\cap V_k$ with $d^*d\eta=d^*\gamma$.
\end{lemma}

\begin{proof}
	The components $\eta$ and $\alpha$ lie in $V^{k-1}$ and $V^k$ by definition.
	The second equation of \eqref{eq: Hodge Laplace problem} says
	$\langle\alpha,d\theta\rangle=\langle\eta,\theta\rangle$ for all
	$\theta\in V^{k-1}$, so $\alpha\in V_k$ with $d^*\alpha=\eta$. Hence $\eta$
	lies in the range of $d^*$, so $\eta\in V_{k-1}$ with $d^*\eta=d^*d^*\alpha=0$.
	The third equation says $\langle d\alpha,d\beta\rangle=
	\langle\gamma-\kappa-d\eta,\beta\rangle$ for all $\beta\in V^k$, so
	$d\alpha\in V_{k+1}$ with $d^*d\alpha=\gamma-\kappa-d\eta$, and
	$d\alpha\in V^{k+1}$ as $dd\alpha=0$. Testing the third equation with
	$\beta\in\mathfrak H^k$, for which $d\beta=0$ and
	$\langle d\eta,\beta\rangle=0$, gives $\kappa=P_{\mathfrak H}\gamma$. Now let
	$\gamma\in V_k$. The identity $d^*d\alpha=\gamma-\kappa-d\eta$ says that
	$\gamma=d\eta+\kappa+d^*d\alpha$ is the Hodge decomposition of $\gamma$
	\cite[Section~3.2.1]{arnold2010finite}. Here $\kappa$ is harmonic and
	$d^*d\alpha$ lies in the range of $d^*$, so both lie in $V_k$ and are
	annihilated by $d^*$. Hence $d\eta=\gamma-\kappa-d^*d\alpha$ lies in $V_k$
	with $d^*d\eta=d^*\gamma$, and $d\eta\in V^k$ as $dd\eta=0$.
\end{proof}

To apply Theorem~\ref{thm: constrained strang}, we view
\eqref{eq: Hodge Laplace problem} as a variational problem on
$V^{k-1}\times V^k\times\mathfrak H^k$, with the bilinear form
\begin{equation}\label{eq: HL form}
	\mathcal B\big((\eta,\alpha,\kappa),(\theta,\beta,\lambda)\big)
	:=\langle\eta,\theta\rangle-\langle\alpha,d\theta\rangle
	+\langle d\eta,\beta\rangle+\langle d\alpha,d\beta\rangle
	+\langle\kappa,\beta\rangle-\langle\alpha,\lambda\rangle.
\end{equation}

The Hodge--Laplace problem is well-posed for every $F\in(V^k)^*$
\cite[Theorem~3.2 and the first remark after it]{arnold2010finite}, with the a priori estimate
\begin{equation}\label{eq: HL a priori}
	\|\eta\|_V+\|\alpha\|_V+\|\kappa\|\le C\|F\|_{(V^k)^*} ,
\end{equation}
where $C$ depends only on the Poincar\'e constant $c_H$
of Lemma~\ref{lem: Poincare inequality}.

\subsection{Approximating Hilbert complexes}\label{sec: approximating hilbert complexes}
This section defines approximations of a Hilbert complex (Definition~\ref{def: approximation of hilbert complex}) and discusses the various error terms associated to them.

\subsubsection{Approximations of Hilbert complexes}
An approximation replaces a Hilbert complex by finite-dimensional complexes,
linked to it by maps in both directions that commute with the differentials.
\begin{definition}\label{def: approximation of hilbert complex}
    An \emph{approximation of a Hilbert complex $(H^\bullet, V^\bullet,d)$ indexed by $\mathcal{I}$} consists, for each $h\in \mathcal{I}$, of a finite-dimensional Hilbert complex $(H_h^\bullet,V_h^\bullet, d_h)$ and chain maps of the domain complexes
    $$
        i_h\colon (V_h^\bullet, d_h)\to (V^\bullet, d)
        \qquad\text{and}\qquad
        \pi_h \colon (V^\bullet, d)\to (V_h^\bullet, d_h).
    $$

    These maps are required to satisfy $\pi_h \circ i_h=\id_{V^\bullet_h}$, and to be bounded uniformly in $h$ for the graph norms: there are constants $C_i$ and $C_\pi$ with $\opnorm{i_h}{V_h^\bullet}{V^\bullet}\le C_i$ and $\opnorm{\pi_h}{V^\bullet}{V_h^\bullet}\le C_\pi$ for all $h\in\mathcal I$.
\end{definition}

\begin{notation}
	We write $\|\cdot\|_h$ and $\langle\cdot,\cdot\rangle_h$ for the norm and
	inner product of $H_h^\bullet$.
\end{notation}

For subcomplexes $V_h^\bullet\subseteq V^\bullet$ with the restricted inner
products and $i_h$ the inclusion, this is the setting of Arnold, Falk and
Winther, in which $\pi_h$ is a bounded cochain projection
\cite[Section~3.3]{arnold2010finite}.  Holst and Stern replaced the
subcomplex by a complex with inner products of its own, linked to
$(H^\bullet,V^\bullet,d)$ by a pair of chain maps $i_h$ and $\pi_h$, bounded
on the Hilbert spaces, with $\pi_h\circ i_h=\id$ \cite[Section~3.1]{Holst2012}.
Definition~\ref{def: approximation of hilbert complex} differs from theirs in
two respects: $\pi_h$ is a chain map of the domain complexes only, as the
cochain projections of Arnold, Falk and Winther are, and the uniform bounds on
$i_h$ and $\pi_h$ are in the graph norms.

The discrete complexes are bounded and closed, being finite dimensional.
In general $i_h$ does not preserve inner products, and $i_h\pi_h$ is not the
orthogonal projection onto $i_h(V_h^\bullet)$. We
approximate the Hodge--Dirac and Hodge--Laplace problems on
$(H^\bullet,V^\bullet,d)$ by the corresponding problems on
$(H^\bullet_h,V_h^\bullet,d_h)$. If $i_h$ were isometric, these would be the
problems on the subcomplex $i_h(V_h^\bullet)$ with the inner products of
$H^\bullet$, the setting of Arnold, Falk and Winther. When it is not, the
discrete problems use inner products of their own, and they form a
criminal approximation of the continuous ones
(Definition~\ref{def: criminal}).

\subsubsection{The metric defect}
Two inner products are available on $V_h^\bullet$: its own,
$\langle\cdot,\cdot\rangle_h$, and the one it inherits from $H^\bullet$
through $i_h$. Their difference is the Hermitian form
\[
\langle x_h,y_h\rangle_h-\langle i_hx_h,i_hy_h\rangle
=\langle(\mathrm{id}-i_h^*i_h)x_h,y_h\rangle_h .
\]
We call its norm on the domain complex the \emph{metric defect}:
\begin{equation}\label{eq: graph metric defect}
	\metricdefect
	:=\sup_{0\ne x_h,y_h\in V_h^\bullet}
	\frac{|\langle(\mathrm{id}-i_h^*i_h)x_h,y_h\rangle_h|}
	{\|x_h\|_{V_h}\,\|y_h\|_{V_h}} .
\end{equation}
Measured on $H_h^\bullet$ instead, the
same form has norm $\opnorm{\mathrm{id}-i_h^*i_h}{H_h^\bullet}{H_h^\bullet}$, the quantity used in
\cite{Holst2012}. Since $\|\cdot\|_h\le\|\cdot\|_{V_h}$,
\begin{equation}\label{eq: metric defect operator norm}
	\metricdefect\ \le\ \opnorm{\mathrm{id}-i_h^*i_h}{H_h^\bullet}{H_h^\bullet} ,
\end{equation}
so a bound on the operator norm, which we will have for finite element
approximations, bounds $\metricdefect$ as well. The converse fails: for the inner product of
discrete exterior calculus on surfaces the operator norm does not tend to zero,
while $\metricdefect$ does \cite{WassermanDEC}.

The maps $i_h$ are uniformly bounded once the metric defect is uniformly
small:

\begin{lemma}\label{lem: metric defect graph bounds}
	Let $(H_h^\bullet,V_h^\bullet,d_h)$ be a finite-dimensional Hilbert
	complex with a chain map $i_h\colon(V_h^\bullet,d_h)\to
	(V^\bullet,d)$ whose metric defect
	\eqref{eq: graph metric defect} satisfies $\metricdefect<\frac12$.  Then
	there is a constant $C$, depending only on $\metricdefect$ and
	nondecreasing in it, such that for all $x_h\in V_h^\bullet$
	\[
		\|x_h\|_{V_h}\le C\,\|i_hx_h\|_V
		\qquad\text{and}\qquad
		\|i_hx_h\|_V\le C\,\|x_h\|_{V_h} .
	\]
\end{lemma}

\begin{proof}
	Let $x_h\in V_h^\bullet$. As
	$i_hd_hx_h=di_hx_h$ and $\langle i_h^*i_hy_h,y_h\rangle_h=\|i_hy_h\|^2$
	for every $y_h\in H_h^\bullet$,
	\[
		\|x_h\|_{V_h}^2-\|i_hx_h\|_V^2
		=\langle(\id-i_h^*i_h)x_h,x_h\rangle_h
		+\langle(\id-i_h^*i_h)d_hx_h,d_hx_h\rangle_h .
	\]
	By \eqref{eq: graph metric defect} the first term is at most
	$\metricdefect\|x_h\|_{V_h}^2$ in absolute value and the second at most
	$\metricdefect\|d_hx_h\|_{V_h}^2$.  As $d_hd_hx_h=0$, the graph norm of
	$d_hx_h$ is $\|d_hx_h\|_h$, which is at most $\|x_h\|_{V_h}$.  So the
	difference is at most $2\metricdefect\|x_h\|_{V_h}^2$ in absolute value, and as
	$2\metricdefect<1$ both inequalities follow, with $C$ depending only on
	$\metricdefect$ and nondecreasing in it.
\end{proof}

\subsubsection{The interpolation defect}
The failure of $i_h\pi_h$ to be the identity at $\omega\in V^\bullet$ is
measured by $\|\omega-i_h\pi_h\omega\|$, which we call the \emph{interpolation
defect} of $\omega$, and by its graph norm $\|\omega-i_h\pi_h\omega\|_V$.

A small metric defect and a bound $\|u-i_h\pi_hu\|_V\le C\|u\|_V$ on the
interpolation defect in the graph norm, equivalently a bound on
$\|i_h\pi_hu\|_V$, make $\{i_h,\pi_h\}_{h\in\mathcal I}$ an approximation, with constants that see the
discrete inner products only through the metric defect.

\begin{lemma}\label{lem: defects approximation}
	Let $(H^\bullet,V^\bullet,d)$ be a Hilbert complex.  For each
	$h\in\mathcal I$, let $(H_h^\bullet,V_h^\bullet,d_h)$ be a finite-dimensional
	Hilbert complex, with chain maps
	$i_h\colon V_h^\bullet\to V^\bullet$ and
	$\pi_h\colon V^\bullet\to V_h^\bullet$ with
	$\pi_h\circ i_h=\id$.  Suppose that the metric defect satisfies
	$\metricdefect\le\nu_0$ for some $\nu_0<\frac12$, and that $\pi_h$ is
	bounded when $V_h^\bullet$ carries the graph norm transported from
	$V^\bullet$ by $i_h$, that is, $\|i_h\pi_hu\|_V\le C_0\|u\|_V$ for
	$u\in V^\bullet$, both uniformly in $h$.  Then
	$\{i_h,\pi_h\}_{h\in\mathcal I}$ is an approximation of
	$(H^\bullet,V^\bullet,d)$, with $C_i$ depending only on $\nu_0$ and
	$C_\pi$ only on $\nu_0$ and $C_0$.
\end{lemma}

\begin{proof}
	The maps are as Definition~\ref{def: approximation of hilbert complex}
	requires, and it remains to bound them for the graph norms.
	Lemma~\ref{lem: metric defect graph bounds}, whose constant is nondecreasing in the metric defect, bounds $i_h$ by a constant $C$ depending only on $\nu_0$.  For
	$u\in V^\bullet$ its first inequality, applied to $x_h=\pi_hu$, gives
	$\|\pi_hu\|_{V_h}\le C\|i_h\pi_hu\|_V\le CC_0\|u\|_V$.
\end{proof}

\subsubsection{The adjoint failure}
From now on, $(H^\bullet,V^\bullet,d)$ is a closed Hilbert complex and
$\{i_h\colon (H_h^\bullet,V_h^\bullet,d_h)\rightleftarrows(H^\bullet,V^\bullet,d):\pi_h\}_{h\in \mathcal{I}}$
an approximation of it.

The failure of $\pi_h$ to be the adjoint of
$i_h$ is measured by $\pi_h-i_h^*$. Its size is bounded by the metric and interpolation defects:

\begin{lemma}\label{lem: adjoint failure}
	For all $h\in\mathcal I$, $\omega\in V^\bullet$ and $y_h\in V_h^\bullet$,
	\begin{equation}\label{eq: non adjoint bound}
		|\langle(\pi_h - i_h^*)\omega,y_h\rangle_h|
		\le \big(C_\pi\metricdefect\,\|\omega\|_V
		+ C_i\|\omega - i_h \pi_h\omega \|\big)\|y_h\|_{V_h} .
	\end{equation}
\end{lemma}

\begin{proof}
	Subtracting and adding the pairing $\langle i_h\pi_h\omega,i_hy_h\rangle$
	gives
	\begin{align*}
		\langle (\pi_h - i_h^*) \omega, y_h \rangle_h
		&= \langle \pi_h \omega, y_h \rangle_h
		- \langle i_h \pi_h \omega, i_h y_h\rangle
		+ \langle i_h \pi_h \omega - \omega, i_h y_h \rangle\\
		&= \langle (\mathrm{id} - i_h^* i_h) \pi_h \omega, y_h \rangle_h
		+ \langle i_h \pi_h \omega - \omega, i_h y_h \rangle .
	\end{align*}
	By \eqref{eq: graph metric defect} and $\|\pi_h\omega\|_{V_h}\le C_\pi\|\omega\|_V$
	the first summand is at most
	$C_\pi\metricdefect\|\omega\|_V\|y_h\|_{V_h}$, and by the
	Cauchy--Schwarz inequality and
	$\|i_hy_h\|\le\|i_hy_h\|_V\le C_i\|y_h\|_{V_h}$ the second is at most
	$C_i\|\omega-i_h\pi_h\omega\|\,\|y_h\|_{V_h}$.
\end{proof}

\subsubsection{The best approximation error}
For $\omega\in V^\bullet$ we call
$\inf_{\beta_h\in V_h^\bullet}\|\omega-i_h\beta_h\|_V$, the graph-norm distance
from $\omega$ to $i_h(V_h^\bullet)$, the \emph{best approximation error} of
$\omega$. The interpolation defect of $\omega$, even in the graph norm, is
bounded by a multiple of it, a property sometimes called
\emph{quasi-optimality}: since $i_h\pi_h$ is a projection onto $i_h(V_h^\bullet)$, one checks that $\omega-i_h\pi_h\omega=(\mathrm{id}-i_h\pi_h)(\omega-i_h\beta_h)$ for every $\beta_h\in V_h^\bullet$, so that
\begin{equation}\label{eq:quasi_V}
    \|\omega - i_h \pi_h \omega\|_V \le C \inf_{\beta_h \in V_h^\bullet} \|\omega - i_h \beta_h\|_V
\end{equation}
for all $h\in\mathcal I$ and $\omega\in V^\bullet$, with $C$ depending only on $C_i$ and $C_\pi$ \cite[(5.7)]{Arnold2018}.

\subsubsection{The harmonic gap}
In both problems the solution is made unique by requiring it to be orthogonal
to the harmonic forms. Because the interpolator and the reconstructor need not
respect $d^*$, they need not preserve harmonic forms.
How far $i_h\pi_h$ is from
fixing them is measured by the \emph{harmonic gap}
\begin{equation}\label{eq: harmonic gap}
	\mu_h:=\sup_{q\in\mathfrak H^\bullet,\ \|q\|=1}\|q-i_h\pi_hq\| .
\end{equation}
For a subcomplex, where $i_h$ is the inclusion, $\mu_h$ is the quantity $\mu$
in the error bound of
\cite[Theorem~3.9]{arnold2010finite}; in the form \eqref{eq: harmonic gap} it
enters the error bound of \cite[Corollary~3.11]{Holst2012}.
The following lemmas compare the harmonic forms with the discrete harmonic
forms in both directions, as needed to compare the orthogonality constraints
of the continuous and the discrete problems. They use the orthogonal decomposition
$Z_h^\bullet=B_h^\bullet\oplus\mathfrak H_h^\bullet$ of the closed discrete
complex, where $B_h^\bullet$ and $Z_h^\bullet$ are the image and the kernel of
$d_h$.

The first lemma shows that a harmonic form $p$ is reconstructed by the discrete
harmonic part $\lambda_h$ of $\pi_hp$. It is analogous to
\cite[Theorem~4.9]{Dodziuk1976}.

\begin{lemma}\label{lem: discrete harmonic projection defect}
	Let $P_{\mathfrak{H}_h}$ be the orthogonal projection of $H_h^\bullet$ onto
	$\mathfrak H_h^\bullet$.  There is a constant $C$, depending only on
	$C_i$ and $C_\pi$, such that for every $p\in\mathfrak H^\bullet$, writing
	$\lambda_h=P_{\mathfrak{H}_h}\pi_h p$,
	\begin{equation}\label{eq:exact_defect_bound}
		\|\pi_h p -\lambda_h\|_h\le C\Big( \metricdefect\,\|p\|
		+\inf_{q_h\in V_h^\bullet}\|p-i_hq_h\|_V\Big)
	\end{equation}
	and
	\begin{equation}\label{eq: harmonic reconstruction bound}
		\|p-i_h\lambda_h\|\le C\Big( \metricdefect\,\|p\|
		+\inf_{q_h\in V_h^\bullet}\|p-i_hq_h\|_V\Big) .
	\end{equation}
\end{lemma}
\begin{proof}
	Let $p\in\mathfrak H^\bullet$.  Since $d_h\pi_hp=\pi_hdp=0$, we have that
	$\pi_hp\in Z_h^\bullet$. Hence $\pi_hp-\lambda_h\in B_h^\bullet=\im d_h$, so
	$\pi_hp-\lambda_h=d_h\eta_h$ for some $\eta_h\in V_h^\bullet$.  As $\lambda_h$
	is orthogonal to $B_h^\bullet$, and $p$ is orthogonal to
	$B^\bullet\ni d(i_h\eta_h)=i_hd_h\eta_h$, we get:
	\begin{align*}
		\|\pi_hp-\lambda_h\|_h^2
		=\langle\pi_hp,d_h\eta_h\rangle_h
		&=\langle(\pi_h-i_h^*)p,d_h\eta_h\rangle_h+\langle p,i_hd_h\eta_h\rangle\\
		&=\langle(\pi_h-i_h^*)p,\pi_hp-\lambda_h\rangle_h .
	\end{align*}
	Dividing by $\|\pi_hp-\lambda_h\|_h$ and applying
	\eqref{eq: non adjoint bound}, with $\|p\|_V=\|p\|$ and
	$\|\pi_hp-\lambda_h\|_{V_h}=\|\pi_hp-\lambda_h\|_h$ (both elements are
	closed), gives
	$$
	\|\pi_h p -\lambda_h\|_h
	\le C\big( \metricdefect\,\|p\|
	+ \|p-i_h\pi_hp\| \big),
	$$
	with $C$ depending only on $C_i$ and $C_\pi$.
	Here $p-i_h\pi_hp$ is closed, so its norm is its graph norm, and
	\eqref{eq:quasi_V} bounds that by a multiple of
	$\inf_{q_h\in V_h^\bullet}\|p-i_hq_h\|_V$. This gives
	\eqref{eq:exact_defect_bound}.

	For \eqref{eq: harmonic reconstruction bound}, the element
	$w_h=\pi_hp-\lambda_h$ is closed, so
	$\|i_hw_h\|\le\|i_hw_h\|_V\le C_i\|w_h\|_{V_h}=C_i\|w_h\|_h$. With the
	triangle inequality, \eqref{eq:exact_defect_bound}, and \eqref{eq:quasi_V}
	for the closed form $p$, this gives
	\[
	\|p-i_h\lambda_h\|
	\le\|p-i_h\pi_hp\|+C_i\|\pi_hp-\lambda_h\|_h
	\le C\Big(\inf_{q_h\in V_h^\bullet}\|p-i_hq_h\|_V
	+\metricdefect\,\|p\|\Big) ,
	\]
	and \eqref{eq: harmonic reconstruction bound} follows after enlarging $C$.
\end{proof}

The infimum in \eqref{eq:exact_defect_bound} and
\eqref{eq: harmonic reconstruction bound} is at most $\mu_h\|p\|$.
Reconstructing a discrete harmonic form gives a closed form, but not necessarily a harmonic
one.  Its exact component can be bounded as follows.

\begin{lemma}\label{lem: exact part of reconstruction}
	If $q_h\in\mathfrak H_h^\bullet$, then $i_hq_h$ is closed. Write
	$i_hq_h=b+z$, with $b\in B^\bullet$ and $z\in\mathfrak H^\bullet$, for its
	strong Hodge decomposition. Then
	\begin{equation}\label{eq: exact part bound}
		\|b\|\le C\big(\mu_h+\metricdefect\big)\,\|q_h\|_h ,
	\end{equation}
	with $C$ depending only on $C_i$ and $C_\pi$.
\end{lemma}

\begin{proof}
	Note $d(i_hq_h)=i_hd_hq_h=0$, and $Z^\bullet$ is the
	orthogonal sum of $B^\bullet$ and $\mathfrak H^\bullet$, which gives the
	decomposition.  The estimate adapts the proof of
	\cite[Theorem~3.5]{arnold2010finite} (see also \cite[Theorem~13]{LeopardiStern2016}) to a reconstructor $i_h$ that is not isometric.

	Since $\pi_hi_h=\mathrm{id}$ we have $i_h\pi_h(i_hq_h)=i_hq_h$, and so
	$b-i_h\pi_hb=-(z-i_h\pi_hz)$.  Since $b\in B^\bullet$, the element
	$\pi_hb$ lies in $B_h^\bullet$ and is orthogonal to $q_h$, while
	$i_h\pi_hb$ lies in $B^\bullet$ and is orthogonal to $z$; hence
	\[
	\langle i_h\pi_hb,b\rangle
	=\langle i_h\pi_hb,i_hq_h\rangle
	=\langle \pi_hb,i_h^*i_hq_h\rangle_h
	=-\langle \pi_hb,(\mathrm{id}-i_h^*i_h)q_h\rangle_h.
	\]
	Using this together with $\|\pi_hb\|_{V_h}\le C_\pi\|b\|$ and
	$\|q_h\|_{V_h}=\|q_h\|_h$ gives
	\[
	\|b\|^2=\langle b-i_h\pi_hb,b\rangle+\langle i_h\pi_hb,b\rangle
	\le\|z-i_h\pi_hz\|\,\|b\|
	+C_\pi\metricdefect\,\|b\|\,\|q_h\|_h .
	\]
	The form $z$ is harmonic and is an orthogonal component of $i_hq_h$, so the
	definition of $\mu_h$ gives
	$\|z-i_h\pi_hz\|\le\mu_h\|z\|\le\mu_h\|i_hq_h\|\le C_i\mu_h\|q_h\|_h$,
	where the last step uses that $q_h$ is closed, so that
	$\|i_hq_h\|\le\|i_hq_h\|_V\le C_i\|q_h\|_{V_h}=C_i\|q_h\|_h$.
	Dividing by $\|b\|$, which we may assume nonzero, gives
	\eqref{eq: exact part bound}.
\end{proof}

We can now bound how far the interpolant of an element orthogonal to
$\mathfrak H^\bullet$ is from satisfying the discrete constraint, orthogonality
to $\mathfrak H_h^\bullet$.

\begin{prop}\label{prop: constraint defect}
	There is a constant $C$, depending only on $C_i$ and
	$C_\pi$, with the following property.  If $\alpha \in V^\bullet$ is
	orthogonal to $\mathfrak{H}^\bullet$, then for every
	$q_h\in\mathfrak H_h^\bullet$
	\begin{multline}\label{eq:constraint_pairing_bound}
		|\langle \pi_h \alpha, q_h \rangle_h|
		\le C\Big( \metricdefect\,\|\alpha\|_V
		+ \|\alpha - i_h \pi_h\alpha\|\\
		+ \mu_h
		\inf_{v_h\in V_h^\bullet}\|P_B\alpha-i_hv_h\|_V\Big)\,
		\|q_h\|_h .
	\end{multline}
\end{prop}

\begin{proof}
	Let $q_h\in\mathfrak H_h^\bullet$.  We split
	\begin{equation}\label{eq: pairing split}
		\langle \pi_h \alpha, q_h \rangle_h
		= \langle (\pi_h - i_h^*) \alpha, q_h \rangle_h
		+ \langle \alpha, i_h q_h \rangle ,
	\end{equation}
	and \eqref{eq: non adjoint bound} bounds the first term by the first two
	terms of \eqref{eq:constraint_pairing_bound}, since
	$\|q_h\|_{V_h}=\|q_h\|_h$.

	It remains to bound $\langle\alpha,i_hq_h\rangle$. Because $\alpha$ is orthogonal to $\mathfrak H^\bullet$, the strong Hodge decomposition gives $\alpha-P_B\alpha\in B_\bullet$, which is orthogonal to $Z^\bullet$.  As $i_h$ commutes with the differentials, $i_hq_h$ lies
	in $Z^\bullet$.  Hence
	\begin{equation}\label{eq: pairing with exact part}
		\langle\alpha,i_hq_h\rangle=\langle P_B\alpha,i_hq_h\rangle.
	\end{equation}
	Subtracting and adding $i_h\pi_hP_B\alpha$ gives
	\[
	\langle P_B\alpha,i_hq_h\rangle
	=\langle(\mathrm{id}-i_h\pi_h)P_B\alpha,i_hq_h\rangle
	+\langle i_h\pi_hP_B\alpha,i_hq_h\rangle .
	\]
	In the last term we introduce the metric defect by splitting
	\[
	\langle i_h\pi_hP_B\alpha,i_hq_h\rangle
	=\langle\pi_hP_B\alpha,q_h\rangle_h
	-\langle(\mathrm{id}-i_h^*i_h)\pi_hP_B\alpha,q_h\rangle_h .
	\]
	The first of these two terms vanishes.  Indeed $P_B\alpha$ is exact and
	$\pi_h$ is a chain map, so $\pi_h P_B\alpha$ lies in $B_h^\bullet$,
	which is orthogonal to $\mathfrak{H}_h^\bullet$, where $q_h$ lives. The second term is at most
	$C_\pi\metricdefect\,\|\alpha\|_V\,\|q_h\|_h$, which is the first term of
	\eqref{eq:constraint_pairing_bound}.  This uses
	\eqref{eq: graph metric defect} together with three facts: the
	form $P_B\alpha$ is closed, so $\|P_B\alpha\|_V=\|P_B\alpha\|$ and hence
	$\|\pi_hP_B\alpha\|_{V_h}\le C_\pi\|P_B\alpha\|$; the projection is a
	contraction, so $\|P_B\alpha\|\le\|\alpha\|_V$; and $q_h$ is closed, so
	$\|q_h\|_{V_h}=\|q_h\|_h$.

	This leaves us with $\langle(\mathrm{id}-i_h\pi_h)P_B\alpha,i_hq_h\rangle$.
	Write $P_B\alpha=dx$ with $x\in V^\bullet$. Then
	$i_h\pi_hP_B\alpha=d(i_h\pi_hx)$, so both $P_B\alpha$ and
	$i_h\pi_hP_B\alpha$ lie in $B^\bullet$, and
	$(\mathrm{id}-i_h\pi_h)P_B\alpha$ does too.  Write $i_hq_h=b+z$ with
	$b\in B^\bullet$ and $z\in\mathfrak H^\bullet$, as in
	Lemma~\ref{lem: exact part of reconstruction}.  Then
	$(\mathrm{id}-i_h\pi_h)P_B\alpha$ is orthogonal to $z$, and only $b$ contributes,
	\[
	\langle(\mathrm{id}-i_h\pi_h)P_B\alpha,i_hq_h\rangle
	=\langle(\mathrm{id}-i_h\pi_h)P_B\alpha,b\rangle ,
	\]
	and by the Cauchy--Schwarz inequality it suffices to bound the two factors. Being a difference of closed forms, $(\mathrm{id}-i_h\pi_h)P_B\alpha$ has norm equal to its
	graph norm, so \eqref{eq:quasi_V} bounds it by a multiple of
	$\inf_{v_h\in V_h^\bullet}\|P_B\alpha-i_hv_h\|_V$,
	and \eqref{eq: exact part bound} bounds $\|b\|$ by
	a multiple of $(\mu_h+\metricdefect)\|q_h\|_h$.  In the product of the two, the
	part carrying $\metricdefect$ is at most a multiple of
	$\metricdefect\,\|P_B\alpha\|\,\|q_h\|_h\le\metricdefect\,\|\alpha\|_V\,\|q_h\|_h$, since
	the infimum is at most $\|P_B\alpha\|_V=\|P_B\alpha\|$ (take $v_h=0$), so
	it joins the first term of \eqref{eq:constraint_pairing_bound} and $\mu_h$
	alone is left against the infimum.  Collecting the terms gives
	\eqref{eq:constraint_pairing_bound}.
\end{proof}

\subsection{Abstract error bounds}\label{sec: abstract error bounds}
We now bound the errors of the Hodge--Dirac and Hodge--Laplace problems on an approximation by applying Theorem~\ref{thm: constrained strang}, with the error terms of Section~\ref{sec: approximating hilbert complexes}. The error in the discrete source enters through the consistency error of that theorem, and the metric defect through every term that involves the discrete inner products. The harmonic forms enter through how well $i_h\pi_h$ reproduces them, the harmonic gap. Sections~\ref{sec: abstract regularity} and~\ref{sec: abstract rates} turn the bounds into rates.

\subsubsection{Approximating the abstract Hodge--Dirac problem}
The first step for the Hodge--Dirac problem is the following lemma: its discrete bilinear form is inf-sup stable, uniformly in $h$. For a subcomplex with the
restricted inner product, where $i_h$ is the inclusion, it is
\cite[Theorem~10]{LeopardiStern2016}.

\begin{lemma}\label{lem: discrete Hodge Dirac stability}
	Let $(H^\bullet,V^\bullet,d)$ be a closed Hilbert complex with
	Poincar\'e constant $c_H$ and let
	$\{i_h\colon (H_h^\bullet,V_h^\bullet,d_h)\rightleftarrows(H^\bullet,V^\bullet,d):\pi_h\}_{h\in\mathcal I}$
	be an approximation.  Then for every $h\in\mathcal I$ the form
	\begin{equation}\label{eq: discrete HD form}
		\mathcal D_h\big((w_h,r_h),(v_h,q_h)\big):=	\langle d_hw_h,v_h\rangle_h+\langle w_h,d_hv_h\rangle_h
		+\langle r_h,v_h\rangle_h-\langle w_h,q_h\rangle_h
	\end{equation}
	on $V_h^\bullet\times\mathfrak H_h^\bullet$, normed with the graph norm of $V_h^\bullet$, is inf-sup stable with a constant
	$\gamma_0>0$ depending only on $c_H$, $C_i$ and $C_\pi$.
\end{lemma}

\begin{proof}
	By \cite[Corollary~3.6]{Holst2012} the Poincar\'e constant of the closed discrete complex is bounded in terms of that of the domain complex, and the latter may be taken to be $c_H$:
	\[
	\|v_h\|_{V_h}\le c_H\opnorm{\pi_h^k}{V^k}{V_h^k}\,\opnorm{i_h^{k+1}}{V_h^{k+1}}{V^{k+1}}\,\|d_hv_h\|_h
	\]
	for $v_h\in V^k_h$ orthogonal to $Z_h^k$; orthogonality to $Z_h^k$ is the same for both inner products, since $Z_h^k$ consists of closed elements. The two operator norms are at most $C_\pi$ and $C_i$. Moreover $C_\pi C_i\ge\opnorm{\pi_h}{V^\bullet}{V_h^\bullet}\,\opnorm{i_h}{V_h^\bullet}{V^\bullet}\ge\opnorm{\pi_hi_h}{V_h^\bullet}{V_h^\bullet}=1$ and $c_H\ge1$, so $c_HC_\pi C_i\ge1$, as Lemma~\ref{lem: Hodge Dirac stability} requires. The degree-wise bound gives the same bound on all of $V_h^\bullet$. So Lemma~\ref{lem: Hodge Dirac stability} gives the first inequality of \eqref{eq: infsup} for the sum norm, with a constant depending only on $c_H$, $C_i$ and $C_\pi$. The second inequality follows from the first, as $\big|\mathcal D_h\big((v_h,q_h),(w_h,r_h)\big)\big|=\big|\mathcal D_h\big((w_h,-r_h),(v_h,-q_h)\big)\big|$.
\end{proof}

Before applying Theorem~\ref{thm: constrained strang}, we estimate its
consistency error for the Hodge--Dirac problem. Testing the continuous problem against reconstructed discrete
forms bounds it by the source error, the metric defect and the interpolation
defect of the solution.

\begin{lemma}\label{lem: Hodge Dirac consistency}
	Let $F\in(V^\bullet)^*$, let $(\alpha,\kappa)$ solve the Hodge--Dirac problem
	\eqref{eq: Hodge Dirac problem} on $(H^\bullet,V^\bullet,d)$ with
	source $F$, let $F_h\in(V_h^\bullet)^*$ and let
	$\lambda_h=P_{\mathfrak{H}_h}\pi_h\kappa$.  Then
	\begin{multline}\label{eq: HD consistency}
		\sup_{0\ne\beta_h\in V_h^\bullet}
		\frac{\big|F_h(\beta_h)
			-\mathcal D_h\big((\pi_h\alpha,\lambda_h),(\beta_h,0)\big)\big|}
		{\|\beta_h\|_{V_h}}\\
		\le\|F_h-F\circ i_h\|_{(V_h^\bullet)^*}
		+ C\Big(\metricdefect\big(\|\alpha\|_V+\|\kappa\|\big)
		+\inf_{q_h\in V_h^\bullet}\|\kappa-i_hq_h\|_V
		+\|\alpha-i_h\pi_h\alpha\|_V\Big),
	\end{multline}
	where the dual norm on the right is taken with respect to
	$\|\cdot\|_{V_h}$, and $C$ depends only on $C_i$ and $C_\pi$.
\end{lemma}

\begin{proof}
	Let $\beta_h\in V_h^\bullet$ be nonzero. Recalling the definition \eqref{eq: discrete HD form} of $\mathcal D_h$, we have
	\begin{multline}\label{eq: HD consistency expanded}
		F_h(\beta_h)
		-\mathcal D_h\big((\pi_h\alpha,\lambda_h),(\beta_h,0)\big)\\
		=F_h(\beta_h)
		-\langle d_h\pi_h\alpha,\beta_h\rangle_h
		-\langle \pi_h\alpha,d_h\beta_h\rangle_h
		-\langle\lambda_h,\beta_h\rangle_h .
	\end{multline}
	Testing the second equation of
	\eqref{eq: Hodge Dirac problem} with $\beta=i_h\beta_h$, and using
	$d(i_h\beta_h)=i_hd_h\beta_h$, gives
	\[
	\langle i_h^*d\alpha,\beta_h\rangle_h
	+\langle i_h^*\alpha,d_h\beta_h\rangle_h
	+\langle i_h^*\kappa,\beta_h\rangle_h
	-F(i_h\beta_h)=0 .
	\]
	Adding this to the right-hand side of \eqref{eq: HD consistency expanded}, and using
	$d_h\pi_h\alpha=\pi_hd\alpha$, gives
	\[
	F_h(\beta_h)-F(i_h\beta_h)
	-\langle(\pi_h-i_h^*)d\alpha,\beta_h\rangle_h
	-\langle(\pi_h-i_h^*)\alpha,d_h\beta_h\rangle_h
	-\langle\lambda_h-i_h^*\kappa,\beta_h\rangle_h ,
	\]
	whose first two summands are by definition at most
	$\|F_h-F\circ i_h\|_{(V_h^\bullet)^*}\|\beta_h\|_{V_h}$.  Split the last summand as
	\[
	\langle\lambda_h-i_h^*\kappa,\beta_h\rangle_h
	=\langle\lambda_h-\pi_h\kappa,\beta_h\rangle_h
	+\langle(\pi_h-i_h^*)\kappa,\beta_h\rangle_h ,
	\]
	and bound the first of these by $\|\lambda_h-\pi_h\kappa\|_h\|\beta_h\|_h$,
	which by \eqref{eq:exact_defect_bound} with $\kappa$ in the role of $p$ is at
	most $C\big(\metricdefect\|\kappa\|
	+\inf_{q_h\in V_h^\bullet}\|\kappa-i_hq_h\|_V\big)\|\beta_h\|_{V_h}$.

	The three remaining pairings are $\langle(\pi_h-i_h^*)\omega,y_h\rangle_h$
	with $(\omega,y_h)$ equal to $(d\alpha,\beta_h)$, $(\alpha,d_h\beta_h)$ and
	$(\kappa,\beta_h)$. Inequality \eqref{eq: non adjoint bound} bounds each by a
	metric defect term and an interpolation term, times
	$\|y_h\|_{V_h}\le\|\beta_h\|_{V_h}$; for $y_h=d_h\beta_h$ this holds because
	$d_h\beta_h$ is closed, so
	$\|d_h\beta_h\|_{V_h}=\|d_h\beta_h\|_h\le\|\beta_h\|_{V_h}$.  The metric
	defect terms the three pairings contribute are multiples of $\metricdefect\|d\alpha\|_V$,
	$\metricdefect\|\alpha\|_V$ and $\metricdefect\|\kappa\|_V$, and
	$\|d\alpha\|_V=\|d\alpha\|\le\|\alpha\|_V$ and $\|\kappa\|_V=\|\kappa\|$
	because $d\alpha$ and $\kappa$ are closed.  Their interpolation terms are
	multiples of $\|d\alpha-i_h\pi_hd\alpha\|$, $\|\alpha-i_h\pi_h\alpha\|$ and
	$\|\kappa-i_h\pi_h\kappa\|$; the first is $\|d(\alpha-i_h\pi_h\alpha)\|$,
	hence at most $\|\alpha-i_h\pi_h\alpha\|_V$, because $i_h$ and $\pi_h$ commute with
	the differentials, and the last is the graph norm of $\kappa-i_h\pi_h\kappa$,
	since $\kappa$ is closed, hence at most
	a multiple of $\inf_{q_h\in V_h^\bullet}\|\kappa-i_hq_h\|_V$ by
	\eqref{eq:quasi_V}.

	Every summand is therefore at most $\|\beta_h\|_{V_h}$ times one of the terms
	on the right of \eqref{eq: HD consistency}.  Dividing by $\|\beta_h\|_{V_h}$
	and taking the supremum gives \eqref{eq: HD consistency}.
\end{proof}

We can now apply Theorem~\ref{thm: constrained strang} to the Hodge--Dirac
problem. The source error in the bound, $\|F_h-F\circ i_h\|_{(V_h^\bullet)^*}$, is measured in the norm dual to
$\|\cdot\|_{V_h}$.

\begin{theorem}\label{thm: Hodge Dirac abstract approximation theorem}
	Let $(H^\bullet,V^\bullet,d)$ be a closed Hilbert complex with
	Poincar\'e constant $c_H$ and let
	$\{i_h\colon (H_h^\bullet,V_h^\bullet,d_h)\rightleftarrows(H^\bullet,V^\bullet,d):\pi_h\}_{h\in\mathcal I}$
	be an approximation.
	Let $F\in(V^\bullet)^*$ and, for each $h$, let
	$F_h\in(V_h^\bullet)^*$.  Then the Hodge--Dirac problem
	\eqref{eq: Hodge Dirac problem} on $(H_h^\bullet,V_h^\bullet,d_h)$
	with source $F_h$ has a unique solution
	$(\alpha_h,\kappa_h)\in V_h^\bullet\times\mathfrak H_h^\bullet$.
	Furthermore, if $(\alpha,\kappa)\in V^\bullet\times\mathfrak H^\bullet$
	solves \eqref{eq: Hodge Dirac problem} on $(H^\bullet,V^\bullet,d)$
	with source $F$, then
	\begin{multline}\label{eq: Hodge Dirac bound}
		\|\alpha - i_h \alpha_h\|_V + \|\kappa - i_h \kappa_h\|
		\le C \bigg(
		\inf_{\beta_h \in V_h^\bullet} \|\alpha - i_h \beta_h\|_V
		+ \inf_{q_h \in V_h^\bullet} \|\kappa - i_h q_h\|_V\\
		+ \mu_h
		\inf_{v_h \in V_h^\bullet} \|P_B\alpha - i_h v_h\|_V
		+ \|F_h - F\circ i_h\|_{(V_h^\bullet)^*}
		+ \metricdefect\, \|F\|_{(V^\bullet)^*}
		\bigg),
	\end{multline}
	where $\mu_h$ is the harmonic gap \eqref{eq: harmonic gap}
	and $C$ depends only on $c_H$, $C_i$ and $C_\pi$.
\end{theorem}

\begin{proof}
	Fix $h\in\mathcal I$ and put $\lambda_h:=P_{\mathfrak{H}_h}\pi_h\kappa$. We apply Theorem~\ref{thm: constrained strang}. Take $V^\bullet$ and $V_h^\bullet$ with their graph inner products in the roles of $V$ and $V_h$, and the constraint spaces $Q=\mathfrak H^\bullet$ and $Q_h=\mathfrak H_h^\bullet$. The form $(w,v)\mapsto\langle dw,v\rangle+\langle w,dv\rangle$ plays the role of $B$, and the same expression in $d_h$ and $\langle\cdot,\cdot\rangle_h$ that of $B_h$ [$F$ and $F_h$ play their own roles]. Then the constrained problems \eqref{eq: constrained} and \eqref{eq: discrete constrained} are the Hodge--Dirac problems \eqref{eq: Hodge Dirac problem} on the continuous and the discrete complex. Indeed, only their constraint and multiplier terms could differ, and in each of these one argument is harmonic, hence closed. For a closed $z_h$ and any $x_h$ in $V_h^\bullet$ we have $\langle x_h,z_h\rangle_{V_h}=\langle x_h,z_h\rangle_h+\langle d_hx_h,d_hz_h\rangle_h=\langle x_h,z_h\rangle_h$, and likewise in $V^\bullet$.

	The hypotheses of Theorem~\ref{thm: constrained strang} on $\pi_h$ and $i_h$ are those of an approximation. The form $\mathcal A_h$ is $\mathcal D_h$, inf-sup stable by Lemma~\ref{lem: discrete Hodge Dirac stability} with a constant	$\gamma_0$ depending only on $c_H$, $C_i$ and $C_\pi$. The discrete	problem therefore has a unique solution $(\alpha_h,\kappa_h)$, which is the	first assertion.

	Theorem~\ref{thm: constrained strang} provides the bound \eqref{eq: constrained strang bound}. In it $P_{Q_h}\pi_h\kappa=\lambda_h$, because orthogonality to $\mathfrak H_h^\bullet$ is the same for $\langle\cdot,\cdot\rangle_{V_h}$ and $\langle\cdot,\cdot\rangle_h$. By the identity above, $\langle\lambda_h,\beta_h\rangle_{V_h}=\langle\lambda_h,\beta_h\rangle_h$, and $\langle\pi_h\alpha,q_h\rangle_{V_h}=\langle\pi_h\alpha,q_h\rangle_h$ for $q_h\in\mathfrak H_h^\bullet$. Three terms involve closed elements only, so their graph norms are the norms of $H^\bullet$ and $H_h^\bullet$: those of $\kappa-i_h\lambda_h$, of $\kappa-i_h\kappa_h$ and of $q_h$. The estimate becomes
	\begin{multline}\label{eq: HD strang instance}
		\|\alpha - i_h \alpha_h\|_V + \|\kappa - i_h \kappa_h\|
		\le C\bigg(\|\alpha - i_h \pi_h \alpha\|_V + \|\kappa - i_h \lambda_h\| \\
		+ \sup_{0\ne\beta_h\in V_h^\bullet}
		\frac{\big|F_h(\beta_h)
			-\mathcal D_h\big((\pi_h\alpha,\lambda_h),(\beta_h,0)\big)\big|}
		{\|\beta_h\|_{V_h}}
		+ \sup_{0\ne q_h\in\mathfrak H_h^\bullet}
		\frac{|\langle\pi_h\alpha,q_h\rangle_h|}{\|q_h\|_h}\bigg).
	\end{multline}
	Each of the four terms is bounded by one of the preceding results.  The
	first is at most a multiple of $\inf_{\beta_h}\|\alpha-i_h\beta_h\|_V$ by
	\eqref{eq:quasi_V}. The second is bounded by
	\eqref{eq: harmonic reconstruction bound}, with $\kappa$ in the role of $p$.
	The third is bounded by
	\eqref{eq: HD consistency} from Lemma~\ref{lem: Hodge Dirac consistency}.  The fourth is bounded by
	\eqref{eq:constraint_pairing_bound}, which applies because $\alpha$ is
	orthogonal to $\mathfrak H^\bullet$ by the first equation of
	\eqref{eq: Hodge Dirac problem}, and which is valid for every
	$q_h\in\mathfrak H_h^\bullet$.

	Combining the four bounds, the terms carrying $\metricdefect$ alone come
	with the factor $\|\alpha\|_V+\|\kappa\|$, which is at most a multiple of
	$\|F\|_{(V^\bullet)^*}$ by \eqref{eq: HD a priori}, and $\mu_h$ appears only as a
	factor of the best approximation error of $P_B\alpha$.  This is
	\eqref{eq: Hodge Dirac bound}, with a constant built from $C_i$, $C_\pi$
	and the constants of the cited results, hence depending only on $c_H$, $C_i$ and
	$C_\pi$.
\end{proof}

\begin{remark}\label{rem: graph metric defect}
	By \eqref{eq: metric defect operator norm}, the bound
	\eqref{eq: Hodge Dirac bound} also holds with the operator norm
	$\opnorm{\mathrm{id}-i_h^*i_h}{H_h^\bullet}{H_h^\bullet}$ in place of $\metricdefect$;
	Section~\ref{sec: dec application} explains why the smaller quantity
	matters.
\end{remark}

\subsubsection{Approximating the abstract Hodge--Laplace problem}
The Hodge--Laplace problem is treated in the same way as the Hodge--Dirac
problem: inf-sup stability of the discrete form, a bound on the consistency
error, and Theorem~\ref{thm: constrained strang}.

\begin{lemma}[{\cite[Theorem~3.9]{Holst2012}}]\label{lem: discrete Hodge Laplace stability}
	Let $(H^\bullet,V^\bullet,d)$ be a closed Hilbert complex with
	Poincar\'e constant $c_H$ and let
	$\{i_h\colon (H_h^\bullet,V_h^\bullet,d_h)\rightleftarrows(H^\bullet,V^\bullet,d):\pi_h\}_{h\in\mathcal I}$
	be an approximation.  Then for every $h$ the form $\mathcal B_h$
	obtained from $\mathcal B$ of \eqref{eq: HL form}
	by replacing $d$, $\langle\cdot,\cdot\rangle$ and $\mathfrak H^k$ by
	$d_h$, $\langle\cdot,\cdot\rangle_h$ and $\mathfrak H_h^k$, on
	$V_h^{k-1}\times V_h^k\times\mathfrak H_h^k$ with the Hilbert norm built
	from the graph norms, is inf-sup
	stable, with a constant $\gamma_0>0$ depending only on $c_H$, $C_i$ and
	$C_\pi$.
\end{lemma}
\begin{proof}
	By \cite[Corollary~3.6]{Holst2012} the discrete complex has Poincar\'e
	constant at most $c_HC_\pi C_i$, as in the proof of
	Lemma~\ref{lem: discrete Hodge Dirac stability}; that corollary uses only
	the graph norms of $i_h$ and $\pi_h$. By \cite[Theorem~2.10]{Holst2012},
	which is \cite[Theorem~3.2]{arnold2010finite}, applied to the closed complex
	$(H_h^\bullet,V_h^\bullet,d_h)$, the form $\mathcal B_h$ satisfies the first
	inequality of \eqref{eq: infsup}, for the sum norm, with a constant depending
	only on this Poincar\'e constant. The first
	inequality holds for the Hilbert norm by the argument of
	Remark~\ref{rem: product norms}, and the second follows from the first, as
	$\big|\mathcal B_h\big((\theta_h,\beta_h,\lambda_h),(\eta_h,\alpha_h,\kappa_h)\big)\big|
	=\big|\mathcal B_h\big((\eta_h,-\alpha_h,\kappa_h),(\theta_h,-\beta_h,\lambda_h)\big)\big|$.
\end{proof}

The consistency error is estimated as for the Hodge--Dirac problem:
\begin{lemma}\label{lem: Hodge Laplace consistency}
	Let $F\in(V^k)^*$, let $(\eta,\alpha,\kappa)$ solve
	\eqref{eq: Hodge Laplace problem} with source $F$, let
	$F_h\in(V_h^k)^*$ and let $\lambda_h=P_{\mathfrak{H}_h}\pi_h\kappa$.
	Then
	\begin{multline}\label{eq: HL consistency}
		\sup_{(\theta_h,\beta_h)\ne(0,0)}
		\frac{\big|F_h(\beta_h)
			-\mathcal B_h\big((\pi_h\eta,\pi_h\alpha,\lambda_h),
			(\theta_h,\beta_h,0)\big)\big|}
		{\|(\theta_h,\beta_h)\|}\\
		\le\|F_h-F\circ i_h\|_{(V_h^k)^*}
		+C\Big(\metricdefect
		\big(\|\eta\|_V+\|\alpha\|_V+\|\kappa\|\big)
		+\inf_{q_h\in V_h^k}\|\kappa-i_hq_h\|_V\\
		+\|\eta-i_h\pi_h\eta\|_V+\|\alpha-i_h\pi_h\alpha\|_V\Big),
	\end{multline}
	where $(\theta_h,\beta_h)$ runs over $V_h^{k-1}\times V_h^k$, the dual
	norm on the right is taken with respect to $\|\cdot\|_{V_h}$, and $C$
	depends only on $C_i$ and $C_\pi$.
\end{lemma}

\begin{proof}
	Let $(\theta_h,\beta_h)\in V_h^{k-1}\times V_h^k$ be nonzero.  Test the second equation of \eqref{eq: Hodge Laplace problem} with $\theta=i_h\theta_h$ and the third with $\beta=i_h\beta_h$, using $d i_h=i_hd_h$.  Adding the results to the expansion of $F_h(\beta_h)-\mathcal B_h\big((\pi_h\eta,\pi_h\alpha,\lambda_h),(\theta_h,\beta_h,0)\big)$, and using $d_h\pi_h=\pi_hd$, gives
	\begin{multline*}
		F_h(\beta_h)-F(i_h\beta_h)
		-\langle(\pi_h-i_h^*)\eta,\theta_h\rangle_h
		+\langle(\pi_h-i_h^*)\alpha,d_h\theta_h\rangle_h\\
		-\langle(\pi_h-i_h^*)d\eta,\beta_h\rangle_h
		-\langle(\pi_h-i_h^*)d\alpha,d_h\beta_h\rangle_h
		-\langle\lambda_h-i_h^*\kappa,\beta_h\rangle_h .
	\end{multline*}
	The terms are bounded as in the proof of Lemma~\ref{lem: Hodge Dirac consistency}.  The first two are at most $\|F_h-F\circ i_h\|_{(V_h^k)^*}\|\beta_h\|_{V_h}$.  The last is split and bounded as in that proof, using \eqref{eq:exact_defect_bound}, \eqref{eq: non adjoint bound} and \eqref{eq:quasi_V}.  The remaining four are of the form $\langle(\pi_h-i_h^*)\omega,y_h\rangle_h$, with $(\omega,y_h)$ equal to $(\eta,\theta_h)$, $(\alpha,d_h\theta_h)$, $(d\eta,\beta_h)$ and $(d\alpha,d_h\beta_h)$, and \eqref{eq: non adjoint bound} applies to each.  In every case $\|y_h\|_{V_h}\le\|(\theta_h,\beta_h)\|$, because $d_h\theta_h$ and $d_h\beta_h$ are closed.  Moreover $\|d\eta\|_V\le\|\eta\|_V$ and $\|d\eta-i_h\pi_hd\eta\|\le\|\eta-i_h\pi_h\eta\|_V$, because $d\eta$ is closed and $i_h$ and $\pi_h$ commute with the differentials; likewise for $\alpha$.  Dividing by $\|(\theta_h,\beta_h)\|$ and taking the supremum gives \eqref{eq: HL consistency}.
\end{proof}

As for the Hodge--Dirac problem, Theorem~\ref{thm: constrained strang} now gives the error bound.
\begin{theorem}\label{thm:Hodge Laplace bound}
	Let $(H^\bullet,V^\bullet,d)$ be a closed Hilbert complex with
	Poincar\'e constant $c_H$ and let
	$\{i_h\colon (H_h^\bullet,V_h^\bullet,d_h)\rightleftarrows(H^\bullet,V^\bullet,d):\pi_h\}_{h\in \mathcal{I}}$
	be an approximation.
	Let $F\in(V^k)^*$ and, for each $h$, let $F_h\in(V_h^k)^*$.  Then the
	Hodge--Laplace problem \eqref{eq: Hodge Laplace problem} on
	$(H_h^\bullet,V_h^\bullet,d_h)$ with source $F_h$ has a unique
	solution $(\eta_h,\alpha_h,\kappa_h)$.  Furthermore, if
	$(\eta,\alpha,\kappa)$ solves \eqref{eq: Hodge Laplace problem} on
	$(H^\bullet,V^\bullet,d)$ with source $F$, then
	\begin{multline}\label{eq: Hodge Laplace bound}
		\|\eta -i_h \eta_h\|_V+\|\alpha-i_h\alpha_h\|_V+\|\kappa-i_h\kappa_h\|
		\le C\Big(\inf_{\theta_h\in V_h^{k-1}}\|\eta-i_h\theta_h\|_V\\
		+\inf_{\beta_h\in V_h^{k}}\|\alpha-i_h\beta_h\|_V
		+\inf_{q_h\in V_h^k}\|\kappa-i_hq_h\|_V\\
		+\mu_h\inf_{v_h\in V_h^k}\|P_B\alpha-i_hv_h\|_V
		+ \|F_h-F\circ i_h\|_{(V_h^k)^*}
		+ \metricdefect\, \|F\|_{(V^k)^*}\Big),
	\end{multline}
	where $\mu_h$ is the harmonic gap \eqref{eq: harmonic gap} and $C$ depends
	only on $c_H$, $C_i$ and $C_\pi$.
\end{theorem}

\begin{proof}
	Fix $h\in\mathcal I$ and put $\lambda_h:=P_{\mathfrak H_h}\pi_h\kappa$.  We apply Theorem~\ref{thm: constrained strang} in the following case.  Take $V^{k-1}\times V^k$ with its graph inner products in the role of $V$, and $V_h^{k-1}\times V_h^k$ in the role of $V_h$.  The constraint spaces are $Q=\{0\}\times\mathfrak H^k$ and $Q_h=\{0\}\times\mathfrak H_h^k$.  The form
	\[
	\big((\eta,\alpha),(\theta,\beta)\big)\mapsto\langle\eta,\theta\rangle-\langle\alpha,d\theta\rangle+\langle d\eta,\beta\rangle+\langle d\alpha,d\beta\rangle
		\]
	plays the role of $B$, and the same expression in $d_h$ and $\langle\cdot,\cdot\rangle_h$ that of $B_h$.  The functionals are $(\theta,\beta)\mapsto F(\beta)$ and $(\theta_h,\beta_h)\mapsto F_h(\beta_h)$.  The interpolator and the reconstructor act componentwise by $\pi_h$ and $i_h$, so their norms are at most $C_\pi$ and $C_i$ [on the constraint space $\mathfrak H_h^k$, whose elements are closed, the $H_h$ norm is the graph norm].  Identify $(0,\kappa)$ with $\kappa$.  Then $\mathcal A$ is $\mathcal B$ and $\mathcal A_h$ is $\mathcal B_h$, because in the multiplier and constraint terms one argument is harmonic, hence closed, as in the proof of Theorem~\ref{thm: Hodge Dirac abstract approximation theorem}.  So the two constrained problems are \eqref{eq: Hodge Laplace problem} on the two complexes.  By Lemma~\ref{lem: discrete Hodge Laplace stability} the form $\mathcal B_h$ is inf-sup stable with a constant depending only on $c_H$, $C_i$ and $C_\pi$, and the discrete problem has a unique solution $(\eta_h,\alpha_h,\kappa_h)$.

	Since $P_{Q_h}$ is the restriction of the orthogonal projection onto $\{0\}\times\mathfrak H_h^k$, we have $P_{Q_h}\pi_h(0,\kappa)=(0,\lambda_h)$.  The left-hand side of \eqref{eq: Hodge Laplace bound} is at most $\sqrt2$ times that of \eqref{eq: constrained strang bound}, and the four terms on the right of \eqref{eq: constrained strang bound} are bounded as follows.  The approximation error of $(\eta,\alpha)$ is at most $\|\eta-i_h\pi_h\eta\|_V+\|\alpha-i_h\pi_h\alpha\|_V$.  The term $\|\kappa-i_h\lambda_h\|$ is bounded by \eqref{eq: harmonic reconstruction bound}.  The consistency error is bounded by \eqref{eq: HL consistency}.  The constraint term is $\sup_{q_h}|\langle\pi_h\alpha,q_h\rangle_h|/\|q_h\|_h$, which \eqref{eq:constraint_pairing_bound} bounds, since $\alpha$ is orthogonal to $\mathfrak H^k$ by the first equation of \eqref{eq: Hodge Laplace problem}.  Finally \eqref{eq: HL a priori} bounds $\|\eta\|_V+\|\alpha\|_V+\|\kappa\|$ by a multiple of $\|F\|_{(V^k)^*}$, so the terms carrying $\metricdefect$ alone are at most a multiple of $\metricdefect\|F\|_{(V^k)^*}$, and \eqref{eq:quasi_V} bounds the two interpolation defects by the infima in \eqref{eq: Hodge Laplace bound}.
\end{proof}

\begin{remark}\label{rem: Hodge Laplace bound comparison}
	Theorem~\ref{thm:Hodge Laplace bound} recovers \cite[Corollary~3.11]{Holst2012}, with a sharper metric term; see Section~\ref{sec: related abstract}.  Our constant
	depends only on $c_H$, $C_i$ and $C_\pi$.

	Holst and Stern factor their estimate through the intermediate problem \cite[(6)]{Holst2012}, the discrete problem with $i_h^*i_h$ inserted, which is the problem on the subcomplex $i_h(V_h^\bullet)$ with the inner product of $H^\bullet$. We directly apply Theorem~\ref{thm: constrained strang}.
\end{remark}

\subsubsection{Regularity}\label{sec: abstract regularity}
The error bounds of Theorems~\ref{thm: Hodge Dirac abstract approximation theorem} and~\ref{thm:Hodge Laplace bound} contain best approximation errors of the solution in the graph norm, which involves its differential. To turn them into rates we need a smaller space whose elements are approximated to the order of a scale parameter, the mesh size in our application. For a single Hilbert space this approach goes back to Aubin \cite{Aubin1967memoire}.
\begin{definition}\label{def: regularity scale}
	Let $(H^\bullet,V^\bullet,d)$ be a closed Hilbert complex and, for each
	$h\in\mathcal I$, let $(V_h^\bullet,d_h)$ be a finite-dimensional complex with
	chain maps $i_h\colon V_h^\bullet\to V^\bullet$ and
	$\pi_h\colon V^\bullet\to V_h^\bullet$.  A \emph{regularity scale}
	$(\mathcal W^\bullet,\|\cdot\|_{\mathcal W})$ for
	$\{i_h,\pi_h\}_{h\in\mathcal I}$, with \emph{scale parameter}
	$h\mapsto\ell_h\in(0,\infty)$ and constant $C_{\mathrm{app}}$, consists of a subspace
	$\mathcal W^k\subseteq V^k$ with a norm $\|\cdot\|_{\mathcal W}$ for each $k$,
	such that the interpolation defect satisfies
	\begin{equation}\label{eq:abstract_approx_property}\|\omega - i_h \pi_h \omega\|_H \le C_{\mathrm{app}}\, \ell_h \|\omega\|_{\mathcal{W}} \quad \forall \omega \in \mathcal{W}^k,\ h\in\mathcal I.\end{equation}
\end{definition}
The definition involves the maps $i_h$ and $\pi_h$ and the norms of the
continuous complex, and no inner product on the discrete complexes, so a
regularity scale can be verified without knowing the discrete inner products.

We will use the following consequence for the best approximation error of having a regularity scale:
\begin{prop}\label{prop:abstract_best_approx}
	Let $(\mathcal{W}^\bullet, \|\cdot\|_{\mathcal{W}})$ be a regularity scale for chain maps $\{i_h,\pi_h\}_{h\in\mathcal I}$ as in Definition~\ref{def: regularity scale}. If $\omega \in \mathcal{W}^k$ and $d\omega \in \mathcal{W}^{k+1}$, then
	\begin{equation}\label{eq:abstract_graph_bound}
		\inf_{\beta_h \in V_h^k} \|\omega - i_h \beta_h\|_V \le \|\omega - i_h \pi_h \omega\|_V \le C_{\mathrm{app}}\, \ell_h \left( \|\omega\|_{\mathcal{W}} + \|d\omega\|_{\mathcal{W}} \right).
	\end{equation}
\end{prop}

\begin{proof}
	The first inequality holds because $\beta_h=\pi_h\omega$ is admissible.  For
	the second, the graph norm of $\omega-i_h\pi_h\omega$ is at most
	$\|\omega-i_h\pi_h\omega\|+\|d\omega-i_h\pi_hd\omega\|$, because $i_h$ and
	$\pi_h$ are chain maps, and \eqref{eq:abstract_approx_property}, applied to
	$\omega\in\mathcal W^k$ and to $d\omega\in\mathcal W^{k+1}$, bounds the two terms.
\end{proof}

To apply Proposition~\ref{prop:abstract_best_approx} to the solutions of the
Hodge--Dirac and Hodge--Laplace problems, we need these solutions and their
differentials to lie in the scale. Lemma~\ref{lem: abstract mixed regularity}
gives a criterion. It asks that every element in the domains of both $d$ and
$d^*$ lie in the scale, with a bound. For the de Rham complex and the Sobolev
space $H^1$ this bound is Gaffney's inequality \cite{Gaffney1955}.

\begin{definition}\label{def: gaffney estimate}
	Let $(H^\bullet,V^\bullet,d)$ be a closed Hilbert complex, and for each
	$k$ let $\mathcal W^k\subseteq V^k$ be a subspace with a norm
	$\|\cdot\|_{\mathcal W}$, as in Definition~\ref{def: regularity scale}.
	The spaces $(\mathcal W^\bullet,\|\cdot\|_{\mathcal W})$ satisfy a \emph{Gaffney estimate} if there is a constant $C_G$ such that for every $k$ every $\omega\in V^k\cap V_k$ lies in
	$\mathcal W^k$ and satisfies
	\begin{equation}\label{eq: abstract gaffney}
		\|\omega\|_{\mathcal W}\le C_G\big(\|\omega\|+\|d\omega\|+\|d^*\omega\|\big).
	\end{equation}
\end{definition}

We norm $\mathcal W^\bullet=\bigoplus_k\mathcal W^k$ by
$\|\omega\|_{\mathcal W}^2:=\sum_k\|\omega^k\|_{\mathcal W}^2$, where $\omega^k$
is the component of $\omega$ in degree $k$. If \eqref{eq: abstract gaffney} holds for every component $\omega^k$, then it
holds for $\omega$ with $C_G$ multiplied by $\sqrt3$. The following lemma, which
expands on \cite[(7.29)]{Arnold2006}, places the solutions of both problems in
the scale, and, when $\gamma\in V_\bullet$, also $d\alpha$ for the
Hodge--Dirac problem and $d\eta$ for the Hodge--Laplace problem.

\begin{lemma}\label{lem: abstract mixed regularity}
	Let $(H^\bullet,V^\bullet,d)$ be a closed Hilbert complex, and let
	subspaces $(\mathcal W^\bullet,\|\cdot\|_{\mathcal W})$ as in
	Definition~\ref{def: gaffney estimate} satisfy a Gaffney estimate.

	Let $\gamma\in H^\bullet$ and let $(\alpha,\kappa)$ solve the Hodge--Dirac
	problem \eqref{eq: Hodge Dirac problem} with source $\gamma$.  Then
	$\alpha\in\mathcal W^\bullet$ with $\|\alpha\|_{\mathcal W}\le C\|\gamma\|$.
	Furthermore, if $\gamma\in V_\bullet$, then $d\alpha\in\mathcal W^\bullet$
	with $\|d\alpha\|_{\mathcal W}\le C(\|\gamma\|+\|d^*\gamma\|)$.

	Let $\gamma\in H^k$ and let $(\eta,\alpha,\kappa)$ solve the mixed
	Hodge--Laplace problem \eqref{eq: Hodge Laplace problem} with source
	$\gamma$.  Then $\eta=d^*\alpha$, and $\eta$, $\alpha$ and $d\alpha$ lie in
	$\mathcal W^\bullet$, each with norm at most $C\|\gamma\|$.  Furthermore,
	if $\gamma\in V_k$, then $d\eta\in\mathcal W^k$ with
	$\|d\eta\|_{\mathcal W}\le C(\|\gamma\|+\|d^*\gamma\|)$.

	Harmonic elements lie in $\mathcal W^\bullet$, with
	$\|q\|_{\mathcal W}\le C_G\|q\|$ for $q\in\mathfrak H^\bullet$.  The
	constant $C$ depends only on $C_G$ and on the Poincar\'e constant $c_H$.
\end{lemma}

\begin{proof}
	A harmonic element $q$ lies in $V^\bullet\cap V_\bullet$ with $dq=0$ and
	$d^*q=0$, so using \eqref{eq: abstract gaffney} gives
	$\|q\|_{\mathcal W}\le C_G\|q\|$.

	For the other assertions, Lemmas~\ref{lem: HD solutions in intersection} and~\ref{lem: HL solutions in intersection} show that each element $\omega$ concerned lies in $V^\bullet\cap V_\bullet$, and give $d^*\omega$. So \eqref{eq: abstract gaffney} applies to $\omega$. By the a priori bounds \eqref{eq: HD a priori} and \eqref{eq: HL a priori}, the norms $\|\omega\|$, $\|d\omega\|$ and $\|d^*\omega\|$ are at most a multiple of $\|\gamma\|$, or of $\|\gamma\|+\|d^*\gamma\|$ in the two ``furthermore'' assertions, where $d^*\omega=d^*\gamma$.
\end{proof}

This lemma immediately gives a rate for the harmonic gap
\eqref{eq: harmonic gap}. For subcomplexes of piecewise polynomial forms
containing all polynomials of degree at most $r$, Arnold, Falk and Winther
state that $\mu_h$ is at most a multiple of the $(r+1)$-st power of the mesh size \cite[Section~3.5]{arnold2010finite}.
\begin{cor}\label{cor: abstract harmonic gap rate}
	Let $(\mathcal W^\bullet,\|\cdot\|_{\mathcal W})$ be a regularity scale for
	chain maps $\{i_h,\pi_h\}_{h\in\mathcal I}$ as in
	Definition~\ref{def: regularity scale}, and suppose that it satisfies a
	Gaffney estimate. Then the harmonic gap \eqref{eq: harmonic gap} satisfies
	\begin{equation}\label{eq: abstract harmonic gap rate}
		\mu_h\le C_{\mathrm{app}}C_G\,\ell_h .
	\end{equation}
\end{cor}

\begin{proof}
		Let $q\in\mathfrak H^\bullet$ with $\|q\|=1$.  By
	Lemma~\ref{lem: abstract mixed regularity} each component $q^k$ lies in
	$\mathcal W^k$ with $\|q^k\|_{\mathcal W}\le C_G\|q^k\|$.  By
	\eqref{eq:abstract_approx_property},
	$\|q^k-i_h\pi_hq^k\|\le C_{\mathrm{app}}C_G\,\ell_h\|q^k\|$.  Squaring and summing over $k$ gives $\|q-i_h\pi_hq\|\le C_{\mathrm{app}}C_G\,\ell_h$.  Taking the supremum over $q$
	gives \eqref{eq: abstract harmonic gap rate}.
\end{proof}

Combining Proposition~\ref{prop:abstract_best_approx} and
Lemma~\ref{lem: abstract mixed regularity} gives rates for the best
approximation errors of the solutions of the Hodge--Dirac and Hodge--Laplace
problems.

\begin{cor}\label{cor: abstract best approximation rates}
	Let $(\mathcal W^\bullet,\|\cdot\|_{\mathcal W})$ be a regularity scale for
	$\{i_h,\pi_h\}_{h\in\mathcal I}$ that satisfies a Gaffney estimate. If $(\alpha,\kappa)$ solves the Hodge--Dirac
	problem \eqref{eq: Hodge Dirac problem} with source $\gamma\in V_\bullet$, then
	\begin{equation}\label{eq: abstract HD best approximation}
		\inf_{\beta_h\in V_h^\bullet}\|\alpha-i_h\beta_h\|_V
		\le C\,\ell_h\big(\|\gamma\|+\|d^*\gamma\|\big).
	\end{equation}
	If $(\eta,\alpha,\kappa)$ solves the Hodge--Laplace problem \eqref{eq: Hodge Laplace problem} with source $\gamma\in V_k$, then
	\begin{equation}\label{eq: abstract HL best approximation}
		\inf_{\theta_h\in V_h^{k-1}}\|\eta-i_h\theta_h\|_V
		+\inf_{\beta_h\in V_h^k}\|\alpha-i_h\beta_h\|_V
		\le C\,\ell_h\big(\|\gamma\|+\|d^*\gamma\|\big).
	\end{equation}
	In both, $C$ depends only on $C_{\mathrm{app}}$,
	$C_G$ and $c_H$.
\end{cor}

\begin{proof}
Consider the Hodge--Dirac problem. Both $\alpha$ and $d\alpha$ lie in $\mathcal W^\bullet$, by Lemma~\ref{lem: abstract mixed regularity}, with norms at most a multiple of $\|\gamma\|+\|d^*\gamma\|$. The infimum in \eqref{eq: abstract HD best approximation} is at most the value at $\beta_h=\pi_h\alpha$. Proposition~\ref{prop:abstract_best_approx} bounds $\|\alpha^k-i_h\pi_h\alpha^k\|_V$ by
	$C_{\mathrm{app}}\ell_h(\|\alpha^k\|_{\mathcal W}+\|(d\alpha)^{k+1}\|_{\mathcal W})$. Squaring and summing over $k$ bounds $\|\alpha-i_h\pi_h\alpha\|_V$ by $\sqrt2\,C_{\mathrm{app}}\ell_h(\|\alpha\|_{\mathcal W}+\|d\alpha\|_{\mathcal W})$.

	The proof of the statement for the Hodge--Laplace problem is similar, with Proposition~\ref{prop:abstract_best_approx} applied once to $\alpha$ and once to $\eta$.
\end{proof}

For the rates of Section~\ref{sec: abstract rates}, the two harmonic terms of
the error bounds need only be bounded by $\mu_h$ times the norm of the source,
and for solutions of the two problems this is immediate:
\begin{cor}\label{cor: harmonic error term coarsening}
	If $(\alpha,\kappa)$ solves the
	Hodge--Dirac problem with source $\gamma\in H^\bullet$, then
	\begin{equation}\label{eq: abstract harmonic terms}
		\inf_{q_h\in V_h^\bullet}\|\kappa-i_hq_h\|_V
		+\mu_h\inf_{v_h\in V_h^\bullet}\|P_B\alpha-i_hv_h\|_V
		\le C\,\mu_h\|\gamma\| .
	\end{equation}
	The same holds for the solution $(\eta,\alpha,\kappa)$ of the mixed
	Hodge--Laplace problem with source $\gamma\in H^k$, with the infima taken
	over $V_h^k$.  The constant $C$ depends only on $c_H$.
\end{cor}

\begin{proof}
	The element $\kappa$ is harmonic, so $\kappa-i_h\pi_h\kappa$ is closed and $\inf_{q_h}\|\kappa-i_hq_h\|_V\le\|\kappa-i_h\pi_h\kappa\|\le\mu_h\|\kappa\|$, by the definition \eqref{eq: harmonic gap} of $\mu_h$.  Taking $v_h=0$ gives
	$\inf_{v_h}\|P_B\alpha-i_hv_h\|_V\le\|P_B\alpha\|_V=\|P_B\alpha\|\le\|\alpha\|$, as
	$P_B\alpha$ is closed.  By \eqref{eq: HD a priori}, respectively
	\eqref{eq: HL a priori}, $\|\alpha\|+\|\kappa\|$ is at most a multiple of
	$\|\gamma\|$.
\end{proof}

\subsubsection{Rates}\label{sec: abstract rates}

We take the discrete source $F\circ i_h$, for which the source error in the bounds of Theorems~\ref{thm: Hodge Dirac abstract approximation theorem} and~\ref{thm:Hodge Laplace bound} vanishes.  For a square-integrable source $\gamma$ this is $\langle i_h^*\gamma,\cdot\rangle_h$.  The remaining terms are of two kinds.  The best approximation errors and the harmonic terms depend on the solution, and the metric defect depends on the discrete inner products.  As shown in Section~\ref{sec: abstract regularity}, regularity of the
solutions gives rates for the terms of the first kind. We further add the hypothesis that the
source lies in the domain of $d^*$, which places the differential of the solution in
the scale. Indeed, that differential, $d\alpha$ for the Hodge--Dirac problem and
$d\eta$ for the Hodge--Laplace problem, is the exact part of the source. It is
closed, and its codifferential is $d^*\gamma$ by
Lemmas~\ref{lem: HD solutions in intersection}
and~\ref{lem: HL solutions in intersection}, so the Gaffney estimate applies to
it. This gives the following theorem.

\begin{theorem}\label{thm: abstract rates}
	Let $(H^\bullet,V^\bullet,d)$ be a closed Hilbert complex with an
	approximation $\{i_h,\pi_h\}_{h\in\mathcal I}$.  Let
	$(\mathcal W^\bullet,\|\cdot\|_{\mathcal W})$ be a regularity scale for
	$\{i_h,\pi_h\}_{h\in\mathcal I}$ with scale parameter $\ell_h$, and suppose
	that it satisfies a Gaffney estimate.

	Let $\gamma\in V_\bullet$, and let
	$(\alpha,\kappa)$ and $(\alpha_h,\kappa_h)$ solve the Hodge--Dirac problems
	on $(H^\bullet,V^\bullet,d)$ and on $(H_h^\bullet,V_h^\bullet,d_h)$ with
	sources $\gamma$ and $i_h^*\gamma$.  Then
	\begin{equation}\label{eq: abstract HD rate}
		\|\alpha-i_h\alpha_h\|_V+\|\kappa-i_h\kappa_h\|
		\le C\big(\ell_h+\metricdefect\big)\big(\|\gamma\|+\|d^*\gamma\|\big).
	\end{equation}

	Let $\gamma\in V_k$, let $(\eta,\alpha,\kappa)$ and
	$(\eta_h,\alpha_h,\kappa_h)$ solve the Hodge--Laplace problems in degree $k$
	with sources $\gamma$ and $i_h^*\gamma$, and write $\mathrm{err}_h$ for the
	left-hand side of \eqref{eq: Hodge Laplace bound}.  Then
	\begin{equation}\label{eq: abstract HL rate}
		\mathrm{err}_h\le C\big(\ell_h+\metricdefect\big)\big(\|\gamma\|+\|d^*\gamma\|\big).
	\end{equation}
	The constant $C$ depends only on $c_H$, $C_i$, $C_\pi$, $C_{\mathrm{app}}$ and
	$C_G$.
\end{theorem}

\begin{proof}
	Theorem~\ref{thm: Hodge Dirac abstract approximation theorem} gives
	\eqref{eq: Hodge Dirac bound}, whose source error vanishes and whose
	metric defect term is at most $\metricdefect\|\gamma\|$, as $\langle\gamma,\cdot\rangle$ has norm at
	most $\|\gamma\|$.  Its second and third terms are at most
	$C\mu_h\|\gamma\|$ by Corollary~\ref{cor: harmonic error term coarsening}, and
	$\mu_h\le C_{\mathrm{app}}C_G\,\ell_h$ by
	Corollary~\ref{cor: abstract harmonic gap rate}.  Its first term is at most
	$C\ell_h(\|\gamma\|+\|d^*\gamma\|)$ by
	Corollary~\ref{cor: abstract best approximation rates}.  This gives
	\eqref{eq: abstract HD rate}.  The Hodge--Laplace case is the same, with
	Theorem~\ref{thm:Hodge Laplace bound} in place of
	Theorem~\ref{thm: Hodge Dirac abstract approximation theorem}.
\end{proof}

The right-hand sides of the estimates in Theorem~\ref{thm: abstract rates} no longer depend on the
solution. The metric defect is the one term that has to be estimated for the approximation at hand. If
$\gamma\notin V_\bullet$, Corollary~\ref{cor: abstract harmonic gap rate} still
bounds $\mu_h$ in \eqref{eq: Hodge Dirac bound} and
\eqref{eq: Hodge Laplace bound} by $C_{\mathrm{app}}C_G\,\ell_h$, but the best
approximation error of the solution need not have a rate.

The error still tends to zero, without a rate, for every source
$F\in(V^\bullet)^*$ as the scale parameter and the metric defect tend to zero.
It is at most a multiple of $\|F\|$, uniformly in $h$, and it tends to zero on
the dense set of sources of Theorem~\ref{thm: abstract rates}.

\begin{cor}\label{cor: abstract L2 convergence}\label{cor: abstract dual convergence}
	Let $(H^\bullet,V^\bullet,d)$, $\{i_h,\pi_h\}_{h\in\mathcal I}$ and
	$(\mathcal W^\bullet,\|\cdot\|_{\mathcal W})$ be as in
	Theorem~\ref{thm: abstract rates}, and let $F\in(V^\bullet)^*$.  Let
	$(\alpha,\kappa)$ solve the Hodge--Dirac problem
	\eqref{eq: Hodge Dirac problem} with source $F$ and, for
	$h\in\mathcal I$, let $(\alpha_h,\kappa_h)$ solve the discrete
	Hodge--Dirac problem on $(H_h^\bullet,V_h^\bullet,d_h)$ with source
	$F\circ i_h$.  Then for every $\varepsilon>0$ there is $\delta>0$ such that
	\[
		\|\alpha-i_h\alpha_h\|_V+\|\kappa-i_h\kappa_h\|\le\varepsilon
		\qquad\text{for all $h\in\mathcal I$ with }
		\ell_h+\metricdefect\le\delta .
	\]
	The same holds for the mixed Hodge--Laplace problem
	\eqref{eq: Hodge Laplace problem} in degree $k$, with
	$F\in(V^k)^*$ and the left-hand side of \eqref{eq: Hodge Laplace bound} in
	place of this error.
\end{cor}

\begin{proof}
We give the proof for the Hodge--Dirac problem.  For the Hodge--Laplace
problem replace the Hodge--Dirac results by their Hodge--Laplace
counterparts, $V_\bullet$ by $V_k$, and the target space below by
$V^{k-1}\times V^k\times H^k$.

Let $(h_n)$ be a sequence in $\mathcal I$ with
$\ell_{h_n}+\metricdefect[h_n]\to0$.  For $F\in(V^\bullet)^*$ let
$T_nF:=(\alpha-i_{h_n}\alpha_{h_n},\kappa-i_{h_n}\kappa_{h_n})\in
V^\bullet\times H^\bullet$, with $(\alpha,\kappa)$ and
$(\alpha_{h_n},\kappa_{h_n})$ the solutions with sources $F$ and
$F\circ i_{h_n}$.  The problems are linear and their solutions unique, so
$T_n$ is linear.  It is bounded uniformly in $n$.  Indeed,
$\|\alpha\|_V+\|\kappa\|$ is at most a multiple of $\|F\|_{(V^\bullet)^*}$ by \eqref{eq: HD a priori}.  The functional
$(v_h,q_h)\mapsto F(i_hv_h)$ has norm at most $C_i\|F\|_{(V^\bullet)^*}$, so
Theorem~\ref{thm: babuska} gives
$\|\alpha_h\|_{V_h}+\|\kappa_h\|_h\le\sqrt2\,\gamma_0^{-1}C_i\|F\|_{(V^\bullet)^*}$.
Finally $\|i_h\alpha_h\|_V\le C_i\|\alpha_h\|_{V_h}$, and
$\|i_h\kappa_h\|\le C_i\|\kappa_h\|_h$ because discrete harmonic forms are
closed.

For $u\in V_\bullet$ we have $T_n\langle u,\cdot\rangle\to0$, by
Theorem~\ref{thm: abstract rates} [with source $u$; the discrete source
$\langle u,\cdot\rangle\circ i_{h_n}$ of $T_n$ is the pairing with
$i_{h_n}^*u$].  The functionals $\langle u,\cdot\rangle$, $u\in V_\bullet$, form a dense subspace of
$(V^\bullet)^*$.  Indeed, $V_\bullet$ is dense in $H^\bullet$, because $d^*$,
the adjoint of the closed densely defined operator $d$, is densely defined
\cite[Section~3.1.3]{arnold2010finite}; the map $\gamma\mapsto\langle\gamma,\cdot\rangle$ is
bounded; and its image is dense in $(V^\bullet)^*$
\cite[Chapter~1, (2.39)]{LionsMagenes1972} [with $V^\bullet$ in the role of
$V$ and $H^\bullet$ in the role of $H$; $V^\bullet$ is dense in $H^\bullet$
with continuous injection, and the injection of $H$ into $V'$ there is
$\gamma\mapsto\langle\gamma,\cdot\rangle$].

So \cite[Proposition~1.19]{vanNeerven2022} applies.
It shows that $T_nF$ converges for every $F$, to $TF$ for a bounded linear
operator $T$.  The operator $T$ vanishes on the dense subspace, hence
everywhere, and so $T_nF\to0$ for every $F$.  As the sequence was arbitrary,
the assertion follows.
\end{proof}

\section{Finite element exterior calculus on closed manifolds}\label{sec: feec on manifolds}Throughout this section we fix a closed oriented $n$-dimensional Riemannian manifold $(M,g)$, that is, compact and without boundary; we may assume $M$ connected, treating its components separately. We discuss the approximation of the Hilbert complex of square-integrable forms on $M$ by finite element exterior calculus (FEEC) complexes on smooth triangulations $K$ of $M$. Section~\ref{sec: triangulations} introduces the exact notion of triangulation used here, Section~\ref{sec: existence of triangulations} provides them, and we review the FEEC complexes in Section~\ref{sec:FEECforms}. The approximation is indexed by a family of triangulations, and its scale parameter is the mesh size $\ell_K$ of $K$, the largest diameter of a simplex of $K$ (Section~\ref{sec: triangulations}).

We prove convergence by applying Theorem~\ref{thm: abstract rates}, which requires checking three hypotheses. The finite element complexes must form an approximation in the sense of
Definition~\ref{def: approximation of hilbert complex},
and we also need its metric defect to be of order $\ell_K$
(Theorem~\ref{thm: feec approximation} below shows both). The space
$H^1\Lambda^\bullet(M)$ of forms with one Sobolev derivative
(Section~\ref{sec: l2 forms}) must be a regularity scale for its maps, with
scale parameter $\ell_K$; this follows from the approximation estimate of
Section~\ref{sec: feec approximation estimates}
(Corollary~\ref{cor: feec regularity scale}). And this scale must satisfy a
Gaffney estimate, which follows from elliptic regularity for $d+d^*$ and
concerns $M$ alone (Lemma~\ref{lem: feec dirac regularity}); by
Lemma~\ref{lem: abstract mixed regularity} it gives the solutions of both
problems one Sobolev derivative, bounded by the source. The first two hypotheses depend on
the family of triangulations and can fail for poorly chosen families; we prove
them for the \emph{uniformly regular} families of
Section~\ref{sec: uniformly regular families}. The convergence theorems are in Section~\ref{sec: convergence on manifolds}.

\subsection{Triangulations of closed Riemannian manifolds}\label{sec: triangulations}

This section introduces the triangulations we work with and the uniformly
regular families along which our estimates hold, and the piecewise flat metric
they carry, which is within order $\ell_K$ of $g$ along such families.

We work with the smooth triangulations of \cite[Section~3]{Whitehead1940}, in
the intrinsic formulation of \cite[Section~4.1]{Christiansen2007} and
\cite[Section~3.1]{Licht2023smoothed}. Write $\Delta^d\subset\R^d$ for the
reference $d$-simplex, the convex hull of the origin and the unit vectors $e_1,\dots,e_d$, with its vertices ordered $0,e_1,\dots,e_d$. A
\emph{smooth $d$-simplex} in $M$ is a smooth embedding
$\imath_T\colon\Delta^d\to M$, with image $T$. A smooth $l$-simplex $S$ is a
\emph{smooth subsimplex} of $T$ if $\imath_S$ is the restriction of $\imath_T$
to an $l$-dimensional face of $\Delta^d$, that face being identified with
$\Delta^l$ by the affine map that preserves the order of the vertices.

\begin{definition}\label{def: smooth triangulation}
	A \emph{smooth triangulation} of $M$ is a locally finite collection $K$ of smooth
	simplices such that every point of $M$ lies in an $n$-simplex of $K$, every
	smooth subsimplex of a simplex of $K$ is in $K$, and any two simplices of
	$K$ are disjoint or meet in a common smooth subsimplex.
\end{definition}
The $n$-simplices of $K$ cover $M$, and two distinct ones meet only in
simplices of lower dimension, whose union is a null set.  So a metric or a
form given on each $n$-simplex defines one almost everywhere on $M$.
As $M$ is compact, $K$ is finite.  We write $h_T$ for the $g$-diameter of a
simplex $T$ and $\ell_K=\max_{T}h_T$ for the largest of these.

\subsubsection{Uniformly regular families}
\label{sec: uniformly regular families}%
The estimates for the finite element approximation are proved one simplex at a
time: we pull $k$-forms back along the parametrisation $\imath_T$, estimate them
on the reference simplex, and push them forward to $T$. Pulling back and pushing
forward introduce at most factors of $|d\imath_T|^k$ and $|d\imath_T^{-1}|^k$,
where $|d\imath_T|$ is the operator norm of the differential of $\imath_T$, for
the Euclidean metric on the reference simplex and the metric $g$ on $M$. This
means that we need control over the norm of the differential at the scale of the
diameter $h_T$ of $T$ and, for the metric error and the approximation estimate,
over its higher covariant derivatives. We also need control over the number of
simplices that meet, and over how the diameter of a simplex compares with those
of its faces. The next definition, which follows
\cite[Section~3.2]{Licht2023smoothed}, makes this precise.

\begin{definition}\label{def: uniformly regular}
	Let $K$ be a smooth triangulation of $M$.
	\begin{itemize}
		\item The \emph{shape constant} $C_\mu$ of $K$ is the smallest number
		such that every simplex $T$ of $K$ of positive dimension $d$ satisfies
		$\sup_{\Delta^d}|d\imath_T|\le C_\mu h_T$ and
		$\sup_T|d\imath_T^{-1}|\le C_\mu h_T^{-1}$.
		\item For $m\ge2$, the $m$-th \emph{higher shape constant} $C_{\mu,m}$ of
		$K$ is the smallest number such that every simplex $T$ of $K$ of
		positive dimension $d$ satisfies
		$\sup_{\Delta^d}|\nabla^{m-1}d\imath_T|\le C_{\mu,m}h_T^m$, where $\nabla$
		is the flat connection on $\Delta^d$ combined with the pullback of the
		Levi-Civita connection.
		\item The \emph{adjacency bound} $C_\sharp$ of $K$ is the largest number
		of simplices of $K$ that intersect a single simplex of $K$.
		\item The \emph{face ratio} $C_\triangle$ of $K$ is the smallest number
		such that $h_T\le C_\triangle h_S$ whenever $S$ is a subsimplex of a
		simplex $T$ of $K$, where, when $S$ is a vertex, $h_S$ stands for the
		smallest $h_{T'}$ among the simplices $T'$ of positive dimension that
		contain $S$.
	\end{itemize}
	A family $\mathcal K$ of smooth triangulations of $M$ with
	$\inf_{K\in\mathcal K}\ell_K=0$ is called \emph{uniformly regular} if each of
	the constants $C_\mu$, $C_{\mu,m}$ for all $m\ge2$, $C_\sharp$ and
	$C_\triangle$ is bounded independently of $K\in\mathcal K$.
\end{definition}

The shape constant says that $\imath_T$ distorts lengths by a factor comparable
to $h_T$ in both directions.  So $T$ is a copy of the reference simplex scaled
by its diameter $h_T$, up to a distortion bounded by $C_\mu$.  The constant $C_{\mu,2}$ bounds the second
derivative of $\imath_T$ by $h_T^2$, so that the pulled-back metric
$\imath_T^*g$ varies over $\Delta^d$ by a factor $1+O(h_T)$.  This is what makes
the metric error of the order of the mesh size
(Proposition~\ref{prop: piecewise flat metric error}).
Our estimates use the constants $C_{\mu,m}$ only for $m\le2$.  Those with
$m\ge3$ enter approximation estimates of higher order, in which the
derivatives of the pullback of a form along $\imath_T$ involve the higher
derivatives of $\imath_T$ \cite[Theorem~7.6]{Licht2023smoothed}.

The last two conditions in Definition~\ref{def: uniformly regular} concern how the simplices fit together.  The adjacency bound says that a simplex has boundedly many neighbours.  The face ratio $C_\triangle$ says that a simplex and its faces have comparable diameters, and hence so do neighbouring simplices.  They enter through Licht's smoothed projection and his approximation estimate (Proposition~\ref{prop: feec smoothed projection} and Theorem~\ref{thm: feec l2 approximation}), both of which combine information from the neighbours of a simplex. Over a domain in Euclidean space the shape constant bounds $C_\sharp$ and $C_\triangle$.  On a	manifold it need not \cite[Remark~3.3]{Licht2023smoothed}, so they are conditions of their own.

\subsubsection{The piecewise flat metric}
The finite element complexes will carry the inner products of the following
metric, which makes each simplex flat, with the edge lengths it has in $(M,g)$.

\begin{definition}\label{def: piecewise flat metric}
	Let $K$ be a smooth triangulation of $M$.  For an edge $e$ of $K$ let
	$\ell_e$ be the $g$-length of the curve $e\subset M$.  We call $K$
	\emph{realisable} if for every $n$-simplex $T$ of $K$ the numbers $\ell_e$,
	for $e$ an edge of $T$, are the edge lengths of a non-degenerate simplex
	$\sigma_T$ in Euclidean space $\R^n$.  The \emph{piecewise flat metric}
	$g_K$ of a realisable triangulation is then defined on each $n$-simplex $T$
	as follows.  Let $A_T\colon\Delta^n\to\sigma_T$ be the affine isomorphism
	that matches the vertices in the way the edge lengths prescribe.  Then
	$\imath_T^*g_K$ is the pullback along $A_T$ of the Euclidean metric.
	The \emph{metric error} is the smallest $\epsilon_K\ge0$ such that for a.e.\ $x\in M$ and all $a\in T_xM $ we have
	\begin{equation}\label{eq: feec metric error}
		(1-\epsilon_K)\,g(a,a)\le g_K(a,a)\le(1+\epsilon_K)\,g(a,a).
	\end{equation}
\end{definition}

The simplex $\sigma_T$ is unique up to an isometry of $\R^n$, so $g_K$ does
not depend on its choice.  The metric $g_K$ is constant on each $n$-simplex in
the coordinates of $\Delta^n$, gives every edge of $K$ its $g$-length, and
across a face only its restriction to that face is continuous.  For $n=2$ a triangulation is realisable exactly
when the edge lengths of every triangle satisfy the strict triangle
inequality.  For $n\ge3$ the triangle inequalities on the faces do not
suffice; the condition is that a Gram matrix formed from the squared edge
lengths is positive definite
\cite[Theorem~1]{Schoenberg1935}, \cite[Section~5.1]{DyerVegterWintraecken2015}.  Every sufficiently fine
triangulation in a uniformly regular family is realisable
(Proposition~\ref{prop: piecewise flat metric error}).

\begin{remark}
	Metrics that are Euclidean on each simplex of a simplicial decomposition were
	introduced by Regge \cite[Section~4]{Regge1961}; such a metric is determined by
	the lengths of its edges \cite[Section~1]{Cheeger1984}.
\end{remark}

The metric error $\epsilon_K$ of Definition~\ref{def: piecewise flat metric}
enters every error bound below, through the metric defect
(Proposition~\ref{prop: feec metric defect}), and it is defined only when $K$ is
realisable. The next proposition shows that a sufficiently fine triangulation is
realisable, with metric error of the order of its mesh size. How fine it must be
depends only on $n$ and the shape constants $C_\mu$ and $C_{\mu,2}$, so one threshold serves a
whole uniformly regular family.

\begin{prop}\label{prop: piecewise flat metric error}
	Let $K$ be a smooth triangulation of $M$.  There are constants $C$ and
	$\ell_0>0$, depending only on $n$ and on the shape constants $C_\mu$ and
	$C_{\mu,2}$ of $K$, such that if
	$\ell_K\le\ell_0$, then $K$ is realisable and
	\[
	\epsilon_K\le C\,\ell_K\le\tfrac12 .
	\]
\end{prop}

\begin{proof}
	Throughout the proof $C$ denotes a constant depending only on $n$, $C_\mu$ and
	$C_{\mu,2}$, whose value may change from one occurrence to the next.
	Fix an $n$-simplex $T$ of $K$, and write $h:=h_T\le\ell_K$.  The
	quantity $\epsilon_K$ is the supremum over the $n$-simplices of the corresponding
	quantity for $T$, and pulling both metrics back along $\imath_T$ does not change
	their comparison.  It therefore suffices to show that the edge lengths of $T$
	are those of a non-degenerate Euclidean simplex (for realisability), and to compare
	$\hat g:=\imath_T^*g$ with the constant metric $\imath_T^*g_K$ on $\Delta^n$ (for the bound).
Both go through one intermediate object: the value $G:=\hat g_{x_0}$ of
$\hat g$ at a fixed point $x_0\in\Delta^n$, regarded as a constant metric on
$\Delta^n$.

On one side, $\hat g$ is nearly constant, so it is close to $G$.  Measured
against the Euclidean metric of $\Delta^n$, the coefficients of $\hat g$ are
of size $h^2$, by the bound on $d\imath_T$: the definition of the shape
constant gives $|d\imath_T|\le C_\mu h$ and $|d\imath_T^{-1}|\le C_\mu h^{-1}$,
so for $x\in\Delta^n$ and $v\in\R^n$
\begin{equation}\label{eq: pf metric bounds}
    C_\mu^{-2}h^2\,|v|^2\le\hat g_x(v,v)\le C_\mu^2h^2\,|v|^2 .
\end{equation}
The coefficients of the derivative of $\hat g$ are of size $h^3$, by the bound
on $\nabla d\imath_T$ in terms of $C_{\mu,2}$.  To see this, let
$u,v,w\in\R^n$.  Then $\hat g(v,w)=g\big(d\imath_T(v),d\imath_T(w)\big)$ is a
function on $\Delta^n$, and we differentiate it in the direction $u$.  The
Levi-Civita connection is compatible with $g$, and so is its pullback
$\nabla$ along $\imath_T$.  Hence
\[
    \partial_u\big(\hat g(v,w)\big)
    =g\big(\nabla_u(d\imath_T(v)),d\imath_T(w)\big)
    +g\big(d\imath_T(v),\nabla_u(d\imath_T(w))\big) .
\]
The vector $v$ is constant, so $\nabla_u(d\imath_T(v))=(\nabla d\imath_T)(u,v)$, and
likewise for $w$.  The bound on $\nabla d\imath_T$ gives
$|(\nabla d\imath_T)(u,v)|_g\le C_{\mu,2}h^2|u||v|$, and the bound on $d\imath_T$ gives
$|d\imath_T(w)|_g\le C_\mu h|w|$.  So
\[
    \big|\partial_u\big(\hat g(v,w)\big)\big|\le Ch^3\,|u||v||w| .
\]
Consequently $\hat g_x$ differs from $G$ by a relative error of order $h$.
Indeed, integrating along the segment from $x_0$ to $x$, whose length is at
most $\operatorname{diam}(\Delta^n)$, we get
$|\hat g_x(v,w)-G(v,w)|\le Ch^3|v||w|$ for every $x\in\Delta^n$.  By the lower
bound in \eqref{eq: pf metric bounds} for $G$ we have $h^2|v|^2\le C\,G(v,v)$, so
$|\hat g_x(v,v)-G(v,v)|\le Ch\,G(v,v)$.  That is
\begin{equation}\label{eq: pf frozen}
    (1-Ch)\,G\le\hat g_x\le(1+Ch)\,G
    \qquad\text{on }\Delta^n .
\end{equation}

On the other side, the edge lengths of $T$ are nearly those of the Euclidean
simplex that $G$ describes, and a perturbation result for Euclidean simplices
turns this into realisability together with a comparison of
$\imath_T^*g_K$ with $G$.  For an edge $e$ of $\Delta^n$ write $v_e$ for its
edge vector, the difference of its endpoints in the order of the vertices,
so that $|v_e|_G$ is the $G$-length of $e$.  Its $\hat g$-length is $\ell_e$,
as $\imath_T$ maps $e$ onto the corresponding edge of $T$, and because $\hat g$
is close to $G$ the two lengths are close.  By \eqref{eq: pf frozen},
$\big||v|_{\hat g_x}-|v|_G\big|\le Ch\,|v|_G$ for $x\in\Delta^n$ and
$v\in\R^n$ [as $(1-Ch)^{1/2}\ge1-Ch$ and $(1+Ch)^{1/2}\le1+Ch$], and we
suppose $h$ so small that $Ch<1$.  By \cite[Lemma~6.13]{vonDeylen2015}
\[
    \big|\ell_e-|v_e|_G\big|\le Ch\,|v_e|_G .
\]
We now apply \cite[Proposition~3.16]{vonDeylen2015}, with the Euclidean
simplex that $G$ describes as its given simplex and the numbers $\ell_e$ as its
desired edge lengths.  The given simplex satisfies the hypothesis of that
proposition.  By
\eqref{eq: pf metric bounds} its edges have length
$|v_e|_G\le\sqrt2\,C_\mu h<2C_\mu h$, and its volume $(\det G)^{1/2}/n!$ is at
least $C_\mu^{-n}h^n/n!$, so it is $(\vartheta,2C_\mu h)$-small in the sense of
\cite[Definition~3.3]{vonDeylen2015} for a $\vartheta>0$ depending only on $n$
and $C_\mu$.  The lengths $\ell_e$ satisfy the hypothesis with
$\varepsilon=Ch$, and we suppose $h$ so small that $Ch\le\frac12$.  The
proposition then gives a Euclidean simplex whose edge lengths are the numbers
$\ell_e$ and whose metric $\bar g$ over $\Delta^n$ satisfies
$|(G-\bar g)(v,v)|\le Ch\,G(v,v)$.  In particular $\bar g$ is
positive definite, so this simplex is non-degenerate.  As this holds for
every $n$-simplex $T$, the triangulation $K$ is realisable.  The metric
$\bar g$ is the pullback of the Euclidean metric along the affine map onto
the simplex that matches the vertices, so it is $\imath_T^*g_K$ by
Definition~\ref{def: piecewise flat metric}.  That is,
\[
    (1-Ch)\,G\le\imath_T^*g_K\le(1+Ch)\,G .
\]

As $\hat g$ and $\imath_T^*g_K$ are both close to $G$, they are close to each
other.  Combining \eqref{eq: pf frozen} with the bound just proved, and using
$Ch\le\frac12$, gives
$(1-Ch)\,\hat g_x\le\imath_T^*g_K\le(1+Ch)\,\hat g_x$ on $\Delta^n$ for a larger
$C$, so the quantity for $T$ is at most $C\ell_K$.  Let $\ell_0$ be so small
that $C\ell_0\le\frac12$ for every constant $C$ in this proof.  If
$\ell_K\le\ell_0$, then the smallness assumptions on $h$ hold for every
$n$-simplex $T$, and taking the supremum over them gives the assertion.
\end{proof}

\begin{remark}\label{rem: piecewise flat sharper}
	The rate $O(\ell_K)$ is all that the error bounds need, since the
	interpolation defect is of the same order.  For constructions with a better rate, and why we do not use them, see Section~\ref{sec: related triangulations}.
	The canonical lowest-order Regge interpolant of $g$ \cite[(2.4b)]{Li2018} satisfies the same bound \cite[Theorem~2.5]{Li2018} as the piecewise flat metric;
	we use edge lengths because they depend only on the curves $e\subseteq M$
	and not on how the edges are parametrised.
\end{remark}

\subsection{Existence of triangulations}\label{sec: existence of triangulations}
Every convergence statement of this paper assumes a uniformly regular family of
smooth triangulations (Definition~\ref{def: uniformly regular}).
This section shows that every closed Riemannian manifold has such a family.
Section~\ref{sec: initial triangulation} fixes an initial triangulation, of the
kind the theorem of Cairns and Whitehead provides.
Section~\ref{sec: freudenthal subdivisions} subdivides it repeatedly by
Freudenthal's rule and checks, constant by constant, that the subdivisions are
uniformly regular
(Theorem~\ref{thm: freudenthal subdivisions uniformly regular}).

\subsubsection{The initial triangulation}\label{sec: initial triangulation}
The subdivision starts from one smooth triangulation $K_0$ of $M$, which exists
by the theorem of Cairns and Whitehead \cite[Section~2]{Cairns1935},
\cite[Theorem~7]{Whitehead1940}, \cite[Section~20]{MilnorStasheff1974}.  Any two
of its simplices are disjoint or meet in a single common subsimplex, so the vertex sets of its simplices form
a finite simplicial complex.  We also write $K_0$ for a realisation of this
complex in a Euclidean space, and order the vertices of each of its simplices
$\sigma$ as those of $\imath_\sigma$.  The parametrisations then glue to a
homeomorphism $\chi\colon|K_0|\to M$ that agrees on each closed simplex
$\sigma$ with $\imath_\sigma$ composed with the inverse of the order-preserving
affine map of $\Delta^d$ onto $\sigma$ \cite[Remark~3.1]{Licht2023smoothed}.
The restriction of $\chi$ to $\sigma$ is smooth and of maximal rank: it extends
to a smooth embedding of an open neighbourhood of $\sigma$ in its affine hull,
as $\imath_\sigma$ extends to one of an open neighbourhood of $\Delta^d$
\cite[Section~3.1]{Licht2023smoothed}.  So $(K_0,\chi)$ is a smooth
triangulation in the sense of Milnor and Stasheff
\cite[Section~20]{MilnorStasheff1974}, the smooth analogue of Whitehead's
$C^1$-triangulations \cite[Sections~2 and~3]{Whitehead1940}.

\subsubsection{Freudenthal subdivisions}\label{sec: freudenthal subdivisions}
Let $\tau$ be a $d$-simplex in a Euclidean space with ordered vertices
$q_0<\dots<q_d$.  Its \emph{edge vectors} are $v_i:=q_i-q_{i-1}$ for
$1\le i\le d$.  A \emph{conjugate} of $\tau$ is a simplex with ordered vertices
$q_0,\ q_0+v_{\varpi(1)},\ \dots,\ q_0+v_{\varpi(1)}+\dots+v_{\varpi(d)}$ for a permutation
$\varpi$ of $\{1,\dots,d\}$, that is, a simplex whose edge vectors are those of
$\tau$ in some order \cite[Section~1]{Freudenthal1942}.

For $0\le i\le j\le d$ put $q_{ij}:=\frac12(q_i+q_j)$, so that $q_{ii}=q_i$ and
$q_{ij}$ is the midpoint of the edge from $q_i$ to $q_j$ when $i<j$.  The
\emph{Freudenthal subdivision} of $\tau$ has the points $q_{ij}$ as vertices.  Its
$d$-simplices are spanned by the sequences of $d+1$ of these points that begin at
$q_{0u}$ for some $0\le u\le d$, and in which each $q_{ij}$ is followed either by
$q_{i+1,j}$, allowed when $i<u$, or by $q_{i,j+1}$, allowed when $j<d$.  In dimension two this
divides a triangle into four by the segments joining the midpoints of its edges.

\begin{figure}[ht]
\centering
\begin{tikzpicture}[scale=0.85,
    vertex/.style={circle,fill,inner sep=1.2pt},
    every node/.style={font=\footnotesize}]
  \coordinate (q0) at (0,0);
  \coordinate (q1) at (3,0);
  \coordinate (q2) at (3.8,2.4);
  \coordinate (q01) at (1.5,0);
  \coordinate (q02) at (1.9,1.2);
  \coordinate (q12) at (3.4,1.2);
  \fill[black!12] (q01) -- (q02) -- (q12) -- cycle;
  \draw (q0) -- (q1) -- (q2) -- cycle;
  \draw (q01) -- (q02) -- (q12) -- cycle;
  \foreach \p in {q0,q1,q2,q01,q02,q12} \node[vertex] at (\p) {};
  \node[below left] at (q0) {$q_{00}$};
  \node[below] at (q01) {$q_{01}$};
  \node[below right] at (q1) {$q_{11}$};
  \node[left] at (q02) {$q_{02}$};
  \node[right] at (q12) {$q_{12}$};
  \node[above] at (q2) {$q_{22}$};
\end{tikzpicture}\hspace{1.5em}%
\begin{tikzpicture}[scale=0.85,
    vertex/.style={circle,fill,inner sep=1.2pt},
    every node/.style={font=\footnotesize}]
  \coordinate (q0) at (0,0);
  \coordinate (q1) at (3,0);
  \coordinate (q2) at (3.8,2.4);
  \coordinate (p) at (0.8,2.4);
  \fill[black!12] (q0) -- (p) -- (q2) -- cycle;
  \draw (q0) -- (q1) -- (q2) -- (p) -- cycle;
  \draw (q0) -- (q2);
  \foreach \x in {q0,q1,q2,p} \node[vertex] at (\x) {};
  \node[below left] at (q0) {$q_0$};
  \node[below right] at (q1) {$q_0+v_1$};
  \node[above right] at (q2) {$q_0+v_1+v_2$};
  \node[above left] at (p) {$q_0+v_2$};
  \node[below] at (1.5,0) {$v_1$};
  \node[right] at (3.45,1.2) {$v_2$};
\end{tikzpicture}
\caption{Left: the Freudenthal subdivision of a triangle $\tau$ with
vertices $q_0<q_1<q_2$.  The three unshaded triangles, with vertices ordered
along their sequences $q_{00}q_{01}q_{02}$, $q_{01}q_{11}q_{12}$ and
$q_{02}q_{12}q_{22}$, have edge vectors $\frac12v_1,\frac12v_2$ and are
translates of $\frac12\tau$; the shaded triangle $q_{01}q_{02}q_{12}$ has
edge vectors $\frac12v_2,\frac12v_1$ and is a translate of $\frac12$ times the
conjugate of $\tau$.  Right: $\tau$ (unshaded) and its conjugate (shaded)
tile the parallelogram spanned by the edge vectors $v_1,v_2$ at $q_0$.}
\label{fig: freudenthal subdivision}
\end{figure}
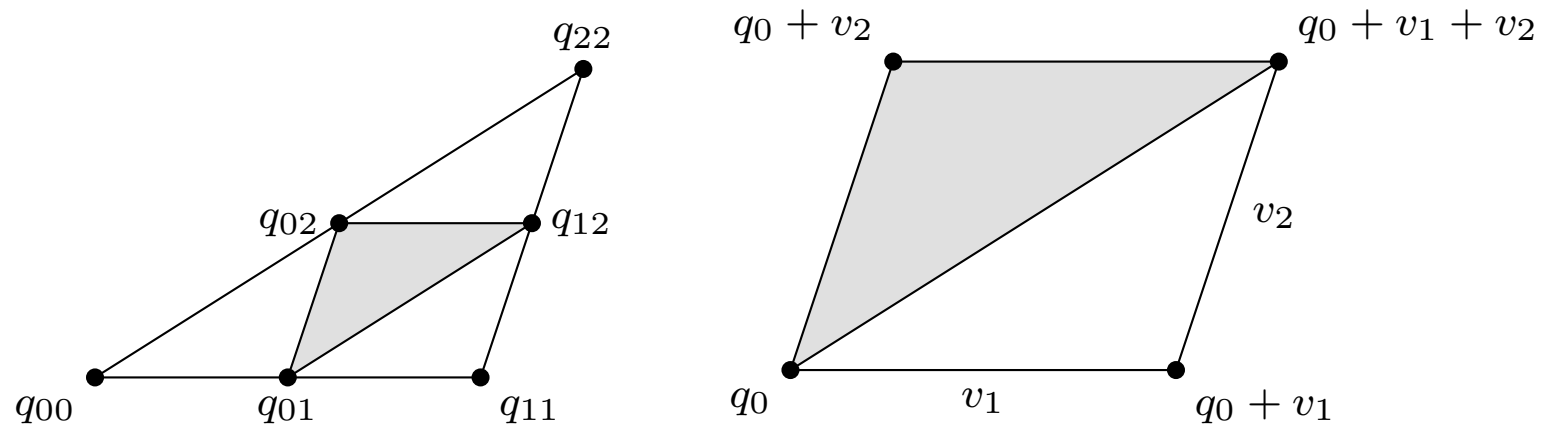

The vertices of each simplex of the subdivision are ordered along its sequence
\cite[Sections~3 and~4]{Freudenthal1942}.  With this order, each $d$-simplex of
the subdivision is a translate of $\frac12$ times a conjugate of $\tau$
\cite[Sections~2 and~3]{Freudenthal1942}; Figure~\ref{fig: freudenthal subdivision} shows the case $d=2$.

A simplicial complex whose simplices carry vertex orders that agree on common
faces is subdivided by subdividing each of its simplices with its order.  The Freudenthal
subdivision of a face, with the restricted order, consists of simplices of the
Freudenthal subdivision of the simplex
\cite[Face Lemma, with $k=2$]{EdelsbrunnerGrayson2000}, so these subdivisions
agree on common faces and together form a simplicial complex.  Its vertices are the
midpoints $q_{ij}$ of the edges and the vertices of the original complex, and
within each simplex of the original complex we order them
lexicographically in the pair $(q_i,q_j)$ \cite[Section~4]{Freudenthal1942}.  On each
simplex this is the order along its sequence, since both steps of a sequence
increase the pair lexicographically, and these orders agree on common faces.
Using them, the construction can be repeated.

For the rest of this subsection $K_0$ and $\chi$ are as in
Section~\ref{sec: initial triangulation}.  We write $K_j$ for the complex
obtained from $K_0$ after $j$ steps; these are the \emph{Freudenthal
subdivisions} of $K_0$.  Each $K_j$ is a simplicial complex with $|K_j|=|K_0|$,
and each of its simplices lies in a simplex of $K_0$.  If $\tau$ is a
$d$-simplex of $K_j$ with vertices $q_0<\dots<q_d$, let
$a_\tau\colon\Delta^d\to\tau$ be the affine map that sends $0$ to $q_0$ and the
$i$-th unit vector $e_i$ to $q_i$, and set $\imath_\tau:=\chi\circ a_\tau$.
We write $h_\tau$ for the $g$-diameter of $\chi(\tau)$, the image of
$\imath_\tau$.

The rest of this subsection shows that the Freudenthal subdivisions form a
uniformly regular family
(Theorem~\ref{thm: freudenthal subdivisions uniformly regular}).  The
conditions of Definition~\ref{def: uniformly regular} are verified one at a
time, and all of them rest on two facts.  The linear part of each
affine map $a_\tau$ is $2^{-j}$ times one of finitely many linear maps, which
depend only on $K_0$
(Lemma~\ref{lem: freudenthal linear parts}).  The map $\chi$ changes distances
within a simplex of $K_0$ by at most a bounded factor
(Lemma~\ref{lem: freudenthal distance comparison}).  From these follow, in
turn, the diameters of the simplices and with them the face ratio, the bounds
on the derivatives of $\imath_\tau$ and of its inverse, and the adjacency
bound.

The constants of Definition~\ref{def: uniformly regular} are defined for smooth
triangulations, so we first check that the parametrisations $\imath_\tau$ make
each $K_j$ one.

\begin{lemma}\label{lem: freudenthal smooth triangulation}
    For every $j\ge0$ the Freudenthal subdivision $K_j$, with the
    parametrisations $\imath_\tau=\chi\circ a_\tau$, is a smooth triangulation
    of $M$.
\end{lemma}

\begin{proof}
Each $\imath_\tau$ is a smooth embedding, as the composite of the affine
embedding $a_\tau$ into an $n$-simplex $\sigma$ of $K_0$ and the smooth
embedding $\chi|_\sigma$.  Vertex orders agree on common faces, so the restriction of $\imath_\tau$ to a face $\rho$ of
$\tau$ is $\imath_\rho$, and the smooth subsimplices of $\tau$ are its faces
with these parametrisations.  The conditions of
Definition~\ref{def: smooth triangulation} therefore follow from those of the
finite simplicial complex $K_j$, as $\chi$ is a homeomorphism and every simplex
of $K_j$ is a face of an $n$-simplex.
\end{proof}

To differentiate $\imath_\tau$ we write it as a fixed smooth map, one for each $n$-simplex of $K_0$, composed with the affine map $a_\tau$.  This is
possible because every simplex $\tau$ of $K_j$ lies in an $n$-simplex $\sigma$ of $K_0$. We use the smooth extension $F_\sigma$ of $\chi|_\sigma$ to an open neighbourhood of $\sigma$ in its affine hull, and write $V_\sigma$ for the space of vectors parallel to $\sigma$, with the Euclidean norm.  The differential of $F_\sigma$ is invertible at every point of $\sigma$, because $\chi|_\sigma$ has maximal rank.  If $\tau$ lies in $\sigma$, then
$\imath_\tau=F_\sigma\circ a_\tau$ on $\Delta^d$, because $F_\sigma=\chi$ on $\sigma$.

The derivatives of $\imath_\tau=F_\sigma\circ a_\tau$ are derivatives of
$F_\sigma$ evaluated on the linear part of $a_\tau$.  The next lemma shows that
these linear parts shrink at the rate $2^{-j}$ without degenerating.

\begin{lemma}\label{lem: freudenthal linear parts}
    There are constants $0<\beta_-\le\beta_+$, depending only on $K_0$, with the following property.  For every $j\ge0$ and
    every $d$-simplex $\tau$ of $K_j$, write $2^{-j}B_\tau$ for the linear part
    of $a_\tau$.  Then
    \begin{equation}\label{eq: subdivision linear parts}
        \beta_-|x|\le|B_\tau x|\le\beta_+|x|
        \qquad\text{for all } x\in\R^d .
    \end{equation}
\end{lemma}

\begin{proof}
Fix an $n$-simplex $\sigma$ of $K_0$, and let $\tau$ be a $d$-simplex of $K_j$
that lies in $\sigma$.  Then $\tau$ is a translate of $2^{-j}$ times a face of a
conjugate of $\sigma$, with the vertex order of that face.  Indeed, $\tau$ is a
face of an $n$-simplex $\tau'$ of $K_j$ that lies in $\sigma$, and with its
vertex order that $\tau'$ is a translate of $2^{-j}$ times a conjugate of $\sigma$.

Let $v_1,\dots,v_n\in V_\sigma$ be the edge vectors of $\sigma$, let $\varpi$ be
the permutation of that conjugate, and let $I=\{i_0<\dots<i_d\}$ index the
vertices of that face.  The linear part of $a_\tau$ sends $e_r$ to
$a_\tau(e_r)-a_\tau(0)=2^{-j}(v_{\varpi(i_0+1)}+\dots+v_{\varpi(i_r)})$, so
$B_\tau=B_\varpi L_I$, where $B_\varpi\colon\R^n\to V_\sigma$ sends $e_i$ to
$v_{\varpi(1)}+\dots+v_{\varpi(i)}$ and $L_I\colon\R^d\to\R^n$ sends $e_r$ to
$e_{i_r}-e_{i_0}$, with the convention $e_0:=0$.  The maps $B_\tau$ therefore range over a finite set of
injective linear maps that does not depend on $j$.  Indeed, there are finitely
many $n$-simplices $\sigma$ of $K_0$, permutations $\varpi$ and index sets $I$, and
every map $B_\varpi L_I$ that occurs is injective, because $a_\tau$ is.  So we may
take for $\beta_-$ and $\beta_+$ the smallest and the largest value of $|Bx|$
over unit vectors $x$ and over the maps $B$ of this set.
\end{proof}

Lemma~\ref{lem: freudenthal linear parts} measures the simplices of $K_j$ in
the Euclidean space that contains $K_0$.  The map $\chi$ carries these
measurements to $M$, and the next lemma says that within one $n$-simplex of
$K_0$ it is bi-Lipschitz: it changes distances by at most a bounded factor in
either direction.  Only the finitely many fixed maps $\chi|_\sigma$ are
involved, so the bounds come from compactness.

\begin{lemma}\label{lem: freudenthal distance comparison}
    There is a constant $c>0$, depending only on $(M,g)$ and $(K_0,\chi)$,
    such that for every $n$-simplex $\sigma$ of $K_0$
    \[
        c\,|p-q|\le d_g\big(\chi(p),\chi(q)\big)\le c^{-1}|p-q|
        \qquad\text{for all } p,q\in\sigma .
    \]
\end{lemma}

\begin{proof}
Fix an $n$-simplex $\sigma$ of $K_0$ and points $p,q\in\sigma$.  The upper bound
comes from the image of the segment from $p$ to $q$.  Indeed, let $\Lambda$
bound $|dF_\sigma(p')w|_g/|w|$ over $p'\in\sigma$ and nonzero $w\in V_\sigma$,
which is possible because $\sigma$ is compact.  The segment lies in $\sigma$,
because $\sigma$ is convex, and on $\sigma$ the map $\chi$ equals $F_\sigma$.  So
the image of the segment under $\chi$ is a smooth path from $\chi(p)$ to
$\chi(q)$ of $g$-length at most $\Lambda|p-q|$, and
$d_g(\chi(p),\chi(q))\le\Lambda|p-q|$.

For the lower bound we compare $d_g$ with a Euclidean distance in a chart.  By
\cite[Theorem~6.1]{Licht2023smoothed} $M$ has a finite atlas of charts
$\Psi_i\colon U_i\to\R^n$ and a constant $C$ such that
\begin{equation}\label{eq: chart distance}
    |\Psi_i(p')-\Psi_i(q')|\le C\,d_g(p',q')
    \qquad\text{for all }p',q'\in U_i
\end{equation}
\cite[(6.1a)]{Licht2023smoothed}.  By Lebesgue's number lemma there is
$\delta>0$ such that every subset of $M$ of diameter less than $\delta$ lies
in one of the $U_i$.  Licht's charts are restrictions of charts
defined on a neighbourhood of the compact closure of $U_i$, so
$\Psi_i\circ F_\sigma$ is defined and smooth on a neighbourhood of the compact
set $\sigma\cap F_\sigma^{-1}(\overline{U_i})$, on which its differential is
injective.  Let $m_0>0$ be the minimum of $|d(\Psi_i\circ F_\sigma)(p')w|$
over these sets, over unit vectors $w\in V_\sigma$ and over the finitely many
$i$, and let $\Lambda'>0$ bound the second derivative of $\Psi_i\circ F_\sigma$
on the same sets.

Near the diagonal the lower bound follows from Taylor's formula in one chart.
Precisely, put $r:=\min\{\delta/\Lambda,\ m_0/\Lambda'\}$ and suppose
$|p-q|<r$.  The image under $\chi$ of the segment from $q$ to $p$ has
$g$-diameter at most $\Lambda|p-q|<\delta$, so it lies in one chart $U_i$.
With $F:=\Psi_i\circ F_\sigma$, Taylor's formula along the segment gives
\begin{equation*}
    |F(p)-F(q)-dF(q)(p-q)|\le\tfrac12\Lambda'|p-q|^2 ,
\end{equation*}
and therefore $|F(p)-F(q)|\ge\frac12 m_0|p-q|$, as $|p-q|\le m_0/\Lambda'$.
By \eqref{eq: chart distance} for $\chi(p),\chi(q)\in U_i$,
\begin{equation*}
    d_g\big(\chi(p),\chi(q)\big)
    \ge C^{-1}\big|\Psi_i(\chi(p))-\Psi_i(\chi(q))\big|
    =C^{-1}|F(p)-F(q)|
    \ge\tfrac12C^{-1}m_0\,|p-q| .
\end{equation*}

Away from the diagonal the lower bound follows from compactness.  Indeed, on
the compact set of pairs $(p',q')\in\sigma\times\sigma$ with $|q'-p'|\ge r$
the quotient $d_g(\chi(p'),\chi(q'))/|q'-p'|$ is continuous, and it is
positive because $\chi$ is injective, so it has a positive minimum.

Taking for $c$ the smallest of $\frac12C^{-1}m_0$, this minimum and
$\Lambda^{-1}$, over the finitely many $n$-simplices $\sigma$ of $K_0$, gives
the claim.
\end{proof}

Together, the last two lemmas give the diameters of the simplices of $K_j$.
These set the scale against which Definition~\ref{def: uniformly regular}
measures the derivatives of the parametrisations.

\begin{lemma}\label{lem: freudenthal diameters}
    There are constants $0<c_1\le c_2$, depending only on $(M,g)$
    and $(K_0,\chi)$, such that
    $c_12^{-j}\le h_\tau\le c_22^{-j}$ for every $j\ge0$ and every simplex
    $\tau$ of $K_j$.
\end{lemma}

\begin{proof}
For a vertex, $h_\tau$ stands by convention for the smallest diameter of a
simplex of positive dimension containing it (Definition~\ref{def: uniformly regular}), so
it suffices to treat simplices of positive dimension.  Let $\tau$ be a simplex
of $K_j$ of dimension $d\ge1$ that lies in an $n$-simplex $\sigma$ of $K_0$,
and let $c$ be the constant of
Lemma~\ref{lem: freudenthal distance comparison}.  For $x,y\in\Delta^d$ the
points $a_\tau(x)$ and $a_\tau(y)$ lie in $\sigma$ and
$a_\tau(x)-a_\tau(y)=2^{-j}B_\tau(x-y)$, so the upper bound of
Lemma~\ref{lem: freudenthal distance comparison} gives
\begin{equation*}
    d_g\big(\imath_\tau(x),\imath_\tau(y)\big)
    \le c^{-1}|a_\tau(x)-a_\tau(y)|
    = c^{-1}2^{-j}|B_\tau(x-y)|
    \le c^{-1}\beta_+2^{-j}|x-y| ,
\end{equation*}
by \eqref{eq: subdivision linear parts} for the last step.  Hence
$h_\tau\le c_22^{-j}$ for some $c_2$.

The lower bound $c_12^{-j}\le h_\tau$, for some $c_1>0$, follows similarly from
the lower bounds of Lemma~\ref{lem: freudenthal distance comparison} and
\eqref{eq: subdivision linear parts}, applied to $\imath_\tau(0)$ and
$\imath_\tau(e_1)$.

We have $c_1\le c_2$ since both bounds hold for every edge of $K_j$.
\end{proof}

The first condition from Definition~\ref{def: uniformly regular} on the shape constant $C_\mu$, and the conditions on the
higher shape constants, ask for bounds on the derivatives of $\imath_\tau$ by
powers of $h_\tau$.  Since $a_\tau$ is affine, these derivatives are derivatives of
$F_\sigma$ evaluated on $2^{-j}B_\tau$, and
Lemma~\ref{lem: freudenthal diameters} turns powers of $2^{-j}$ into powers of
$h_\tau$.  Licht derives bounds of this form, and the matching bound for the inverse,
by the same chain rule \cite[Section~8.1]{Licht2023smoothed}. His bounds are in terms of the derivatives of $\chi$ and the shapes and diameters of the simplices of $K_j$, and the next two lemmas bound these independently of $j$.

\begin{lemma}\label{lem: freudenthal derivatives}
    For every $m\ge1$ there is a constant $C$, depending only on $(M,g)$,
    $(K_0,\chi)$ and $m$, such that for
    every $j\ge0$ and every $d$-simplex $\tau$ of $K_j$
    \[
        \sup_{\Delta^d}|\nabla^{m-1}d\imath_\tau|\le C\,h_\tau^m ,
    \]
    where $\nabla$ is as in Definition~\ref{def: uniformly regular}.
\end{lemma}

\begin{proof}
Let $\tau$ lie in an $n$-simplex $\sigma$ of $K_0$.  Write
$\nabla^{m-1}dF_\sigma$ for the covariant derivative of $dF_\sigma$, with the
flat connection on $V_\sigma$ and the Levi-Civita connection of $g$, measured as
a multilinear map from $V_\sigma$ to the tangent spaces of $(M,g)$.  It is
bounded on $\sigma$ uniformly in $\sigma$.  Indeed, it is continuous, $\sigma$
is compact, and there are finitely many $n$-simplices $\sigma$.  So there is a
constant $\Lambda_m$, depending only on $(M,g)$, $(K_0,\chi)$ and $m$, with
$|\nabla^{m-1}dF_\sigma(p)|\le\Lambda_m$ for every $n$-simplex $\sigma$ of $K_0$
and every $p\in\sigma$.

Since $a_\tau$ is affine, its higher derivatives vanish, and induction on $m$ gives
\begin{equation*}
    \nabla^{m-1}d\imath_\tau(x)
    =(\nabla^{m-1}dF_\sigma)\big(a_\tau(x)\big)
      \big(2^{-j}B_\tau\,\cdot\,,\dots,2^{-j}B_\tau\,\cdot\,\big)
    \qquad\text{for } x\in\Delta^d .
\end{equation*}

The claim follows by the bound on $\nabla^{m-1}dF_\sigma$,
\eqref{eq: subdivision linear parts} and the lower bound
$c_12^{-j}\le h_\tau$ of Lemma~\ref{lem: freudenthal diameters}.
\end{proof}

The second condition from Definition~\ref{def: uniformly regular} on the shape constant $C_\mu$ asks for a bound on the inverse of
$d\imath_\tau$.  We get this from the lower bound in
\eqref{eq: subdivision linear parts} and the invertibility of $dF_\sigma$ on the
compact simplex $\sigma$.

\begin{lemma}\label{lem: freudenthal inverse derivative}
    There is a constant $C$, depending only on $(M,g)$ and $(K_0,\chi)$, such that
    $\sup_{\chi(\tau)}|d\imath_\tau^{-1}|\le C\,h_\tau^{-1}$ for every $j\ge0$ and
    every simplex $\tau$ of $K_j$.
\end{lemma}

\begin{proof}
Let $\tau$ be a $d$-simplex of $K_j$ that lies in an $n$-simplex $\sigma$ of
$K_0$.  Let $\lambda$ be the minimum of $|dF_\sigma(p)w|_g$ over $p\in\sigma$, unit vectors
$w\in V_\sigma$ and the finitely many $n$-simplices $\sigma$ of $K_0$.  It is
positive, because $dF_\sigma$ is injective at every point of the compact set
$\sigma$.  For $x\in\Delta^d$ and $v\in\R^d$ the chain rule gives
$d\imath_\tau(x)v=dF_\sigma(a_\tau(x))\,2^{-j}B_\tau v$, and so, by the lower
bound in \eqref{eq: subdivision linear parts} and the upper bound
$h_\tau\le c_22^{-j}$ of Lemma~\ref{lem: freudenthal diameters},
\begin{equation*}
    |d\imath_\tau(x)v|_g
    \ge\lambda\,2^{-j}|B_\tau v|
    \ge\lambda\beta_-2^{-j}|v|
    \ge\frac{\lambda\beta_-}{c_2}\,h_\tau\,|v| .
\end{equation*}
This says that $|d\imath_\tau^{-1}|\le c_2/(\lambda\beta_-)\,h_\tau^{-1}$ at
every point of $\chi(\tau)$.
\end{proof}

The remaining constant from Definition~\ref{def: uniformly regular}, the adjacency bound, is combinatorial.  In a triangulation of a
full-dimensional subset of Euclidean space it is bounded in terms of the shapes
of the simplices, by comparing volumes, and
Lemma~\ref{lem: freudenthal linear parts} controls these shapes.

\begin{lemma}\label{lem: freudenthal adjacency}
    The adjacency bounds $C_\sharp$ of the $K_j$, $j\ge0$, are bounded by a
    constant depending only on $K_0$.
\end{lemma}

\begin{proof}
A simplex that meets $T$ meets it in a common face, so it contains a vertex of
$T$. It therefore suffices to bound the number of simplices of $K_j$ that
contain a given vertex $p$.  Such a simplex lies in one of the finitely many
$n$-simplices $\sigma$ of $K_0$ that contain $p$.  The simplices of $K_j$ that
lie in $\sigma$ triangulate $\sigma$, an $n$-simplex in its $n$-dimensional
affine hull.  The shape measure $\operatorname{diam}(\tau)^d/\operatorname{vol}_d(\tau)$
of each $d$-simplex $\tau$ among them, $d\ge1$, is bounded independently of
$j$, since $\tau=a_\tau(\Delta^d)$ and the linear part $2^{-j}B_\tau$ of $a_\tau$
satisfies \eqref{eq: subdivision linear parts} [so
$\operatorname{diam}(\tau)\le2^{-j}\beta_+\operatorname{diam}(\Delta^d)$ and
$\operatorname{vol}_d(\tau)\ge(2^{-j}\beta_-)^d\operatorname{vol}_d(\Delta^d)$].
By \cite[Lemma~II.4.1]{LichtThesis2017}, the number of simplices of
$K_j$ in $\sigma$ that contain $p$ is bounded by a constant depending only on
$n$, $\beta_-$ and $\beta_+$.  Summing over the $\sigma$ that contain $p$ gives
the bound.
\end{proof}

The lemmas together verify every condition of
Definition~\ref{def: uniformly regular} along the Freudenthal subdivisions.

\begin{theorem}\label{thm: freudenthal subdivisions uniformly regular}
    Let $(M,g)$ be a closed Riemannian manifold and $K_0$ a smooth
    triangulation of $M$, with $\chi$ as in
    Section~\ref{sec: initial triangulation}.
    Then the Freudenthal subdivisions $K_j$, $j\ge0$, with the parametrisations
    $\imath_\tau=\chi\circ a_\tau$, form a uniformly regular family of smooth
    triangulations of $M$.  Furthermore, there are constants $0<c_1\le c_2$ such
    that $c_12^{-j}\le h_\tau\le c_22^{-j}$ for every $j$ and every simplex $\tau$
    of $K_j$.
\end{theorem}

\begin{proof}
Each $K_j$ is a smooth triangulation of $M$, by
Lemma~\ref{lem: freudenthal smooth triangulation}.  The bounds
$c_12^{-j}\le h_\tau\le c_22^{-j}$ are provided by
Lemma~\ref{lem: freudenthal diameters}.  They give $\ell_{K_j}\le c_22^{-j}$, so
$\inf_j\ell_{K_j}=0$.  They also give $h_\tau\le c_22^{-j}\le(c_2/c_1)h_\rho$
whenever $\rho$ is a subsimplex of $\tau$, so $C_\triangle\le c_2/c_1$ for every
$K_j$.  The first condition on the shape constant $C_\mu$ follows from
Lemma~\ref{lem: freudenthal derivatives} with $m=1$, and the second from
Lemma~\ref{lem: freudenthal inverse derivative}.  The higher shape constants
$C_{\mu,m}$, $m\ge2$, are bounded by
Lemma~\ref{lem: freudenthal derivatives}, and the adjacency bound $C_\sharp$ by
Lemma~\ref{lem: freudenthal adjacency}.  None of these bounds depends on $j$,
so the $K_j$ form a uniformly regular family.
\end{proof}

\begin{remark}\label{rem: other subdivisions}
	The Freudenthal rule enters the proof of
	Theorem~\ref{thm: freudenthal subdivisions uniformly regular} only
	through Lemma~\ref{lem: freudenthal linear parts}.  The other lemmas use of
	$K_j$ only that it is a subdivision of $K_0$ with ordered vertices whose
	affine maps have linear parts as in \eqref{eq: subdivision linear parts}.
	We expect the proof to apply also to Whitney's standard subdivisions, which
	Dodziuk uses for his convergence theorems \cite[Section~2]{Dodziuk1976}.
\end{remark}

\subsection{Square-integrable forms}\label{sec:FEECforms}
In this section we recall the Hilbert complex of square-integrable forms (Section~\ref{sec: l2 forms}) and the regularity of the solutions of the two problems on it (Section~\ref{sec: regularity}). We then define the FEEC approximation to this complex (Section~\ref{sec: feec complex}), show that it is an approximation whose metric defect is of the order of the mesh size (Section~\ref{sec: feec approximation}), and deduce from an approximation estimate of Licht that $H^1\Lambda^\bullet(M)$ is a regularity scale for it (Section~\ref{sec: feec approximation estimates}).

\subsubsection{The Hilbert complex of square-integrable forms}\label{sec: l2 forms}
We follow the exposition of \cite[Section~4.1]{Holst2012}. On smooth $k$-forms,
the Hodge star $\star\colon\Omega^k(M)\to\Omega^{n-k}(M)$ and integration give
the $L^2$ inner product
$$
\langle \alpha,\beta\rangle=\int_M \alpha\wedge\star\overline{\beta}.
$$
The associated Hilbert space completion $L^2 \Lambda^k(M)$ of $\Omega^k(M)$ \cite{Gaffney1955,Schwarz1995} is the space of \emph{square-integrable $k$-forms}.

On $L^2 \Lambda^\bullet(M)$ we have the \emph{weak exterior derivative} $d$ \cite[Section~4.2]{arnold2010finite}, \cite[(2.3)]{Licht2023smoothed}: a form $\omega\in L^2\Lambda^k(M)$ lies in its domain of definition $V^k$ with $d\omega=\zeta\in L^2\Lambda^{k+1}(M)$ precisely when for every smooth $(n-k-1)$-form $\varphi$ we have
$$
\int_M\zeta\wedge\varphi=(-1)^{k+1}\int_M\omega\wedge d\varphi.
$$
A smooth form lies in $V^\bullet$, with weak exterior derivative its
classical one. The operator $d$ is a closed unbounded operator of degree one on the graded Hilbert space $L^2 \Lambda^\bullet(M)$ \cite[Section~7]{Licht2023smoothed}.

The weak derivative makes $L^2\Lambda^\bullet(M)$ into a Hilbert complex.  We call this
the \emph{Hilbert complex of square-integrable forms}. As $M$ is compact, its differentials
have closed range \cite[Theorem~14]{Kodaira1949}, \cite[Section~7]{Licht2023smoothed}, \cite[Lemma~2.4.3(b)]{Schwarz1995}, so it is a closed Hilbert
complex.

As $d$ is closed and densely defined, it admits a closed and densely defined adjoint $d^*$ with respect to
$\langle\cdot,\cdot\rangle$, with domain $V_\bullet\subset L^2\Lambda^\bullet(M)$.  On
smooth forms this is the \emph{codifferential} $\pm\star d\star$, which is adjoint to $d$ on a closed manifold by Green's formula.

The rates of Section~\ref{sec: convergence on manifolds} need Sobolev
regularity.  We write $H^1\Lambda^k(M)\subset V^k$ for the \emph{Sobolev space
of $k$-forms with square-integrable first covariant derivative}, taken with the Levi-Civita connection of $g$.  Its norm is $\|\omega\|_{H^1}=\|\omega\|+\|\nabla\omega\|$ \cite[(2.7)]{Licht2023smoothed}. For a form $u=\sum_ku^k$ of mixed degree we put
$\|u\|_{H^1}^2:=\sum_k\|u^k\|_{H^1}^2$ and
$H^1\Lambda^\bullet(M):=\bigoplus_kH^1\Lambda^k(M)$.
A solution to the Hodge--Dirac or Hodge--Laplace problem with square-integrable source has this regularity (Section~\ref{sec: regularity}).

Two standard facts about $d$, $d^*$ and $d+d^*$ are used below.  First, on smooth forms, $d$ and $d^*$ are the covariant derivative followed by a bundle map \cite[Chapter~II, Lemma~5.13]{LawsonMichelsohn1989}.  The
bundle map is exterior multiplication, respectively contraction, by the
covectors of an orthonormal frame, so its norm depends only on $n$.  The same holds for $\omega\in H^1\Lambda^k(M)$, because smooth forms are dense in $H^1\Lambda^k(M)$ \cite[Theorem~10.2.36(b)]{Nicolaescu2020} and $d$ and $d^*$ are closed.  Therefore $H^1\Lambda^k(M)$ also lies in the domain of $d^*$, and there
is a constant $C$, depending only on $n$, with
\begin{equation}\label{eq: d dstar bounded on H1}
    \|d\omega\|+\|d^*\omega\|\le C\|\omega\|_{H^1}
    \qquad\text{for }\omega\in H^1\Lambda^k(M).
\end{equation}
Second, $d+d^*$ is the Dirac operator of the Clifford bundle of $(M,g)$
\cite[Chapter~II, Theorem~5.12]{LawsonMichelsohn1989}.  So it is elliptic of order one \cite[Chapter~II, Lemma~5.1]{LawsonMichelsohn1989} and
formally self-adjoint \cite[Chapter~II, Proposition~5.3]{LawsonMichelsohn1989}.

\subsubsection{Regularity of the solutions}\label{sec: regularity}
By Lemma~\ref{lem: abstract mixed regularity}, the regularity of the solutions
follows from a Gaffney estimate for $H^1\Lambda^\bullet(M)$
(Definition~\ref{def: gaffney estimate}).
The next lemma proves the estimate from the elliptic regularity of
the first-order operator $d+d^*$.

\begin{lemma}\label{lem: feec dirac regularity}
	There is a constant $C_G$, depending only on $(M,g)$, such that every
	$\omega\in V^\bullet\cap V_\bullet$ lies in $H^1\Lambda^\bullet(M)$, with
	\[
	\|\omega\|_{H^1}\le C_G\big(\|\omega\|+\|d\omega\|+\|d^*\omega\|\big).
	\]
	In particular $\mathfrak H^\bullet\subset H^1\Lambda^\bullet(M)$, with
	$\|q\|_{H^1}\le C_G\|q\|$ for $q\in\mathfrak H^\bullet$.
\end{lemma}

\begin{proof}
	The operator $d+d^*$ is elliptic of order one and formally self-adjoint
	(Section~\ref{sec: l2 forms}), so elliptic regularity applies to weak solutions of $(d+d^*)u=v$, in the form of
	\cite[Theorem~10.3.11(a)]{Nicolaescu2020}. Let $\omega\in V^\bullet\cap V_\bullet$. By the definition of the adjoints, $\langle\omega,(d+d^*)\varphi\rangle=\langle d\omega+d^*\omega,\varphi\rangle$ for every smooth form $\varphi$, so $\omega$ is a weak solution of $(d+d^*)u=d\omega+d^*\omega$. The theorem gives that $\omega$ has a square-integrable weak covariant derivative \cite[Definition~10.2.33(b)]{Nicolaescu2020}, so $\omega\in H^1\Lambda^\bullet(M)$, with $\|\omega\|_{H^1}$ at most a constant depending only on $(M,g)$ times $\|d\omega+d^*\omega\|+\|\omega\|$.  This proves the first assertion. The second assertion is the case $d\omega=d^*\omega=0$, as $\mathfrak H^\bullet=\ker d\cap\ker d^*$.
\end{proof}

\begin{remark}\label{rem: gaffney}
	The estimate is Gaffney's inequality \cite[(2)]{Gaffney1955}, there stated
	for twice differentiable forms.  For $n\ge2$
	the lemma can also be assembled from
	\cite[Theorem~6.1]{IwaniecScottStroffolini1999}, which gives the weak
	first derivatives of a form in $V^\bullet\cap V_\bullet$, and
	\cite[Proposition~4.10]{Scott1995}, which gives the estimate for forms that
	have them, after passing to real and imaginary parts and identifying their
	chart-defined Sobolev space with $H^1\Lambda^\bullet(M)$.  The argument
	above covers $n=1$ as well.  In the terms of \cite[(7.19)]{Arnold2006}, the
	lemma says that a closed manifold is $1$-regular.  Second-order elliptic
	regularity for the Hodge--Laplace operator, the route of
	\cite[p.~315]{arnold2010finite} on a convex domain, is not needed.
\end{remark}

\subsubsection{The FEEC complex of forms}\label{sec: feec complex}
Choose a finite element de Rham complex
$\mathcal P\Lambda^0(\Delta^n)\xrightarrow{d}\cdots\xrightarrow{d}
\mathcal P\Lambda^n(\Delta^n)$ of Arnold, Falk and Winther, with spaces from the
trimmed and untrimmed families $\mathcal P_r^-\Lambda^k$, $\mathcal P_r\Lambda^k$
\cite[Sections~3.1, 3.3 and~3.5]{Arnold2006}. This choice, and with it the
polynomial degree, is fixed from now on: every finite element complex below is
built from these spaces. Following \cite[Section~4.1]{Licht2023smoothed}, we
transport them to a smooth triangulation along the parametrisations $\imath_T$.

\begin{definition}\label{def: finite element complex}
	Let $K$ be a smooth triangulation of $M$. The \emph{finite element complex}
	$V_K^\bullet$ on $K$ consists of the spaces $V_K^k$ of forms $\omega$ whose
	pullback $\imath_T^*\omega$ to $\Delta^n$ lies in
	$\mathcal P\Lambda^k(\Delta^n)$ for each $n$-simplex $T$, and whose traces
	on common subsimplices agree. If $K$ is realisable, $H_K^k$ denotes $V_K^k$ with the
	inner product
	\[
	\langle\alpha,\beta\rangle_K=\int_M\alpha\wedge\star_K\overline\beta
	\]
	of the piecewise flat metric $g_K$ of
	Definition~\ref{def: piecewise flat metric}, where $\star_K$ is the Hodge
	star operator of $g_K$, and $\|\cdot\|_K$ denotes the associated norm. We will use the symbol $V_K^k$ to denote the same space equipped with the associated graph inner
	product, as in Section~\ref{sec: approximating hilbert complexes}.
\end{definition}

These spaces form a finite-dimensional subcomplex of $V^\bullet$
\cite[Theorem~4.3 and (4.6)]{Licht2023smoothed}.  Every such space contains the lowest-order trimmed space $\mathcal P_1^-\Lambda^k(K)$, which is all that the estimates
below use. The lowest-order trimmed complex is the simplicial cochain complex of $K$ in another guise: it is spanned by the Whitney forms \cite[Section~4.3]{Arnold2006}, \cite[Definition~1.1]{Dodziuk1976}, and sending a cochain to its Whitney form is an isomorphism of complexes $C^\bullet(K)\cong\mathcal P_1^-\Lambda^\bullet(K)$ \cite[(1.4) and (1.7)]{Dodziuk1976}.

Denote by $i_K\colon H_K^\bullet\hookrightarrow L^2\Lambda^\bullet(M)$ the inclusion.
It is a chain map, but in general not an isometry, because $g_K$ differs from $g$.  For
the interpolator we take the smoothed projection
$\pi_K\colon L^2\Lambda^\bullet(M)\rar V_K^\bullet$ of
\cite[Theorem~6.11]{Licht2023smoothed} with $p=2$, which does not depend on $g_K$.

\subsubsection{The finite element complexes as an approximation}\label{sec: feec approximation}
This subsection shows that the finite element complexes form an approximation, in the sense of Definition~\ref{def: approximation of hilbert complex},
of the Hilbert complex of square-integrable forms. We bound the smoothed
projection and the metric defect, and apply
Lemma~\ref{lem: defects approximation}.
The interpolator does not involve the discrete metric, so we record its
properties first.

\begin{prop}\label{prop: feec smoothed projection}
	Let $K$ be a smooth triangulation of $M$ carrying a finite element complex $V_K^\bullet$.  The smoothed
	projection $\pi_K\colon L^2\Lambda^k(M)\to V_K^k$ is linear, satisfies
	$\pi_Kv=v$ for $v\in V_K^k$ and $d\pi_Ku=\pi_Kdu$ for
	$u\in V^k$, and
	$\|\pi_Ku\|\le C\|u\|$ for all $u\in L^2\Lambda^k(M)$.  The constant
	$C$ depends only on $(M,g)$, the polynomial degree and the constants of Definition~\ref{def: uniformly regular} for $K$.
\end{prop}

\begin{proof}
	Our projection $\pi_K$ is the smoothed projection from \cite[Theorem~6.11]{Licht2023smoothed}. This theorem supplies the desired properties; its commutation $d\pi_K=\pi_Kd$ is stated for forms in $L^{1,1}_{\mathrm{loc}}\Lambda^k(M)$, which contains $V^k$ because square-integrable forms are locally integrable. The constant $C$ in Licht's bound $\|\pi_Ku\|\le C\|u\|$ depends only on $M$, the polynomial degree
	and the regularity of $K$.  The regularity is measured there by the constants of Definition~\ref{def: uniformly regular} and by two further quantities, which these constants bound \cite[Section~3.2]{Licht2023smoothed}: the Lipschitz constant of $\imath_T^{-1}$ with respect to $d_g$, which is at most a multiple of $h_T^{-1}$ depending only on $n$, $C_\mu$ and $C_\triangle$ \cite[Theorem~3.4]{Licht2023smoothed}, and the Jacobian determinants of $\imath_T$ and $\imath_T^{-1}$, which are bounded in terms of $C_\mu$ by Hadamard's inequality.
	So Licht's constant depends only on the quantities named in the statement.
\end{proof}

Recall that the bounds from Section~\ref{sec: abstract problems} involve the
metric defect $\metricdefect[K]$ of \eqref{eq: graph metric defect}.  It is
at most the operator norm of $\id-i_K^*i_K$ on $H_K^\bullet$, and that norm is
controlled by the metric error, because the $L^2$ inner products of $g$ and
$g_K$ on forms differ by at most a multiple, depending only on $n$, of the
metric error.  Gawlik and McKee prove this for Regge
metrics, of which $g_K$ is one \cite{GawlikMcKee2024}; for the metric pulled back along a
diffeomorphism the corresponding bound goes back to Holst and Stern
\cite{Holst2012}.  For metrics interpolated by piecewise polynomials from a fixed
chart, Licht draws the same conclusion \cite{Licht2023smoothed}.

\begin{prop}\label{prop: feec metric defect}
	There are constants $C$ and $\epsilon_0>0$, depending only on $n$, such
	that for every realisable smooth triangulation $K$ of $M$ with
	$\epsilon_K\le\epsilon_0$, carrying any finite element complex,
	\begin{equation}\label{eq: feec metric defect}
		\opnorm{\id-i_K^*i_K}{H_K^k}{H_K^k}\le C\,\epsilon_K ,
	\end{equation}
	where $i_K^*\colon L^2\Lambda^k(M)\to H_K^k$ is the adjoint of $i_K$.
\end{prop}

\begin{proof}
	We apply \cite[Theorem~17]{GawlikMcKee2024}.  Its hypothesis and its bound are in terms of
	$\|g-g_K\|_{L^\infty(M,g_K)}$, the essential supremum over $x\in M$ of the
	norm of the symmetric form $(g-g_K)_x$ relative to $(g_K)_x$.  By
	\eqref{eq: feec metric error},
	$|(g-g_K)(a,a)|\le\epsilon_K\,g(a,a)\le\epsilon_K(1-\epsilon_K)^{-1}g_K(a,a)$,
	so this norm is at most $2\epsilon_K$ when $\epsilon_K\le\frac12$.  We choose
	$\epsilon_0\le\frac12$ so small that $2\epsilon_0$ lies below the threshold
	of that theorem, which then gives \eqref{eq: feec metric defect}.
\end{proof}

We collect the results of this subsection in one theorem: the finite element complexes form an approximation whose metric defect is of order $\ell_K$.

\begin{theorem}\label{thm: feec approximation}
    Let $(M,g)$ be a closed oriented Riemannian manifold, let $\mathcal K$ be a
    uniformly regular family of smooth triangulations of $M$, and let each
    $K\in\mathcal K$ carry the finite element complex
    $V_K^\bullet$ with the inner products of the piecewise flat metric $g_K$.
    There is $\ell_0>0$, depending only on $n$ and bounds for the shape
    constants $C_\mu$ and $C_{\mu,2}$ of the family, and at most the constant of
    Proposition~\ref{prop: piecewise flat metric error}, with the following
    property. If $\ell_K\le\ell_0$ for all $K\in\mathcal K$, then every
    $K\in\mathcal K$ is realisable and satisfies
    \begin{equation}\label{eq: feec metric defect rate}
        \opnorm{\id-i_K^*i_K}{H_K^\bullet}{H_K^\bullet}\le C\,\ell_K\le\tfrac14 ,
    \end{equation}
    with $C$ depending only on $n$ and those bounds; and
    $\{i_K,\pi_K\}_{K\in\mathcal K}$ is an approximation of
    $(L^2\Lambda^\bullet(M),V^\bullet,d)$, with $C_i$ and $C_\pi$
    depending only on $(M,g)$, the polynomial degree and bounds for the
    constants of Definition~\ref{def: uniformly regular}.
\end{theorem}

\begin{proof}
Proposition~\ref{prop: piecewise flat metric error}, applied to each
$K\in\mathcal K$, bounds $\epsilon_K$ by a multiple of $\ell_K$ whenever $\ell_K$ is smaller than the constant $\ell_0$, and by uniform regularity the multiple and $\ell_0$ can be chosen
independently of $K$.  Possibly shrinking $\ell_0$, we can further assume that this
multiple of $\ell_0$ is at most the constant $\epsilon_0$ of Proposition~\ref{prop: feec metric defect}, and let $\ell_K\le\ell_0$.  Then
$K$ is realisable with $\epsilon_K\le\epsilon_0$, so
Proposition~\ref{prop: feec metric defect} bounds the operator norm of
$\id-i_K^*i_K$ on $H_K^\bullet$ by $C\ell_K$.  Taking $\ell_0$ also so small
that $C\ell_0\le\frac14$ gives \eqref{eq: feec metric defect rate}.

We use Lemma~\ref{lem: defects approximation} to see that the family is an approximation.  Its
hypotheses hold: $i_K$ is the inclusion of a subcomplex, and $\pi_K$ is a
chain map with $\pi_Ki_K=\id$ by
Proposition~\ref{prop: feec smoothed projection}.  The metric defect satisfies
$\metricdefect[K]\le\opnorm{\id-i_K^*i_K}{H_K^\bullet}{H_K^\bullet}\le\frac14$, by
\eqref{eq: metric defect operator norm} and
\eqref{eq: feec metric defect rate}.  As $i_K$ is the inclusion,
Proposition~\ref{prop: feec smoothed projection} and $d\pi_K=\pi_Kd$ give
$\|i_K\pi_Ku\|_V\le C\|u\|_V$ for $u\in V^\bullet$, with the constant
$C$ of that proposition, which is bounded independently of $K$ because
$\mathcal K$ is uniformly regular.  The lemma gives $C_i$ depending only on $\nu_0=\frac14$, hence absolute, and
$C_\pi$ depending only on $(M,g)$, the polynomial degree and bounds for the
constants of Definition~\ref{def: uniformly regular}.
\end{proof}

From now on $\ell_0$ denotes the constant of Theorem~\ref{thm: feec approximation}.

\subsubsection{The regularity scale}\label{sec: feec approximation estimates}
The rates need one more property of the maps $i_K$ and $\pi_K$: that they
have regularity scale $H^1\Lambda^\bullet(M)$ with scale parameter $\ell_K$, in
the sense of Definition~\ref{def: regularity scale}.  That is,
$i_K\pi_K\omega$ must approximate an $H^1$ form
$\omega$ to order $\ell_K$.  Because $\pi_K$ is bounded on $L^2$ and fixes
$V_K^\bullet$, this reduces to a best approximation estimate: every $H^1$ form is approximated in $L^2$ by the finite
element space to the order of the mesh size.  This estimate follows from the piecewise Bramble--Hilbert lemma of
Licht \cite[Theorem~7.6]{Licht2023smoothed}.

\begin{theorem}\label{thm: feec l2 approximation}
    Let $(M,g)$ be a closed Riemannian manifold and $K$ a smooth triangulation
    of $M$ carrying a finite element complex $V_K^\bullet$.  There is a constant $C$, depending only on $n$, the
    polynomial degree and the constants $C_\mu$, $C_{\mu,2}$, $C_\sharp$ and
    $C_\triangle$ of $K$, such that for all $u\in H^1\Lambda^k(M)$
    \[
        \inf_{v\in V_K^k}\|u-v\|\le C\,\ell_K\,\|u\|_{H^1}.
    \]
\end{theorem}

\begin{proof}
For an $n$-simplex $T$ of $K$ and a form $\omega$ with one square-integrable
covariant derivative on $T$ we write
$\|\omega\|_{H^1(T)}:=\|\omega\|_{L^2(T)}+\|\nabla\omega\|_{L^2(T)}$, the norm
of \cite[(2.7)]{Licht2023smoothed} on $T$.
Let $u\in H^1\Lambda^k(M)$.  We apply
\cite[Theorem~7.6]{Licht2023smoothed}, with $m=1$ and $\ell=0$ in the notation
there.  The theorem gives a linear map
$\Pi\colon L^2\Lambda^k(M)\to V_K^k$ and a constant $C$ such that for every
$n$-simplex $T$ of $K$
\[
    \|u-\Pi u\|_{L^2(T)}
    \le C\sum_{T'\cap T\ne\emptyset}h_T\big(\|u\|_{H^1(T')}+\|du\|_{L^2(T')}\big) ,
\]
where the sum is over the $n$-simplices $T'$ of $K$ that meet $T$
\cite[(7.20)]{Licht2023smoothed}.  The constant $C$ depends only on $n$, the polynomial degree and $C_\mu$, $C_{\mu,2}$,
$C_\sharp$ and $C_\triangle$ (see the proof of
Proposition~\ref{prop: feec smoothed projection}), and so do the constants $C$ below, whose value
may change from one occurrence to the next.

The exterior derivative is
the covariant derivative followed by a bundle map whose norm depends only on
$n$ (Section~\ref{sec: l2 forms}), so
$\|du\|_{L^2(T')}\le C\|u\|_{H^1(T')}$.  With $h_T\le\ell_K$ the bound becomes
$\|u-\Pi u\|_{L^2(T)}\le C\,\ell_K\sum_{T'\cap T\ne\emptyset}\|u\|_{H^1(T')}$. At most $C_\sharp$ simplices meet $T$, so by the inequality between the arithmetic and quadratic means
\[
\|u-\Pi u\|_{L^2(T)}^2
\le C\,\ell_K^2\sum_{T'\cap T\ne\emptyset}\|u\|_{H^1(T')}^2 .
\]

It remains to sum these bounds over $T$. The $n$-simplices cover $M$ and overlap in null sets, and each $T'$ appears in the sum at most $C_\sharp$ times.  Summing over $T$ therefore gives
\[
    \|u-\Pi u\|^2\le C\,\ell_K^2\sum_{T'}\|u\|_{H^1(T')}^2
    \le C\,\ell_K^2\,\|u\|_{H^1}^2 ,
\]
because
$\sum_{T'}\|u\|_{H^1(T')}^2\le2\sum_{T'}\big(\|u\|_{L^2(T')}^2+\|\nabla u\|_{L^2(T')}^2\big)
\le2\|u\|_{H^1}^2$.  As $\Pi u\in V_K^k$, this proves the estimate.
\end{proof}

The regularity scale of Definition~\ref{def: regularity scale} asks for a
bound on the interpolation defect $\|u-i_K\pi_Ku\|$, not on the best
approximation.  The two differ by a constant factor, because $\pi_K$ is
bounded on $L^2$ and fixes $V_K^k$.

\begin{prop}\label{prop: feec interpolation defect}
	Let $K$ and $V_K^\bullet$ be as in
	Proposition~\ref{prop: feec smoothed projection}.
	There is a constant $C$, depending only on $(M,g)$, the polynomial degree
	and the constants of Definition~\ref{def: uniformly regular} for $K$, such
	that for all $u\in H^1\Lambda^k(M)$
	\begin{equation}\label{eq: feec l2 interpolation}
		\|u-i_K\pi_Ku\|\le C\,\ell_K\,\|u\|_{H^1}.
	\end{equation}
\end{prop}

\begin{proof}
	Let $u\in H^1\Lambda^k(M)$ and $v\in V_K^k$.
	Since $i_K\pi_Kv=v$ and $\|i_K\pi_Ku\|\le C\|u\|$ by
	Proposition~\ref{prop: feec smoothed projection},
	\[
	\|u-i_K\pi_Ku\|=\|(\id-i_K\pi_K)(u-v)\|\le(1+C)\|u-v\| .
	\]
	Taking the infimum over $v$ and applying Theorem~\ref{thm: feec l2 approximation}
	gives \eqref{eq: feec l2 interpolation}.
\end{proof}

\begin{cor}\label{cor: feec regularity scale}
	Let $\mathcal K$ be a uniformly regular family of smooth triangulations of
	$M$ with $\ell_K\le\ell_0$ for all $K\in\mathcal K$, each carrying a finite
	element complex. Then $(H^1\Lambda^\bullet(M),\|\cdot\|_{H^1})$ is a
	regularity scale for $\{i_K,\pi_K\}_{K\in\mathcal K}$ with scale parameter
	$\ell_K$, whatever inner products the finite element spaces carry, with
	$C_{\mathrm{app}}$ depending only on $(M,g)$, the polynomial degree and
	bounds for the constants of Definition~\ref{def: uniformly regular}.
\end{cor}

\begin{proof}
	$H^1\Lambda^k(M)\subseteq V^k$ by
	\eqref{eq: d dstar bounded on H1}, and \eqref{eq: feec l2 interpolation},
	with its constant bounded along the family, is
	\eqref{eq:abstract_approx_property}.
\end{proof}

\subsection{Convergence on Riemannian manifolds}\label{sec: convergence on manifolds}

We now combine the abstract error bounds of
Section~\ref{sec: abstract error bounds} with the properties of the finite element
complexes established above.  Theorem~\ref{thm: feec convergence} gives convergence with a rate for sources whose codifferential is square integrable, and Corollary~\ref{cor: feec L2 convergence} gives convergence for every source in the dual of the domain, in particular for every square-integrable source.  Throughout, $K$ is a smooth triangulation of the closed manifold
$M$ carrying the finite element complex and the maps $i_K$ and $\pi_K$ of
Section~\ref{sec: feec complex}; the subscript $h$ of
Section~\ref{sec: approximating hilbert complexes} becomes $K$.

\subsubsection{The discrete source}
As in Section~\ref{sec: abstract rates}, the discrete source is
$F_K=F\circ i_K$, which is defined for every $F\in(V^\bullet)^*$. Since $i_K$ is an
inclusion, the right-hand side of the discrete problem is then
$F(\beta_K)$ for test forms $\beta_K\in V_K^\bullet$: the
right-hand side of the continuous problem, restricted to finite element forms.
\subsubsection{The Hodge--Dirac and Hodge--Laplace problems}
We now apply the abstract theory to the finite element complexes with the
inner products of the piecewise flat metric $g_K$ and the discrete source
$F\circ i_K$.  We treat the two problems together.  For the Hodge--Dirac
problem with source $F$ and solutions $(\alpha,\kappa)$ and
$(\alpha_K,\kappa_K)$ write
\[
    \mathrm{err}_K:=\|\alpha-i_K\alpha_K\|_V+\|\kappa-i_K\kappa_K\| ,
\]
and for the mixed Hodge--Laplace problem in degree $k$ with solutions
$(\eta,\alpha,\kappa)$ and $(\eta_K,\alpha_K,\kappa_K)$
\[
    \mathrm{err}_K:=\|\eta-i_K\eta_K\|_V+\|\alpha-i_K\alpha_K\|_V+\|\kappa-i_K\kappa_K\| .
\]

The hypotheses of Theorem~\ref{thm: abstract rates} hold for finite element
exterior calculus on a uniformly regular family with $\ell_K\le\ell_0$.  The
complex $(L^2\Lambda^\bullet(M),V^\bullet,d)$ is closed
(Section~\ref{sec: l2 forms}). By
Theorem~\ref{thm: feec approximation}, $\{i_K,\pi_K\}_{K\in\mathcal K}$ is an
approximation of $(L^2\Lambda^\bullet(M),V^\bullet,d)$, and its metric defect is
of order $\ell_K$ by
\eqref{eq: feec metric defect rate} and \eqref{eq: metric defect operator norm}.
By Corollary~\ref{cor: feec regularity scale},
$(H^1\Lambda^\bullet(M),\|\cdot\|_{H^1})$ is a regularity scale for it with
scale parameter $\ell_K$. By Lemma~\ref{lem: feec dirac regularity}, this scale
satisfies a Gaffney estimate. This gives convergence at the rate of the mesh
size:

\begin{theorem}\label{thm: feec convergence}
    Let $(M,g)$ be a closed oriented Riemannian manifold, let $\mathcal K$
    be a uniformly regular family of smooth triangulations of $M$ with
    $\ell_K\le\ell_0$ for all $K\in\mathcal K$, where $\ell_0$ is as
    in Theorem~\ref{thm: feec approximation}, and let each $K\in\mathcal K$
    carry the finite element complex $V_K^\bullet$ with the inner products of
    the piecewise flat metric $g_K$.  Consider either the Hodge--Dirac
    problem \eqref{eq: Hodge Dirac problem} with source $\gamma\in V_\bullet$,
    or the mixed Hodge--Laplace problem \eqref{eq: Hodge Laplace problem} in
    degree $k$ with source $\gamma\in V_k$, and for each $K$ the discrete
    source $i_K^*\gamma$.  Then the error $\mathrm{err}_K$ between the solution of the problem
    on $(L^2\Lambda^\bullet(M),V^\bullet,d)$ with source $\gamma$ and the
    solution of the problem on $(H_K^\bullet,V_K^\bullet,d_K)$ with source
    $i_K^*\gamma$ satisfies
    \begin{equation}\label{eq: feec rate}
        \mathrm{err}_K\le C\,\ell_K\big(\|\gamma\|+\|d^*\gamma\|\big).
    \end{equation}
    The constant $C$ depends only on $(M,g)$, the polynomial degree and
    bounds for the constants of Definition~\ref{def: uniformly regular}, and
    not on $K$ or on $\gamma$.
\end{theorem}

\begin{proof}
We apply Theorem~\ref{thm: abstract rates}. Since $\metricdefect[K]\le C\ell_K$ by \eqref{eq: feec metric defect rate}, its error bounds \eqref{eq: abstract HD rate} and \eqref{eq: abstract HL rate} give \eqref{eq: feec rate}. Its constant
depends only on $c_H$, $C_i$, $C_\pi$, $C_{\mathrm{app}}$, $C_G$ and the
constant of \eqref{eq: feec metric defect rate}.  The Poincar\'e constant $c_H$
and the constant $C_G$ of Lemma~\ref{lem: feec dirac regularity} depend only on
$(M,g)$, and the others only on $(M,g)$, the polynomial degree and bounds for
the constants of Definition~\ref{def: uniformly regular}, by
Theorem~\ref{thm: feec approximation} and
Corollary~\ref{cor: feec regularity scale}.
\end{proof}

\begin{remark}\label{rem: feec regularity sharp}
	The hypothesis $\gamma\in V_k$ places the exact part $d\eta$ of the source
	in $H^1\Lambda^k(M)$.  Arnold, Falk and Winther make a hypothesis of the
	same kind for the error in $d\sigma$, the analogue of $d\eta$: on an
	$s$-regular domain in $\R^n$ \cite[(7.19)]{Arnold2006}, their rate for
	$d(\sigma-\sigma_h)$ assumes that the exact part $f_d$ of the source lies
	in $H^t\Lambda^k$ \cite[Theorem~7.10]{Arnold2006}.
\end{remark}

\subsubsection{Sources without a rate}
For a source that is merely square integrable, or only a bounded functional on the domain,
Theorem~\ref{thm: feec convergence} does not apply, but
Corollary~\ref{cor: abstract L2 convergence} still gives convergence, as the
metric defect tends to zero with the mesh size.
The finite element literature deduces convergence for square-integrable sources
from the density of the discrete spaces in $V^\bullet$ in the graph norm
\cite[Section~3.5, (8)]{arnold2010finite},
\cite[Corollary~3.11 and the sentence after it]{Holst2012},
\cite[(4)]{LeopardiStern2016}; Corollary~\ref{cor: abstract L2 convergence}
approximates the source instead.

\begin{cor}\label{cor: feec L2 convergence}
    Let $\mathcal K$ and the inner products be as in
    Theorem~\ref{thm: feec convergence}.  Consider either the Hodge--Dirac
    problem \eqref{eq: Hodge Dirac problem} with source
    $F\in(V^\bullet)^*$, or the mixed Hodge--Laplace problem
    \eqref{eq: Hodge Laplace problem} in degree $k$ with source
    $F\in(V^k)^*$, and for each $K$ the discrete source
    $F\circ i_K$.  Then the error $\mathrm{err}_K$ between the solution of the
    problem on $(L^2\Lambda^\bullet(M),V^\bullet,d)$ with source $F$
    and the solution of the problem on $(H_K^\bullet,V_K^\bullet,d_K)$
    with source $F\circ i_K$ tends to zero as $\ell_K\to0$: for every
    $\varepsilon>0$ there is $\ell>0$ such that $\mathrm{err}_K\le\varepsilon$ for
    all $K\in\mathcal K$ with $\ell_K\le\ell$.
\end{cor}

\begin{proof}
We apply Corollary~\ref{cor: abstract L2 convergence}. Let
$\varepsilon>0$, and let $\delta>0$ be as in that corollary.  Since
$\metricdefect[K]\le C\ell_K$, every $K\in\mathcal K$ with
$\ell_K\le\delta/(1+C)$ satisfies $\ell_K+\metricdefect[K]\le\delta$, and so
$\mathrm{err}_K\le\varepsilon$.
\end{proof}

\addtocontents{toc}{\protect\setcounter{tocdepth}{0}}
\section*{Acknowledgements}
Thanks go to Martin Licht for helpful comments on a draft of this article. The author is grateful to Owen Gwilliam for his support during the preparation of the manuscript. This research was supported by the European Union's Horizon Europe
research and innovation programme under the Marie Sk\l{}odowska-Curie
grant agreement No~101203556 (DisQ FA), and by Aalto University.

\section*{Tool and computational resource disclosure}
Large language models (predominantly Claude) were used throughout the drafting of this paper. They were used for literature search and for drafting and revising prose and proofs. The author has revised and edited all text, and has ensured that every proof is correct, readable and properly attributed.

The author takes full responsibility for the correctness of the arguments and results, and for the completeness and accuracy of the citations.  No automated system is an author of this paper; the core intellectual contributions are the author's.

	\addtocontents{toc}{\protect\setcounter{tocdepth}{2}}
	\bibliographystyle{alpha}
	\bibliography{lib}

\end{document}